\documentclass[a4paper,12pt]{article}
\usepackage{latexsym}
\usepackage[francais]{babel}
\usepackage{amsfonts,amsthm,amsmath}

\let\cal\mathcal

\newcommand {\C}              {{\Bbb C}}

\newcommand {\R}              {{\Bbb R}}

\newcommand {\N}              {{\Bbb N}}

\newcommand {\Z}              {{\Bbb Z}}

\newcommand {\T}              {{\Bbb T}}

\newcommand {\al}             {\alpha}

\newcommand {\be}             {\beta}

\newcommand {\un}             {\underline}

\newcommand {\ov}             {\overline}

\newcommand {\f}             {\hat}

\newcommand {\lov}             {\overleftarrow}
\newcommand {\rov}             {\overrightarrow}

\newcommand {\bc} {\begin{cases}}

\newcommand {\ec} {\end{cases}}

\newcommand {\bcd} {\begin{cases}\displaystyle}

\newcommand {\ecd} {\end{cases}\displaystyle}

\newcommand {\bml} {\begin{multline*}}

\newcommand {\bs} {\begin{small}}

\newcommand {\es} {\end{small}}

\newcommand {\di} {\displaystyle}
\newcommand {\beq}{\begin{equation*}}
\newcommand {\eeq}{\end{equation*}}

\newcommand {\eml} {\end{multline*}}

\theoremstyle{remark}

\newtheorem{remarque}{Remarque}[section]

\theoremstyle{plain}

\newtheorem{theoreme}[remarque]{Th\'eor\`eme}

\newtheorem{lemme}[remarque]{Lemme}

\newtheorem{proposition}[remarque]{Proposition}

\newtheorem{corolaire}[remarque]{Corolaire}

\newtheorem{definition}[remarque]{Definition}

\newcommand {\eqn}{\begin{equation}}

\normalsize

\author{Jean-Baptiste Yvernault
 \\ Universit\'e Pierre et Marie Curie  
\\ UFR 920  BP 172 
\\ 4, place Jussieu 75252 Paris cedex 05  
\\  e-mail :yvernaul@math.jussieu.fr }

\begin{document}

\bibliographystyle{plain}

\centerline{\bf  THEOREME DE CAUCHY GLOBAL POUR }
\centerline{\bf LES EQUATIONS D'EVOLUTIONS NONLINEAIRES}
\baselineskip 15pt
\vspace{5mm}
\centerline{ J.-B YVERNAULT}
\baselineskip 13pt
\centerline{ \it Universit\'e Pierre et Marie Curie }
\centerline{ \it UFR 920  BP 172}
\centerline{ \it 4, place Jussieu 75252 Paris cedex 05 }
\centerline{ \it e-mail :yvernaul@math.jussieu.fr }


\begin{abstract}
Nous \'etudions les  \'equations d'\'evolutions semilin\'eaires \`a conditions initiales
p\'eriodiques. On d\'emontre alors un th\'eor\`eme d'existence globale des
perturbations Hamiltoniennes de KdV sans aucune restriction sur la taille des
conditions initiales. Pour cel\`a on d\'eveloppe une m\'ethode
bas\'ee sur la th\'eorie des graphes. On en d\'eduit aussi,  dans le cas g\'en\'eral, une condition de
non existence de solutions.
\end{abstract}

\baselineskip 15pt
\hfuzz=11pt

\section{Introduction}
\subsection{Le probl\`eme}
Soit $H^s(\T)$, l'espace de hilbert des fonction $2\pi$
p\'eriodique $s$ fois d\'erivable ( \`a valeur r\'eelle), muni de la norme $\parallel u\parallel_s=\sqrt{\sum_{n\in
\N}n^{2s}\f{u}_n\f{u}_{-n}}$.  Soit $W_s(\T)$ espace de Banach des fonction $2\pi$
p\'eriodique $s$ d\'erivable de norme $\mid u\mid_s=\sum_{n\in
\N}n^s\mid \f{u}_n\mid$.  On s'int{\'e}resse aux solutions sur $H^s(\T)$ (
respectivement sur $W^s(\T)$ de l'{\'e}quation:
\begin{equation}{\label{1}}
\begin{cases}\displaystyle
\frac{\partial u}{\partial t} =\sum_{i=0}^{2k+1}\lambda_i
\frac{\partial^{i} u}{\partial x^i}+P(u,\partial_{x}u,...\partial_{x^k}u)\\
u(x,0)=u_0(x)\in H^s(\T) \hbox{ respectivement } W^s(\T).\end{cases}
\end{equation}
 On montre dans quels cas l'\'equation est bien pos\'ee.
Il s'agit d'un analogue des r\'esultats d\'ej\`a connus sur les espace de Sobolev
($H^s(\R)$)( voir  \cite{KPV}), \'etendu au cas des conditions initiales p\'eriodique. J
Bourgain a d\'emontr\'e dans \cite{BOU5}, l'existence de solutions locale puis
globales ( dans le cas Hamiltonien) pour des conditions initiales suffisamment
petites. Nous \'etendons ses r\'esultats {\it sans restriction } sur la ``taille''  des
conditions initiales. Puis nous proposons une g\'en\'eralisation de ses m\'ethodes
\`a des \'equations non hamiltoniennes et d\'emontrons un r\'esultat de non-existence de solutions. Notre m\'ethode repose sur une construction directe de
solutions en utilisant le formalisme des arbres introduit par Gallavoti \cite{GAL}.

Pour comprendre les m{\'e}thodes utilis{\'e}es, consid\'erons le  cas particulier:
\begin{equation}
\begin{cases}
\partial_tu=\partial_x^3u+u^2\partial_xu\\
u(0,x)=v\in \cal {H}^k(\T),
\end{cases}
\end{equation}
o{\`u} $\cal {H}^k(\T)$ est le sous espace de $H^k(\T)$ telles  que pour tout
$v$ de $\cal {H}^k(\T)$, on ait $\int_\T v=0=\f{v}_0$. On remarque d'abord que
$\f{u}_0(t)=0$ et  $\int_\T u^2(t)=\int_\T v^2$ qui sera not\'e  $\Omega_0$.
Si on {\'e}crit cette {\'e}quation {\`a} l'aide de ses coefficients de Fourier, on obtient:
\begin{equation}\label{Eqex}
\begin{cases}\displaystyle
\frac{d\f{u}_n}{dt}=((in)^3+in\Omega_0)\f{u}_n+\sum_{\substack{(p,q)\in \Z^2\\
(p+q)(n-q)(n-p)\neq 0}}(in)\f{u}_{p}\f{u}_{q}\f{u}_{n-p-q},\\
\f{u}_n(0)=\f{v}_n.
\end{cases}
\end{equation}
L'\'equation lin\'earis\'ee de (\ref{Eqex}) est :
\begin{equation}\label{Eqexlin}
\begin{cases}
\frac{d\f{w}_n}{dt}=((in)^3+in\Omega_0)\f{w}_n,\\
\f{w}_n(0)=\f{v}_n.
\end{cases}
\end{equation}
On cherche une application $\phi$ telle que si $w(t)$ est solution de l'\'equation
lin\'earis\'ee (\ref{Eqexlin}), alors $u(t)=\phi(t,w(t))$ est solution de l'\'equation (\ref{Eqex}). Si on note
$\f{u}_n=\f{\phi}_n(t,w(t))$, on cherche les $\f{\phi}_n$ sous la forme suivante
:
\[\f{\phi}_n(t,w)=\f{w}_n+\sum_{\al=2}^\infty \sum_{\un{n}\in\ov{\Z}^\al(n)}\f{\phi}_{\un{n}}(t)\prod_{i=1}^\al\f{w}_{n_i},\]
o{\`u} les $\f{\phi}_{\un{n}}$ sont des fonctions \`a valeur dans $\C$ et
$\ov{\Z}^\al(n)$ est l'ensemble des multi-entiers, modulo les permutations
d'indice et  dont la somme vaut $n$.
La construction formelle de cette application $\phi$ est donn{\'e}e par  une
r{\'e}curence sur l'ordre $\al$  des
$\f{\phi}_{\un{n}}$, avec $\un{n}\in \ov{\Z}^\al(n)$, obtenus en rempla\c cant $\f{\phi}_n(w)$ dans
l'{\'e}quation (\ref{Eqex}). Nous obtenons alors :
\begin{small}
\begin{multline*}
\partial_t\f{\phi}_{\un{n}}(t)+\left(\sum_{j=1}^\al(in_j)^3-(in)^3\right)\f{\phi}_{\un{n}}(t)=\\
\sum_{\substack{\un{m_i}\in
\ov{\Z}^{\al_i}(m_i)\\m_1+m_2+m_3=n\\\al_1+\al_2+\al_3=\al\\
(m_1+m_2)(m_1+m_3)(m_2+m_3)\neq
0}}(in)\f{\phi}_{\un{m_1}}(t)\f{\phi}_{\un{m_2}}(t)\f{\phi}_{\un{m_3}}(t),\end{multline*}
\end{small}
avec pour condition initiale
$\f{\phi}_{\un{n}}(0)=0$.

C'est  une \'equation differentielle du premier ordre avec pour condition
initiale $\f{\phi}_{\un{n}}(0)=0$, qui est toujours r{\'e}solvable, que
$\left(\sum_{j=1}^\al(in_j)^3-(in)^3\right)$ soit nul ou non. On obtient :
\bs
\bml
\f{\phi}_{\un{n}}(t)=e^{-t(\sum_{j=1}^\al(in_j)^3-(in)^3)}in\times\\\sum_{\substack{\un{m_i}\in
\ov{\Z}^{\al_i}(m_i)\\m_1+m_2+m_3=n\\\al_1+\al_2+\al_3=\al\\
(m_1+m_2)(m_1+m_3)(m_2+m_3)\neq
0}}( 
\int_0^te^{s(\sum_{j=1}^\al(in_j)^3-(in)^3)}\f{\phi}_{\un{m_1}}(s)\f{\phi}_{\un{m_2}}(s)\f{\phi}_{\un{m_3}}(s)ds).
\end{multline*}
\es

Pour  montrer l'existence de $\phi$ sur ${\cal H}^k(\T)$ ( $k>2$), on  montre que pour $v\in {\cal H}^k(\T)$ fix{\'e}, il existe $t_0$ (ne d{\'e}pendant que de $v=u(0)$) tel que pour tout $t\leq t_0$ on ait:
\[\sum_{n\in \N^*}\left(\sum_{\al=1}^\infty\sum_{\un{n}\in\ov{\Z}^\al(n)}\left\vert n^k\f{\phi}_{\un{n}}(t)\prod_{j=1}^\al\f{v}_{n_j}\right\vert\right)^2<\infty.\]
Le principal obstacle vient du fait que la perturbation est non born{\'e}e, ce
 qui signifie :
 \[\sup_{u:\parallel u\parallel_k\leq 1}\frac{\parallel
 u^2\partial_xu\parallel_k}{\parallel u \parallel_k}=\infty,\] 
 symbolis{\'e} par la multiplication par $in$ du terme non lin{\'e}aire.

On propose maintenant, une m\'ethode pour lever cet obstacle:
 
\subsection{Cont{\^o}le des perturbations non born{\'e}es}
Les coefficients $\f{\phi}_{\un{n}}$ sont formellement donn{\'e}s par des
produits d'op{\'e}rateurs du type $(\partial_t+D)^{-1}$, o{\`u} les $D$ sont des
nombres complexes imaginaires purs ( ou nuls). Nous nous proposons  de
montrer que la norme $\parallel \phi \parallel_k$, est major{\'e}e pour un temps $t_0$ d{\'e}pendant
de $\parallel v\parallel_k$ par une constante d{\'e}pendant uniquement de $\parallel
v\parallel_k$. 

Au premier ordre de r{\'e}curence, les termes $\f{\phi}_{\un{n}}$ sont donn{\'e}s
par une somme d'{\'e}l{\'e}ments du type :
\[ine^{-t((in_1)^3+(in_2)^3+(in_3)^3-(in)^3)}\int_0^te^{s((in_1)^3+(in_2)^3+(in_3)^3-(in)^3)}ds.\]
Au deuxi{\`e}me ordre, en posant
$D_1=(in_1)^3+(in_2)^3+(in_3)^3+(in_4)^3+(in_5)^3-(in)^3)$ et
$D_2=((in_3)^3+(in_4)^3+(in_5)^3-(i(n_3+n_4+n_5))^3)$, les termes
$\f{\phi}_{\un{n}}$ sont de la forme :
\[in(i(n_3+n_4+n_5))e^{-tD_1}\int_0^te^{s_1D_1}e^{-s_1D_2}\int_0^{s_1}e^{s_2D_2}ds_1ds_2.\]
Nous allons d'abord expliquer la m{\'e}thode sur les termes d'ordre $1$. 

Supposons que $\mid n_1\mid \geq \mid n_2\mid \geq \mid n_3\mid$. On remarque alors que $3\mid n_1\mid\geq \mid n\mid$ et que
\[(n_1)^3+(n_2)^3+(n_3)^3-(n)^3)=3(n_1+n_2)(n_1+n_3)(n_2+n_3).\]
Alors si $2\mid n_2\mid\leq \mid n_1\mid$ et $2\mid n_3\mid\leq \mid n_1\mid$, on a  \[\mid (n_1+n_2)(n_1+n_3)(n_2+n_3)\mid\geq \frac{\mid n_1\mid^2}{4}\geq\frac{\mid n\mid^2}{36}.\]
Donc \[\left\vert ine^{-t((in_1)^3+(in_2)^3+(in_3)^3-(in)^3)}\int_0^te^{s((in_1)^3+(in_2)^3+(in_3)^3-(in)^3)}ds\right\vert\leq K,\] ( $K$ est une constante ind{\'e}pendante de $\un{n}$).

Par contre, si $3\mid n_1\mid\leq \mid n\mid$, on obtient :
\bml
n^k\left\vert
ine^{-t((in_1)^3+(in_2)^3+(in_3)^3-(in)^3)}\int_0^te^{s((in_1)^3+(in_2)^3+(in_3)^3-(in)^3)}ds\right\vert\\\leq
Kt\mid n_1^kn_2\mid.\end{multline*}
Ceci montre alors qu'au premier ordre, 
\[\left\vert n^k\f{\phi}_{(n_1,n_2,n_3)}(t)\f{v}_{n_1}\f{v}_{n_2}\f{v}_{n_3} \right\vert\leq K(1+t)\left\vert n_1^k\f{v}_{n_1}n_2\f{v}_{n_2}\f{v}_{n_3} \right\vert.\]
Donc si on somme tous les termes possibles du premier ordre nous avons bien une
majoration de cette somme par une constante d{\'e}pendant uniquement ( en ce qui
concerne la condition initiale) de $\parallel v\parallel_k$.

A l'ordre $2$, si on suppose $\mid n_1\mid \geq \mid n_2\mid \geq \mid n_3\mid
\geq \mid n_4\mid \geq \mid n_5\mid $, on a de la m\^eme fa\c con que $\mid
D_1\mid $ ( respectivement  $\mid D_2\mid $) est d'ordre de grandeur $n$ (
respectivement $(n_3+n_4+n_5)$), ou bien $\mid n_1\mid$ et $\mid n_2\mid$ sont
d'ordre de grandeur $\mid n\mid$ ( respectivement $\mid n_3$ et $\mid n_4\mid$
sont d'ordre de grandeur $\mid n_3+n_4+n_5\mid$). On le d{\'e}montre   page
\pageref{lemccv0}, lemme \ref{lemccv0}. Il nous suffit alors de montrer que pour
$t\leq 1$, on a 
\begin{multline*}
\left\vert
e^{-tD_1}\int_0^te^{s_1D_1}e^{-s_1D_2}\int_0^{s_1}e^{s_2D_2}ds_1ds_2\right\vert\leq
\\
\inf(t,\frac{1}{\mid D_1 \mid})\inf(t,\frac{1}{\mid D_2
\mid})+\inf(t,\frac{1}{\mid D_1-D_2 \mid})\inf(t,\frac{1}{\mid D_2 \mid}).
\end{multline*}
Cette propri{\'e}t{\'e} est d{\'e}montr\'ee en annexe, page \pageref{lemaqD} lemme \ref{lemaqD}, en remarquant dans l'exemple d'ordre $2$ que,
 si $\mid D_1\mid\geq \frac{1}{t}$ et $\mid D_2\mid\geq \frac{1}{t}$, le
 r{\'e}sultat s'obtient en calculant les int{\'e}grales
\begin{multline*}
\left\vert
e^{-tD_1}\int_0^te^{s_1D_1}e^{-s_1D_2}\int_0^{s_1}e^{s_2D_2}ds_1ds_2\right\vert=
\\
\left\vert\frac{1-e^{-tD_1}}{D_1D_2}-\frac{e^{-tD_1}\int_0^te^{(D_1-D_2)s}}{D_2}ds\right\vert.
\end{multline*}

 Si $\mid D_1\mid< \frac{1}{t}$ il suffit de remarquer que :
\[\left\vert e^{-tD_1}\int_0^te^{s_1D_1}e^{-s_1D_2}\int_0^{s_1}e^{s_2D_2}ds_1ds_2\right\vert\leq t\sup_{s_1\in [0,t]}\left\vert\int_0^{s_1}e^{s_2D_2}ds_2\right\vert.\]
Si $\mid D_1\mid\geq \frac{1}{t}$ et $\mid D_2\mid< \frac{1}{t}$, il suffit
d'int{\'e}grer par partie en posant
$f(s_1)=e^{-s_1D_2}\int_0^{s_1}e^{s_2D_2}ds_2$ :
\[e^{-tD_1}\int_0^te^{s_1D_1}f(s_1)ds_1=\frac{f(t)}{D_1}-e^{-tD_1}\int_0^te^{s_1D_1}f'(s_1)ds_1.\]
En remarquant que $\mid f'(t)\mid\leq \mid D_2\mid \mid f(t)\mid+1$ et que $\mid
f(t)\mid \leq t$, ceci implique alors que :
\[\left\vert e^{-tD_1}\int_0^te^{s_1D_1}e^{-s_1D_2}\int_0^{s_1}e^{s_2D_2}ds_1ds_2\right\vert\leq \frac{t+2t^2}{\mid D_1\mid}.\]

On d{\'e}montre ainsi une in{\'e}galit{\'e} du type:
\[\left\vert n^k\f{\phi}_{\un{n}}(t)\prod_{j=1}^\al\f{v}_{n_j}\right\vert\leq \frac{1}{\mid n\mid^k}\left\vert (\epsilon(t))^\al\prod_{j=1}^\al n_j^k\f{v}_{n_j}\right\vert,\]
o\`u $\epsilon(t)<M \ll 1$ quand $t\leq t_0$.

Pour d{\'e}montrer cette propri{\'e}t{\'e} pour les termes d'ordre sup{\'e}rieur,
nous avons besoin du formalisme des arbres, qui s'av\`ere bien adapt\'e dans ce
contexte.

Nous donnons maintenant les th\'eor\`emes principaux qui font l'objet de cet
article. Ensuite nous exposons le formalisme des arbres que nous allons utiliser
ainsi que les propri\'et\'es classiques qui lui sont li\'es. Dans la quatri\`eme
partie, nous appliquons ce formalisme pour construire les solutions formelles de
notre probl\`emes. La cinqui\`eme partie consiste \`a montrer l'existence de ses
solutions formelles en d\'emontrant la convergence. La d\'emonstration s'appuie
sur des lemmes tr\`es techniques, donn\'es en annexe et sur les d\'efinitions du
formalisme des arbres donn\'e dans la deuxi\`eme partie. Dans la sixi\`eme
partie on d\'emontre les th\'eor\`emes donn\'es dans la partie deux en applicant
les th\'eor\`emes formelles et les th\'eor\`emes de convergence. La septi\`eme partie
est compl\`etement ind\'ependante et permet de contr\^oler certaine normes, en
vue de d\'emontrer l'existence globale.

\section{Th\'eor\`emes principaux}
Etant donn\'e  une fonction $2\pi$ p{\'e}riodique $u$. On note  $\f{u}_n$ son n
i{\`e}me coefficient de Fourier.  On note $H^k(\T)$ l'espace de Sobolev des  
fonctions r{\'e}elles $2\pi$-p{\'e}riodique,  $k$ fois  d{\'e}rivables au sens
distributif sur le tore $\T=\R/\Z$, muni de la 
norme:  
\[\parallel u \parallel_k=\sqrt{\sum_{i=0}^k \int_\T (\frac{\partial^i 
    u}{\partial x^i})^2}=\sqrt{\sum_{i=0}^k \sum_{n\in\N}\f{u}_{n}\f{u}_{-n}n^{2k}}.\] ${\cal H}^k(\T)$  le sous espace hilbertien de $H^k(\T)$ tel que $\f{u}_0=0$. La norme $\mid.\mid_k$ est donn{\'e}e par
\[\mid u \mid_k=\sum_{i=0}^k \sum_{n\in\N}\mid\f{u}_{n}n^{k}\mid,\]
et d\'efinit l'espace de Banach $W^k(\T)$.

On note $\R^{\N^*}$ (respectivement $\R^{\Z^*}$) les suites r{\'e}elles indic{\'e}es dans $\N^*$ ( respectivement dans $\Z^*$). $\N^\N$ d{\'e}signe l'ensemble des multi-entiers naturels et $\Z^\N$ d{\'e}signe l'ensemble des multi-entiers relatifs:
\begin{theoreme}\label{TfnKdVS} 
 
Soit un polyn{\^o}me $P\in \R[X,Y]$ de degr{\'e} $p$, $P(0,0)=\partial_xP(0,0)=\partial_yP(0,0)=0$, un r{\'e}el $\lambda$ et soit  $q$  le degr{\'e}  de $P$
en la variable $Y$. On suppose que $s>3$ ou si $q=1$ $s>2$, alors   pour tout $u_0\in
W^s(\T)$ il existe
$t_0(P,\lambda,\mid u_0\mid_2)$ tel que  pour tout $t\in [0,t_0]$, il existe une
unique solution locale $u(t)$  sur $W^s(\T)$ de l'{\'e}quation d'{\'e}volution:

\begin{equation*}
\bc \di
\frac{\partial u}{\partial 
  t}=\lambda 
\frac{\partial^{3} u}{\partial x^3}+P(u,\partial_{x}u).\\ 
u(x,0)=u_0\in W^s(\T).
\ec
\end{equation*}
De plus si $q=1$  alors $t_0=\frac{M(P,\lambda)}{(\mid u_0\mid_{2})^2(\mid
  u_0\mid_0)^{2p}}$. Et si il existe une solution  $u\in W^s(\T)$ definie sur $[0,t_1]$ ( $t_1\geq t_0$) alors pour tout  $t\in [0, t_1]$, et tout $k\leq 2$: \begin{eqnarray}
\mid u(t)\mid_k&\leq& \mid u_0\mid_k\left(M(P,\lambda)\mid u(0)\mid_0^{2p}+1\right)+ \label{Ine1}\\
&&(M(P,\lambda))^2\sup_{t'\in [0,t_1]}\left(\mid u(t')\mid_0^{2p}\parallel u(t')\parallel_{1}\mid u(t')\mid_k\right)\mid t\mid\nonumber.\end{eqnarray}

\end{theoreme} 
Nous n'imposons  aucune condition sur
la taille de $\mid u_0\mid_s$, {\`a} la
diff{\'e}rence du th{\'e}or{\`e}me de  Bourgain \cite{BOU1} o{\`u} $\mid
u_0\mid_s$, doit {\^e}tre petit.

\begin{corolaire}\label{CfnKdVS} 
 
Soit $s>2$, soit un polyn{\^o}me $P\in \R[X]$, $P(0)=\frac{d P(0)}{dx}=0$. Si $P$ est de degr{\'e} plus petit ou {\'e}gal {\`a} 3 ou si $P$ est de degr{\'e} p impair tel que $\lim_{x\to +\infty} P(x)=-\infty$, alors pour tout $u_0\in
W^s(\T)$ il existe  une
unique solution globale $u(t)$  sur $W^s(\T)$ de l'{\'e}quation d'{\'e}volution:

\begin{equation} \label{Eqkdv}
\bc
\frac{\partial u}{\partial 
  t}=\lambda 
\frac{\partial^{3} u}{\partial x^3}+\partial_{x}P(u).\\ 
u(x,0)=u_0\in W^s(\T).\ec
\end{equation}  
 
\end{corolaire} 

\begin{proof}
 Soit $u$ la  solution de l'{\'e}quation (\ref{Eqkdv}). On remarque que $\parallel
u(t)\parallel_0$ est  une constante ne d{\'e}pendant pas de $t$, ainsi que
$\parallel u\parallel_1^2-2\int_0^1Q(u)dx$ ( avec $Q'=P$). Nous allons montrer
d'abord que pour $P$  de degr{\'e} plus petit ou {\'e}gal que 3 et que pour  $P$
de degr{\'e} $p$ impair de coefficient de plus haut degr{\'e} negatif, alors
$\parallel u\parallel_1$ et $\mid u\mid_0$ sont born{\'e}es. Soit $P$  de
degr{\'e} plus petit ou {\'e}gal {\`a} 3, d'apr{\`e}s l'in{\'e}galit{\'e} de
Schwartz $\mid u\mid_0\leq \parallel u\parallel_{\frac{2}{3}}\sum_{n\in
N^*}\frac{1}{n^{\frac{4}{3}}}$ et $\parallel u\parallel_{\frac{2}{3}}^2\leq
\parallel u\parallel_1\parallel u\parallel_{\frac{1}{2}}$. Enfin $\parallel
u\parallel_{\frac{1}{2}}^2\leq \parallel u\parallel_1\parallel u\parallel_0$. On
a donc:
\[\mid u\mid_0\leq(\parallel u\parallel_1)^{\frac{3}{4}}\parallel u\parallel_0^{\frac{1}{2}}\sqrt{\sum_{n\in N^*}\frac{1}{n^{\frac{4}{3}}}}.\]
Or comme $\parallel u\parallel_1^2-2\int_0^1Q(u)dx=\parallel u_0\parallel_1^2-2\int_0^1Q(u_0)dx$ et $\parallel
u(t)\parallel_0=\parallel
u_0\parallel_0$, on en d{\'e}duit:
\[\parallel u\parallel_1^2\leq \left\vert\parallel u_0\parallel_1^2-2\int_0^1Q(u_0)dx\right\vert+\left\vert 2\int_0^1Q(u)dx\right\vert.\]
On peut donc {\'e}crire $\mid 2\int_0^1Q(u)dx\mid=\mid 2\int_0^1a_1u^3+a_2u^4dx\mid\leq (\mid a_1\mid+\mid a_2\mid)\parallel u_0\parallel_0^2)(\mid u\mid_0+\mid u\mid_0^2)$,  qui implique qu'il existe deux constantes $C$ et $D$ (ind{\'e}pendnantes de $t$) telles que:
\[\mid 2\int_0^1Q(u)dx\mid\leq \left(C\mid u\mid_0(\parallel u\parallel_1)^{\frac{3}{4}}+D\mid u\mid_0(\parallel u\parallel_1)^{\frac{3}{2}}\right).\] Donc $\parallel u\parallel_1$ est born{\'e}e, ce qui implique  que  $\mid u\mid_0$ est aussi born{\'e}e. 

Soit $P$  de degr{\'e} $p$ impair de coefficient de plus haut degr{\'e} n\'egatif. Il existe un intervalle $A$ ferm{\'e} born{\'e} de $\R$ tel que pour tout $x$ en dehors de cet intervalle $Q(x)\leq 0$, donc $\parallel u\parallel_1\leq \sup_A\mid Q(x)\mid +\mid\parallel u(0)\parallel_1-2Q(u(0))\mid$. Donc $\parallel u\parallel_1$ est born{\'e}e, ce qui implique  que  $\mid u\mid_0$ est aussi born{\'e}e. 
 
Si  $\mid u\mid_k$ est born{\'e}e alors il existe une solution globale, car si on note $t_i=\frac{\epsilon M(P,\lambda)}{(\mid u(t_{i-1})\mid_{\frac{3}{2}})(\mid
  u(t_{i-1})\mid_0)^{2p}}$ les $t_i$ sont minor{\'e}s, donc $\sum_{i\in \N}
  t_i=+\infty$ et d'apr{\`e}s le th{\'e}or{\`e}me \ref{TfnKdVS} les solutions
  $u(t)$ sont d{\'e}finis pour tout  $t\leq \sum_{i\in \N} t_i$, donc pour tout
  $t$ de $\R$. Supposons alors qu'il existe $t_1>0$ ( $t_1<\infty$ sinon la
  solution serait globale) tel que $\limsup_{t\to t_1}\mid
  u\mid_k=+\infty$. Alors pour tout $\epsilon$ positif ( $t_1>\epsilon>0$) et
  tout  $t\in [t_1-\epsilon, t_1[$ d'apr{\`e}s l'in{\'e}galit{\'e} (\ref{Ine1}): 
\[\mid u(t)\mid_k\leq \mid u(t_1-\epsilon)\mid_k\left(M(P,\lambda)\mid u(t_1-\epsilon)\mid_0^{2p}+1\right)+ \]\[(M(P,\lambda))^2\sup_{t'\in [t_1-\epsilon,t]}\left(\mid u(t')\mid_0^{2p}\parallel u(t')\parallel_{1}\mid u(t')\mid_k\right)\mid t-t_1+\epsilon\mid.\] 
Ce qui implique, en divisant par $\sup_{t'\in [t_1-\epsilon,t]}(\mid
u(t')\mid_k)$ et en passant {\`a} la limite sup, quand $t$ tend vers $t_1$, que   $\epsilon$  est major{\'e}. C'est absurde et on prouve ainsi le corollaire.
\end{proof}
On g\'en\'eralise  le th\'eor\`eme \ref{TfnKdVS} ainsi:
\begin{theoreme}\label{Tprincgeneral}
Soit un polyn{\^o}me $P\in \R[X_1,...X_k]$ de degr{\'e} $p$, $P(0,0)=\partial_{x_i}P(0)=0$, des r{\'e}els $\lambda_i$ et soit  $q_i$  le degr{\'e}  de $P$
en la variable $X_i$. On suppose que $s>2k+1$. Si $u_0\in
H^s(\T)$ et si dans un voisinage de $u_0$,
\[\lim_{n\to \infty}\Re(\sum_{j=1}^{2k+1}(in)^j(\lambda_{j}+\int_\T
\partial_{x_j}P(u)dx))\leq 0,\]
alors il existe
$t_0(P,\lambda,\parallel u_0\parallel_s)$ tel que  pour tout $t\in [0,t_0]$, il existe une
unique solution locale $u(t)$  sur $H^s(\T)$ de l'{\'e}quation d'{\'e}volution:

\begin{eqnarray*} 
\frac{\partial u}{\partial 
  t}&=&\sum_{i=1}^{2k+1}\lambda_{i} 
\frac{\partial^{i} u}{\partial x^i}+P(u,\partial_{x^1}u,...\partial_{x^k}u).\\ 
u(x,0)&=&u_0\in H^s(\T).
\end{eqnarray*}  
Par contre si \[\lim_{n\to \infty}\Re(\sum_{j=1}^{2k+1}(in)^j(\lambda_{j}+\int_\T
\partial_{x_j}P(u)dx))=+\infty,\] dans un voisinage de $u_0$, alors l'\'equation
est mal pos\'ee. 
\end{theoreme}

Ces r\'esultats se g\'en\'eralise \`a l'espace Hilbertien $H^s(\T)$:
\begin{theoreme}\label{TfnKdVSH} 
 
Soit un polyn{\^o}me $P\in \R[X,Y]$ de degr{\'e} $p$, $P(0,0)=\partial_xP(0,0)=\partial_yP(0,0)=0$, un r{\'e}el $\lambda$ et soit  $q$  le degr{\'e}  de $P$
en la variable $Y$. On suppose que $s>3+\frac{1}{2}$ ou si $q=1$ $s>\frac{5}{2}$, alors   pour tout $u_0\in
H^s(\T)\cap W^{\frac{5}{2}}(\T)$ il existe
$t_0(P,\lambda,\mid u_0\mid_{\frac{5}{2}})$ tel que  pour tout $t\in [0,t_0]$, il existe une
unique solution locale $u(t)$  sur $H^s(\T)$ de l'{\'e}quation d'{\'e}volution:

\begin{equation*}
\bc \di
\frac{\partial u}{\partial 
  t}=\lambda 
\frac{\partial^{3} u}{\partial x^3}+P(u,\partial_{x}u).\\ 
u(x,0)=u_0\in H^s(\T)\cap W^{\frac{5}{2}}(\T).
\ec
\end{equation*}
De plus si $q=1$  alors $t_0=\frac{M(P,\lambda)}{(\mid u_0\mid_{\frac{5}{2}})^2(\mid
  u_0\mid_0)^{2p}}$. Et si il existe une solution  $u\in H^s(\T)$ definie sur $[0,t_1]$ ( $t_1\geq t_0$) alors pour tout  $t\in [0, t_1]$, et tout $k\leq \frac{5}{2}$: \begin{eqnarray}
\mid u(t)\mid_k&\leq& \mid u_0\mid_k\left(M(P,\lambda)\mid u(0)\mid_0^{2p}+1\right)+ \label{Ine1}\\
&&(M(P,\lambda))^2\sup_{t'\in [0,t_1]}\left(\mid u(t')\mid_0^{2p}\parallel u(t')\parallel_{1}\mid u(t')\mid_k\right)\mid t\mid\nonumber.\end{eqnarray}

\end{theoreme} 
Nous n'imposons  aucune condition sur
la taille de $\parallel u_0\parallel_s$, {\`a} la
diff{\'e}rence du th{\'e}or{\`e}me de  Bourgain \cite{BOU1} o{\`u} $\parallel
u_0\parallel_s$, doit {\^e}tre petit.

\begin{corolaire}\label{CfnKdVSH} 
 
Soit $s>\frac{5}{2}$, soit un polyn{\^o}me $P\in \R[X]$, $P(0)=\frac{d P(0)}{dx}=0$. Si $P$ est de degr{\'e} plus petit ou {\'e}gal {\`a} 3 ou si $P$ est de degr{\'e} p impair tel que $\lim_{x\to +\infty} P(x)=-\infty$, alors pour tout $u_0\in
H^s(\T)\cap W^{\frac{5}{2}}(\T)$ il existe  une
unique solution globale $u(t)$  sur $H^s(\T)$ de l'{\'e}quation d'{\'e}volution:

\begin{equation} \label{EqkdvW}
\bc
\frac{\partial u}{\partial 
  t}=\lambda 
\frac{\partial^{3} u}{\partial x^3}+\partial_{x}P(u).\\ 
u(x,0)=u_0\in H^s(\T)\cap W^{\frac{5}{2}}(\T).\ec
\end{equation}  
 
\end{corolaire} 
On remarque, qu'on obtient ici un analogue  des r\'esultats d\'montr\'es par
Kenig, G.Ponce et L. Vega dans \cite{KPV}, pour le cas de
condition initiale sur des espace de type Sobolev.
\section{Formalisme des arbres}
Nous  exposons ici le formalisme  des arbres et leurs 
diff{\'e}rentes propri{\'e}t{\'e}s qui nous servirons plus tard.  Pour 
plus de pr{\'e}cisions se r{\'e}f{\'e}rer \` a \cite{HAR}, \cite{BOLL}. Cette
partie est constitu\'ee d'un grand nombre de d\'efinitions, qu'il est souvent
plus simple de ce repr\'esenter graphiquement. Malgr\`es les apparences, il
s'agit de notions tr\`es simple. Nous commen\c cons par rappeler les
d\'efinitions standard et donner les notations que nous utiliserons dans la suite. 
 \subsection{Rappel de th\'eorie des graphes}
\begin{definition}\label{Darbre} 
 
Soit $V$ un ensemble fini, un graphe G de $V$ est un couple (V,E) o{\`u} 
E est un ensemble de couples d'{\'e}l{\'e}ments de $V$. 
Les {\'e}l{\'e}ments de $V$ sont  appel{\'e}s noeuds. 
Si deux noeuds forment un couple de $E$, on dit qu'ils sont adjacents. On note $E_u$ l'ensemble des couples contenant $u$. On dit qu'un noeud est extr{\'e}mal s'il ne poss{\`e}de qu'un seul noeud 
adjacent. Une branche est un ensemble ordonn{\'e} de noeuds adjacents les uns
aux autres, c'est-{\`a}-dire: 
 
soit $u_i$ des noeuds distincts de $V$, $(u_1,...,u_i,...,u_p)$ forment une branche si  pour tout $i$, $u_i$ et $u_{i+1}$ sont 
adjacents. Cette branche est un cycle si et seulement si $u_1=u_p$. Un graphe est dit connexe si chaque noeud peut \^ etre reli{\'e} par 
une branche. 
\end{definition} 
\begin{definition} 
 
Un sous-graphe $G'$ d'un graphe G, not{\'e}  $G'\subset G$, est un graphe  
v{\'e}rifiant 
$V(G') \subset V(G)$ et 
$E(G') \subset E(G)$.
 
On d{\'e}finit les op{\'e}rateurs bool{\'e}ens: 
\begin{eqnarray*}
G' \cup G&=&\left(V(G')\cup V(G), E(G')\cup E(G)\right)\\ 
G \setminus G'&=&\left(V(G) \setminus V(G'), E(G)\setminus \{vv' \in E(G), v'\in
V(G')\}\right)
\end{eqnarray*}
\end{definition} 
 
\begin{definition}\label{DI2} 
 
Un arbre $T$ est un graphe connexe ne poss{\'e}dant aucun cycle. Ceci 
est {\'e}quivalent {\`a} dire que deux  noeuds quelconques sont reli{\'e}s 
par une et une seule branche. 

\end{definition}

\begin{definition}\label{DI3} 
 
Pour orienter l'arbre $T$ on fixe un noeud $r$ qu'on appelle  racine de T, et
on introduit  la relation  d'ordre partiel suivante : 
\[u \leq_{r} v \Leftrightarrow u \in (r,v),\] 
c'est-{\`a}-dire que pour passer de la racine $r$ au noeud $v$ on passe par 
 $u$: dans ce cas on dit que u est inf{\'e}rieur {\`a} $v$. Si $T$ est muni de
 cette relation, on dit que c'est un arbre orient{\'e}. On note $v$ les noeuds
 de $T$ et $V(T)$ l'ensemble de ses noeuds priv{\'e} de la racine. Soit  $u$ et
 $v$ deux noeuds de $T$, $(u,v)$ d{\'e}signe la branche reliant u 
{\`a} v (c'est-{\`a}-dire la branche d' extr{\'e}mit{\'e}s $u$ et $v$ ). 
\end{definition}

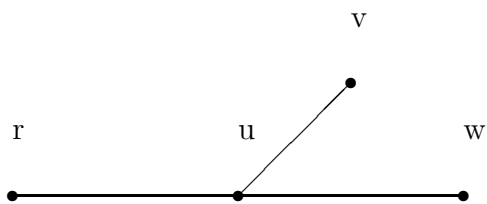
\setlength{\unitlength}{1.5cm} 
\begin{picture}(5,3) 
\put(1,1){\circle*{0.1}} 
\put(3,0){T} 
\put(1,1.5){r} 
\put(1,1){\line(1,0){2}} 
\put(3,1){\circle*{0.1}} 
\put(3,1.5){u} 
\put(3,1){\line(1,1){1}} 
\put(4,2){\circle*{0.1}} 
\put(4,2.5){v} 
\put(3,1){\line(1,0){2}} 
\put(5,1){\circle*{0.1}} 
\put(5,1.5){w} 
\end{picture} 
\\ 
Ici $r$ est la racine et $u \leq_r v$. $T$ muni de cette relation n'est pas totalement ordonn{\'e}. C'est en partie pour palier {\`a} ce d{\'e}faut qu'on introduit d'autres notions.

\begin{definition}\label{DI4} 
 
On appelle {\'e}tiquetage d'un arbre $T$, une bijection $e$ de  $V(T)$ sur
$\{1,...,\be\}$ ( o\`u $\be$ est le nombre d'\'el\'ements de $V(T)$). Le couple $(T,e)$ est 
appel{\'e} arbre {\'e}tiquet{\'e} par $e$. 
\end{definition}  
Dans la suite, les arbres \'etiquet\'es ne servent que pour la construction formelle
de l'application r\'eciproque, et donc pour d\'emontrer la non existence de
solution. Pour d\'emontrer la convergence, on ne travaille que sur des arbres
non-\'etiquet\'e.

On note $A_{v}$ l'ensemble des noeuds adjacents {\`a} $v$ et sup{\'e}rieurs {\`a} $v$ et $\beta_{v}$ d{\'e}signe le cardinal de $A_v$. On dit que $v$ est un noeud multiple si $\beta_v$ est plus grand que 1 et simple si $\beta_v=1$. Si $\be_v=0$, $v$ est un noeud extr{\^e}me. 

\begin{definition}\label{DI5} 
 
Un arbre est dit bien {\'e}tiquet{\'e},  si et seulement si  
 $u \leq_r v \Rightarrow e_u \leq e_v$. 
On {\'e}crit  
$u \leq_e v$  si $e_u \leq e_v$ . 
En d'autres termes l'{\'e}tiquetage doit respecter l'ordre de la 
racine. Dans la suite, on n'utilise que des arbres bien {\'e}tiquet{\'e}s. Pour un arbre {\'e}tiquet{\'e} $T$, on
note $e(T)$ son {\'e}tiquetage. 
\end{definition} 
 
\subsection{Arbres fruitiers}
Nous allons d\'evelopper un formalisme particulier des arbres, adapt\'e
sp\'ecifiquement \`a notre probl\`eme. Nous parlons d'arbres fruitier, ce qui
consite a associ\'e \`a chaque noeuds des entiers relatifs. Nous donnons aussi
les d\'efinitions de sous-arbres fruitiers et de diff\'erentes op\'erations
utiles , sur
ces arbres.

On  d{\'e}signe par $\ov{Z}^p$, l'ensemble des  p-uplet de $\Z$ modulo les permutations sur les indices, et on note $\ov{\Z}^\N=\bigcup_p\ov{Z}^p$.  
\begin{definition}\label{Dfruit} 
 
A chaque noeud, on peut associer un {\'e}l{\'e}ment de $\ov{\Z}^\N$ (ou l'ensemble vide),  
 qu'on appelle fruits.   
On note $\alpha_u$ 
le fruit de $u$ et on d{\'e}signe par sa taille, $\parallel 
\alpha_u\parallel$, l'entier $k$ tel que 
\[\alpha_u \in \ov{\Z}^k \hbox{ o{\`u} } k=\parallel\alpha_u\parallel.\]  
On d{\'e}signe par $\alpha_u^i$ les diff{\'e}rents {\'e}l{\'e}ments (entiers) de 
$\alpha_u$ (le uplet). Soit $\al$ l'entier d{\'e}fini par:

\[\sum_{u\in V(T)}\parallel \al_u\parallel=\al.\]
On d\'esigne par $\al$  le nombre d'\'el\'ements fruits de $T$.
On note $F(T)$  l'ensemble des fruits de $T$.
Le couple $(T,F(T))$ est appel{\'e} arbre fruitier, de fruit $F(T)$ et
$\al_v(T)$ d{\'e}signe le fruit sur $T$ au n\oe{}ud $v$, compos{\'e} des
entiers 
$\al_v^i(T)$.
\end{definition}

\begin{definition}\label{Ddegre}
Le degr{\'e} du noeud $v$, not{\'e}  $\gamma_v$ , est d{\'e}fini par :
\[\gamma_v=\parallel\alpha_v\parallel +\beta_v.\] 
 
 \end{definition}

Un arbre est un tronc s'il est constitu{\'e} que d'une seule branche, ainsi tous
ses noeuds  sont simples et l'{\'e}tiquetage ( s'il existe) est unique.

 D\'esormais, on dit qu'un arbre-fruitiers orient\'e, est  constitu{\'e} de $\beta$ noeuds 
 et $\alpha$ fruits, si $\alpha > \beta$ 
et 
$ \parallel \al_v \parallel \geq 2-\be_v$ et si la racine   $r$ n'a qu'un seul
 et unique noeud adjacent.

On d\'esigne par $T_{\alpha,\beta}$  l'ensemble des arbre-fruitiers orient\'e, non {\'e}tiquet{\'e}s, constitu{\'e}s de $\beta$ noeuds 
 et $\alpha$ fruits. 
  
$T_{\alpha,\beta}^1$ est l'ensemble des troncs fruitiers de $\be$ 
noeuds et $\al$ fruits orient{\'e}s mais non {\'e}tiquet{\'e}s. 

$T_{\alpha,\beta}^e$ est l'ensemble des arbres fruitiers de $\be$ 
noeuds et $\al$ fruits orient{\'e}s et bien {\'e}tiquet{\'e}s. 

$T_{\alpha,\beta}^{e,1}$ est l'ensemble des troncs fruitiers de $\be$ 
noeuds et $\al$ fruits orient{\'e}s et bien {\'e}tiquet{\'e}s. 

Enfin pour un arbre fruitier $(T,F(T))$ on dit que  $F(T)\in F(n)$ (le fruit de
$T$ est dans $F(n)$) si 
$\sum_{u\in V(T)}\sum_i\al_u^i=n$  et on d\'esigne par
$(T_{\alpha,\beta}^e,F(n))$ l'ensemble des  arbres fruitiers de $\be$ 
noeuds et $\al$ fruits orient{\'e}s et bien {\'e}tiquet{\'e}s dont le fruit est dans $F(n)$.

\begin{definition} 
 
On appelle sous-arbre  d'un arbre $T$, un arbre qui est un  
sous-graphe de $T$.  Un sous arbre orient{\'e} est un sous-arbre  
conservant la m{\^e}me orientation; c'est {\`a} dire que si on se donne $T'$ un  
sous arbre   de $T$  ( orient{\'e}), et $(u,v)$  
deux noeuds de $V(T')$, alors $u\leq_{r(T')}v$ si et seulement si  
$u\leq_{r(T)}v$. De m\^eme,  un sous-arbre bien \'etiquet\'e (et orient{\'e})
est un sous arbre bien \'etiquet\'e 
conservant la m{\^e}me orientation et le m\^eme ordre d'\'etiquetage, c'est {\`a} dire que si on se donne $T'$ un  
sous arbre   de $T$  (  orient{\'e}) et $(u,v)$  
deux noeuds de $V(T')$, alors $u\leq_{r(T')}v$ si et seulement si  
$u\leq_{r(T)}v$ et $u\leq_{e(T')}v$ si et seulement si  
$u\leq_{e(T)}v$.

\end{definition} 
Dans la suite on d{\'e}signe par sous arbre, un sous arbre orient{\'e} et
{\'e}tiquet{\'e}. On se propose d'{\'e}tendre la d{\'e}finition de sous-arbre aux arbres
fruitiers.  
 
\begin{definition}  
  
Soit $(T,F(T))$ un arbre fruitier, $(T',F(T'))$ est un sous  
arbres fruitier si et seulemnt si $T'$ est un sous arbre de $T$ et  
$F(T')$ est d{\'e}fini par:

Pour tout noeud $v$  de $T'$,  le fruit $\al_v(T')$ est {\'e}gal {\`a}
\[\al_v(T')=\al_v(T)\cup_{z\in A_v(T)\setminus V(T')}\{\sum_{\substack{u\geq w}\\: u\in V(T)\setminus V(T')}\sum_j\al_u^j\}.\]  
  
\end{definition}  

Pour les d\'emonstration de nos th\'eor\`emes de convergences nous utilisons plus
particuli\`erement certains types de sous-arbres, qui reviennent
r\'eguli\`erement. Nous allons donc leur associer la notation particuli\`ere suivante. 
\begin{definition}\label{DT(v)} 
 
$T(v)$ est le sous-arbre fruitier 
de $T$ constitu{\'e} des noeuds $u$ tels que
$v\leq_r u$.

\end{definition} 
On remarque que l'op{\'e}ration revient {\`a} couper au noeud $v$ puis \`a
relier la racine $r$ au noeud $v$,  \`a  conserver l'ordre de l'{\'e}tiquetage
(si l'arbre est \'etiquet\'e), et \`a associer les fruits pr{\'e}c{\'e}dents  aux m{\^e}mes noeuds. 
\begin{definition}\label{DT(u,v)}
Soit $u$ et $v$ deux noeuds d'un arbre $T$, alors $T(v,u)$ est, si $u>v$ (respectivement $v>u$), le sous arbre de $T$ constitu{\'e} des noeuds de $T$ sup{\'e}rieur ou {\`e}gaux {\`a}   $v$ (respectivement $u$), priv{\'e}  des noeuds  plus grand  ou {\'e}gaux {\`a} $u$ (respectivement $v$).
\end{definition}
Ceci implique que le fruit au noeud  $z$ adjacent inf{\'e}rieur {\`a} $u$ est donn{\'e} par $\al_z(T(v,u))=\al_z(T)\cup \{p_u(T)\}$ et que les fruits des autres noeuds sont identiques {\`a} celui sur $T$.
C'est-{\`a}-dire que:
\begin{equation*}
\bc
V(T(v,u))=\{w\in V(T(v)): w \not \in V(T(u))\},\\
\forall w\in V(T(v,u)): u\not \in A_w: \al_w(T(v,u))=\al_w(T),\\ 
\forall w\in V(T(v,u)):u \in A_w :\al_w(T(v,u))=\al_w(T)\cup\{p_u(T)\}.
\ec
\end{equation*}

 \begin{definition}\label{v(T)} 
 
Soit un arbre fruitier $T$, $v$ un noeud de $T$, alors $v(T)$ est l'ensemble des
sous-arbres fruitier de $T$  donn\'e par: 
\[v(T)=\{T':\exists u\in A_v, T'=T(u)\}.\] 
\end{definition} 
En fait, se sont les arbres extraits ayant pour premiers noeuds ceux
adjacents {\`a} $v$. 

On peut maintenant {\'e}tendre les op{\'e}rations bool{\'e}ennes  aux arbres  
 fruitiers.  Ceci  consite a d\'efinir l'intersection entre deux sous-arbres
comme un sous arbre de m\^eme graphe ( bien s\^ur) et dont le fruit est d\'eduits
des deux arbres dont il est l'intersection. De m\^eme la r\'eunion de deux arbres
consiste \`a attacher les deux arbres pour n'en former plus qu'un et d\'efinir,
ainsi un nouveaux  fruit.

On va aussi d\'efinir des arbres extraits, qui sont d\'efinis comme des arbres
dont on a retir\'e une partie de leurs noeuds. 

Ces op\'erations sont alors tr\`es utiles pour d\'emontrer des propri\'et\'es par
r\'ecurence (sur le nombre de nouds par exemple), ou construire formellement des
objets qui poss\`edent des propri\'et\'es r\'ecurentes. 
 
\begin{definition}\label{Dopboul}  
  
Soit ${\cal T}$ l'ensemble des sous arbres fruitiers  de $T$, on d{\'e}finit l'application $\cap$ de ${\cal T}^2$ dans
${\cal T}$ par:  

$\forall (T_1,T_2)\in{\cal T}^2$,   
$T_1\cap T_2$ est le sous arbre fruitier de $T$ dont le graphe est $G(T_1\cap T_2)=G(T_1)\cap G(T_2)$.

De m{\^e}me $(T_1,F(T_1))\subset (T_2,F(T_2))$ si et seulement  
si $T_1\subset T_2$ et le fruit du sous arbre de $T_2$ de graphe  $(G(T_1)\cap G(T_2))$ est le fruit de $T_1$ ( remarquons que
l'assertion $T_1\subset T_2$ est \'equivalente  si on se place sur ${\cal T}$).  

Enfin, si $T_1$ est un sous arbre de $T$, on d{\'e}finit l'ensemble des arbres fruitiers
$T'=T\setminus T_1$ par:
\begin{description}

\item[(i)]
Si le noeud adjacent  inf{\'e}rieur au premier noeud de $T_1$, $u$, exist:

$T'=T\setminus T_1$ est l'ensemble constitu\'e de l' arbre fruitier dont le graphe est d\'eduit du graphe 
form\'e des noeuds de $V(T)$ priv{\'e} de ceux de $V(T_1)$, o{\`u}  le noeud $u$
adjacent inf{\'e}rieur au premier noeud de $T_1$,   est reli{\'e}
{\`a} tous les noeuds adjacents sup{\'e}rieurs aux noeuds extr{\^e}mes de $T_1$;
 et dont le fruit  pour tout noeud $v$ de $V(T)\setminus V(T_1)$, distinct de $u$ (le
noeud adjacent inf{\'e}rieur au premier noeud de $T_1$), est
$\al_v(T')=\al_v(T)$ alors que le fruit de $u$ est
$\al_u(T')=\al_u(T)\bigcup_{v\in V(T_1)} \{\al_v(T)\}$ ( c'est {\`a} dire le
fruit de $u$ en $T$ union les fruit des noeud de $T_1$ sur $T$).

\item[(ii)] Si  le noeud adjacent  inf{\'e}rieur au premier noeud de $T_1$, $u$ n'existe
pas, alors $T'=T\setminus T_1$ est l'ensemble des sous arbre de $T$
d\'eduit des noeud sup\'erieur \`a ceux de $T_1$ ( dans ce cas il y a plusieurs arbres).
\end{description}
 Ainsi on notera pour tout noeud $v$ de $T$, $T\setminus \{v\}$ comme l'ensemble
des arbres  $T'$ obtenue par $T$ priv\'e du  sous arbre $T''$ de $T$ constitu\'e du
seul noeud $v$. Ce qui signifie que:
\[T\setminus \{v\}=T\setminus T''.\]

Enfin pour  deux sous arbres $T_1$ et $T_2$ de $T$,  s'il
existe, $T_1\cup T_2$ d\'efinit le sous arbre de $T$ constitu\'e du graphe $G(T_1)\cup G(T_2)$. On
remarque donc, que $T_1\cup T_2$ existe si et seulement si le graphe $G(T_1)\cup
G(T_2)$ est un sous graphe connexe de $G(T)$.

\end{definition} 
On peut ainsi extraire d'un arbre $T$ diff\'erent type de sous arbres. On peut
alors choisir de regrouper ceux qui sont ind\'ependants pour former une
d\'ecomposition d'un arbre $T$. Ainsi, il  suffit de montrer une
propri\'et\'e pour chaque sous arbre ``ind\'ependant'', dont la r\'eunion du
graphe forme $T$, pour d\'emontrer cette propri\'et\'e sur $T$. Ce qui nous
am\`ene a d\'efinir la d\'ecomposition d'un arbre $T$.
\begin{definition}\label{DdecompT}
Soit un arbre fruitier $T$, une d{\'e}compositon $G$ de $T$ est un ensemble de sous arbres fruitiers de $T$ distincts (quelques soient $T_1$ et $T_2$ de $G$, $T_1\cap T_2=0$) dont la r{\'e}union des graphes de ces sous arbres est le graphe de $T$ ( $\bigcup_{T'\in G}V(T')=V(T)$). On d{\'e}finit alors $G(\cal T)$ comme l'ensemble des d{\'e}compositions de $T$.
\end{definition}

On remarque que ces transformations s'op{\`e}rent aussi bien sur les 
arbres {\'e}tiquet{\'e}s que non-{\'e}tiquet{\'e}s. De plus, elle est 
compl{\`e}tement ind{\'e}pendante du choix de l'{\'e}tiquetage. Nous allons
maintenant d\'efinir des sous-arbres sp\'ecifiquement pour des arbres
\'etiquet\'es. Dont nous avons besoin pour montrer les th\'eor\`emes de non existence.
  
\begin{definition}\label{DTe(v)} 
 
Soit $(T,F)$ un arbre fruitier bien \'etiquet\'e, $v$ un noeud de $T$, on
d\'efinit $T_e(v)$ comme
 le sous arbre fruitier \'etiquet\'e  $(T',F(T'))$ de  $(T,F)$, dont $V(T')$ est constitu{\'e} des noeuds d'{\'e}tiquette inf{\'e}rieur ou {\'e}gal {\`a} celle de $v$ .

\end{definition} 
Donc, si on note $ A_u^e(v)$  l'ensemble des noeuds adjacent sup{\'e}rieur {\`a}
$u$ d'{\'e}tiquette plus grande que $v$, on a $A_u^e(v)=\{z:z\in A_u, v<_e z\}$ et
 $p_z= \sum_{w\geq_r z}\sum_i\al_w^i$. Alors $T'$ est un sous arbre de $T$ v{\'e}rifiant $V(T')=\{u:  
u\leq_ev\}$, et son ensemble fruit $F(T')$ est donn{\'e} pour tout noeud $u$ par  

\[\al_u(T')=\al_u(T)\bigcup_{z\in A_u^e(v)} \{p_z\}.\]

\begin{definition}\label{DS(T)} 
 
Soit un arbre $T^1\in T_{\al,\be}^{e}$ on dit que: 
$T'\in S(T^1)$, 
si et seulement si, il existe un {\'e}tiquetage $e'$  tel que  
l'on ait
$T^{''}=T^1$, o\`u 
$T^{''}$ est le nouvel arbre d{\'e}duit de $T'$ par ce nouvel {\'e}tiquetage.
 
\end{definition} 

\begin{definition}\label{Tv(e)} 
 
$T_v$ est la famille des sous arbres de $T$ constitu{\'e}e des parties connexes du  
sous graphe $G$ de $T$: 
\[V(G)=V(T)\setminus \{u, e_u \leq e_v\}.\] 
 
\end{definition}

\begin{definition}\label{DS(v(T))} 
 Soit $T$ un arbre d'{\'e}tiquetage $e$, $v$ un noeud de $T$.
$S(v(T))$ est l'ensemble des arbres {\'e}tiquet{\'e}s $T'$ de m{\^e}me graphe que $T$ tel que $v(T')=v(T)$.  

 Soit $T$ un arbre d'{\'e}tiquetage $e$, $v$ un noeud de $T$. On d\'efinit 
$S(T_v)$ comme  l'ensemble des arbres {\'e}tiquet{\'e}s $T'$ de m{\^e}me graphe que $T$ tels que $T^{'}_v=T_v$.  
\end{definition} 
\subsection{D{\'e}finitions  et propri{\'e}t{\'e}s des diviseurs}

 On reprend les notations introduites par
 Chierchia et Falcolini dans \cite{CHIEFA}. On d{\'e}finit deux type de diviseurs: le premier est donn{\'e} pour les arbres {\'e}tiquet{\'e}s  
\begin{definition}\label{Ddiv} 
 
Soit l'entier $p_v$ d{\'e}fini pour un arbre $T$ par
\[p_v(T)=\sum_{u \in V(T) ,v \leq_r u}\sum_i\al_u^i(T)=\sum_{u \in V(T(v))}\alpha_u^i,\] 
 
et soit $\Lambda$ une fonction donn{\'e}e, d{\'e}finie de $\C$ dans $\C$ et  
\[n=\sum_{u\in V(T)}\sum_i\al_u^i.\] 
On d\'efinit alors  $\delta_v^{\Lambda}(T)$, le diviseur de $T$ au noeud $v$
associ{\'e} {\`a} $\Lambda$, comme le nombre complexe: 
\[\delta_v^\Lambda(T)=\sum_{u \in V(T_e(v))}\sum_i\Lambda(\alpha_u^i(T_e(v)))-\Lambda(n).\] 
 
\end{definition}
Le deuxi{\`e}me type de diviseur est donn{\'e} pour les arbres non-{\'e}tiquet{\'e}s.
\begin{definition}\label{DDIVp}
Soit $T$ un arbre fruitier orient{\'e} ( non {\'e}tiquet{\'e}) on d{\'e}finit
au noeuds $v$ pour une application $\Lambda$, un diviseur
$D^\Lambda_v(T)$ par
\[D^\Lambda_v(T)=\sum_{u\geq v}\sum_i\Lambda(\al_u^i)-\Lambda(\sum_{u\geq v}\sum_i\al_u^i).\]
\end{definition} 
On remarque alors que les deux diviseurs sont, dans le cas o{\`u} $T$ est un  tronc de $\be$ noeuds (qui est le cas o{\`u} l'{\'e}tiquetage n'a aucune importance), reli{\'e}s par:
\[D_{u_i}^\Lambda(T)=-\delta_{u_{i-1}}^\Lambda(T)+\delta_{u_\be}^\Lambda(T),\]
o{\`u} $u_i$ est le noeud de $T$ d'{\'e}tiquette $i$ et $\delta_{u_{0}}^\Lambda(T)=0$.

\begin{proposition}{\label{Pdiv1}}

Soit  un tronc $T$. Alors
\[p_v(T)=\sum_{u\geq_r v}\sum_j\alpha_u^j,\] 
et 
\[\delta_v^\Lambda(T)=\sum_i\Lambda (\al_v^i)+\sum_{u\in A_v}\Lambda(p_u(T))+\sum_{u <_r v}\sum_i\Lambda(\alpha_u^i)-\Lambda(n). 
\] 
De plus on a $\delta_u^\Lambda(T)=\delta_u^\Lambda(T_e(u))$. 
 \end{proposition} 
\begin{proposition}\label{Pdutronc}
Soit un arbre $T$ un tronc orient{\'e} et {\'e}tiquet{\'e}, $u$ le dernier
noeud de $T$. On a :
\[\frac{1}{\delta_u^\Lambda(T)}\prod_{v\in V(T),v\neq u}\frac{1}{-\delta_v^\Lambda(T)+\delta_u^\Lambda(T)}=\prod_{v\in V(T)}\frac{1}{D_v^\Lambda(T)}.\]
\end{proposition}
\begin{proposition}\label{Pv(T)}
Soit un arbre $T$ de $\be$ noeud, $w_\be$ le dernier noeud de $T$ ( au sens de
l'{\'e}tiquetage). On a :
\[-\delta_v^\Lambda(T)+\delta_{w_\be}^\Lambda(T)= \sum_{T_i\in T_v} \delta_{w_{\be_i}}^\Lambda(T_i).\]
Et avec $z_i(u)$ le plus grand noeud de $T_i$ (au sens de l'{\'e}tiquetage),
d'{\'e}tiquette sur $T$ inf{\'e}rieur ou {\'e}gale {\`a} $e_u$. On a :
\[-\delta_u^\Lambda(T)+\delta_{w_\be}^\Lambda(T)= \sum_{T_i\in T_v} (\delta_{w_{\be_i}}^\Lambda(T_i)-\delta_{z_i(u)}^\Lambda(T_i) ).\]
\end{proposition}
\begin{proof}
\[\sum_{T_i\in T_v} \delta_{w_{\be_i}}^\Lambda(T_i)=\sum_{u >_e v }\sum_i\Lambda(\al_u^i(T))-\sum_{z\in A_x^e(v): x\leq_e v}\Lambda(p_z(T))\]
\bml=\sum_{u \in V(T)}\sum_i\Lambda(\al_u^i(T))-\Lambda(\sum_{u \in V(T)}\sum_i\al_u^i(T))+\\\Lambda(\sum_{u \in V(T)}\sum_i\al_u^i(T))-\sum_{u \leq_e v }\sum_i\Lambda(\al_u^i(T))\\-\sum_{z\in A_x^e(v): x\leq_e v}\Lambda(p_z(T))\end{multline*}
\[=-\delta_v^\Lambda(T)+\delta_{w_\be}^\Lambda(T).\]
Ainsi 
\[\sum_{T_i\in T_v} (\delta_{w_{\be_i}}^\Lambda(T_i)-\delta_{z_i(u)}^\Lambda(T_i) )=
\delta_{w_\be}^\Lambda(T)-\delta_v^\Lambda(T)-\sum_{T_i\in T_v}\delta_{z_i(u)}^\Lambda(T_i) ).\]
On conclut en  appliquant les d{\'e}finitions \ref{DTe(v)} et \ref{Ddiv} :
\[\sum_{T_i\in v(T)}\delta_{z_i(u)}^\Lambda(T_i) )=(\delta_u^\Lambda(T)-\delta_v^\Lambda(T)).\]
\end{proof}
\begin{proposition}\label{PdivT(u,v)}
Soit un arbre fruitier $T$ et deux noeuds $u>v$. Alors pour tout noeud $w$ de
$T(v,u)$. On a:
\[D_w^{\Lambda}(T(v,u))=D_w^{\Lambda}(T)-D_u^{\Lambda}(T).\]
\end{proposition}
\begin{proof}
En effet \[D_w^{\Lambda}(T)=\sum_i\sum_{z\geq
w}\Lambda(\al_z^i(T))-\Lambda(\sum_i\sum_{z\geq w}\al_z^i(T)), \] donc :
\bml D_w^{\Lambda}(T)-D_u^{\Lambda}(T)=\left(\sum_i\sum_{z\geq w}\Lambda(\al_z^i(T))-\Lambda(\sum_i\sum_{z\geq w}\al_z^i(T))\right)-\\\left(\sum_i\sum_{z\geq u}\Lambda(\al_z^i(T))-\Lambda(\sum_i\sum_{z\geq u}\al_z^i(T))\right),\end{multline*}
ce qui donne 
\bml D_w^{\Lambda}(T)-D_u^{\Lambda}(T)=\sum_i\sum_{z\geq w:z\not \in
V(T(u))}\Lambda(\al_z^i(T))+\\\Lambda(\sum_i\sum_{z\geq
u}\al_z^i(T))-\Lambda(\sum_i\sum_{z\geq w}\al_z^i(T)).\end{multline*}
D'o\`u la proposition.
\end{proof}

Dans la suite,  $\Lambda(x,t)$ est une fonction  de  $\R^2$, donn\'e sous la forme:
\[\Lambda(t,x)=\sum_{i=0}^{k+1}\lambda_i(t)x^i;\]
o\`u les $\lambda_i$  sont des fonctions localement continues dans un voisinage de
$0$. On notera alors :
\begin{eqnarray*}
L(t,x)&=&\int_0^t\Lambda(x,s)ds,\\
\Delta_u^\Lambda&=&\int_0^t\delta_u^\Lambda(s)ds,\\
{\cal D}_u^\Lambda&=&\int_0^t D_u^\Lambda(s)ds.
\end{eqnarray*}

On peut alors d\'efinir la notion de {\it noeud r\'esonant}.

\begin{definition}\label{DRes} 
 
Soit un arbre $T$ tel que $F(T)\in F(n)$, on dit qu'un noeud $v$ est
r{\'e}sonant s'il existe $v'\in A_v$ ( qui est l'ensemble des noeuds adjacents
sup\'erieur \`a $v$) tel que $p_v=p_{v'}$, ou si il existe $i$ tel que $\al_v^i(T)=p_v(T)$ ( definition \ref{Ddiv}). 
Inversement, un noeud u est dit v-r{\'e}sonant si et seulement si $u\in
A_v$ et $p_u=p_v$. Les arbres $T$ constitu{\'e}s   d'un seul noeud tel que $F(T)\in F(n)$ ( c'est {\`a} dire dans $(T_{\al,1}^e,F(n))$) sont dit ``noeuds r{\'e}sonants'', si leur unique noeud est un noeud r{\'e}sonant et on les note $R_p(n)^1$. Par extension, $Rp(F(T))^1=\bigcup_{\al_w\in
  F(T)}\bigcup_iRp(\al_w^i)^1$.
\end{definition} 
\subsection{Propri{\'e}t{\'e}s g{\'e}n{\'e}rales des op\'erateurs de {\it division}}

 Etant donn\'ee  $w$ une fonction de  $C_c(\R)$, on d\'efinit l'op\'erateur
 $S_w$ sur $C_c(\R)$, comme l'application qui a tout $h$ de $C_c(\R)$, associe
 $v\in C_c(\R)$ qui est d\'efinit par :  
 
\[v(t)=\left(\exp(-\int_0^tw(s)ds)\right)\left(\int_{0}^{t}h(x)\exp(\int_0^xw(s)ds)dx\right) =<S_w,h>.\]

La multiplication est alors la  loi de composition 

\[<S_{w_1}S_{w_2},h>=<S{w_1},<S_{w_2},h>>,\] 
et on note : 
\[<S_{w_2},1>=(S_{w_2}).\]
Etant donn\'es  des op\'erateurs  $S_{w_i}$, on note le produit {\`a} 
droite de ces op\'erateurs: 
\[\rov{\prod_{j}}  S_{w_j}=\rov{\prod_{j=1}^n }
  S_{w_j}=S_{w_1}...S_{w_{n-1}}S_{w_n}.\] 
On note aussi
\[(\rov{\prod_{j}}  S_{w_j})=<S_{w_1}...S_{w_{n-1}}S_{w_n},1>,\]
et 
\[\left\vert\rov{\prod_{j}}  S_{w_j}\right\vert=\left\vert<S_{w_1}...S_{w_{n-1}}S_{w_n},1>\right\vert.\]

 Dans la suite pour simplifier les notations et quand il n'y a pas
 d'ambiguit{\'e}, on notera pour $U(T)$ un sous ensemble de $V(T)$ (
 o{\`u} $T$ est un arbre fruitier de fruit $F(T)$) et $w_\be$ le noeud de $T$
 d'{\'e}tiquette la plus grande :
\begin{eqnarray*}\label{NO}
\lov{S}[U(T),\Lambda]&=&\lov{\prod_{v\in
 U(T)}}S_{\delta_v^\Lambda(T)}\\
\lov{{\cal S}}[U(T),\Lambda](u)&=&\lov{\prod_{v\in
 U(T)}}S_{\delta_v^\Lambda(T)}\prod_{v\in
U(T)}f_v\prod_{i}\f{u}_{\al_v^i}\\
\rov{S}[U(T),\Lambda]&=&(\prod_{v\in
 U(T)}\prod_j\exp(L(t,\al_v^j))))\\
&&S_{\delta_{w_\be}}\rov{\prod_{v\in 
U(T), e_v<e_{w_\be}}}S_{-\delta_{v}^\Lambda(T)+\delta_{w_\be}^\Lambda(T)}
\end{eqnarray*}
\begin{eqnarray*}
\rov{{\cal
  S}}[U(T),\Lambda](u)&=&(\prod_{v\in
 U(T)}f_v\prod_j\exp(L(t,\al_v^j))\f{u}_{\al_v^j}))\\
&&S_{\delta_{w_\be}}\rov{\prod_{v\in 
U(T), e_v<e_{w_\be}}}S_{-\delta_{v}^\Lambda(T)+\delta_{w_\be}^\Lambda(T)}
\end{eqnarray*}
Et s'il n'y a pas d'ambiguit{\'e} sur le choix de $u$ on ometera de l'{\'e}crire.


Les solutions formelles de notre probl\`eme sont obtenus sous la forme de poduit
d'op\'erateurs de division, d\'efinis eux m\^eme par des arbres \'etiquet\'es. On propose, au
cours de notre \'etude de regrouper les arbres 
modulo l'\'etiquetage. Nous donnons maintenant les propri\'et\'es de ces
regroupement :

\begin{lemme}\label{LS(T_v)}
Soit un arbre fruitier $T$ de $\be$ noeuds, $w_\be$ le dernier noeud de
$T$. Pour tout noeud $v$ de $V(T)$, on a :
\bml (\sum_{T'\in S(T_v)}\rov{\cal S}[V(T'),\Lambda])=\\<\rov{\cal S}[V(T_e(v)),\Lambda],\prod_{T_j\in T_v}(e^{-L\left(t,\sum_{u\in V(T_j)}\sum_l\al_u^l(T_j)\right)}\rov{\cal S}[V(T_j),\Lambda])>.\end{multline*}
\end{lemme}
\begin{proof}
Supposons que $T_v$ contienne $I$ arbres. Pour chaque $T_i$ de $T_v$ ( cf
d{\'e}finition \ref{Tv(e)} page \pageref{Tv(e)}), on note 
\[X_i^j=\delta_{v_j^i}^{\Lambda}(T_i)-\delta^{\Lambda}_{v_{\be_i}^i}(T_i),\] 
o{\`u} $v_j^i$  est le noeud de $T_i$ d'{\'e}tiquette $j$ sur $T_i$ et $\be_i$
le nombre de noeuds de $T_i$ ( donc l'{\'e}tiquette du dernier noeud). Soit
$u_\tau$ le noeud d'{\'e}tiquette $\tau$ sur $T'\in S(T_v)$ ( cf d{\'e}finition
\ref{DS(v(T))} \pageref{DS(v(T))}), si  $u_\tau\in T_i$ on d{\'e}finit
$j(i,\tau)$ comme la valeur de l'{\'e}tiquette de $u_\tau$ sur $T_i$. Et pour
tout $k\neq i$, $j$ v\'erifit $j(k,\tau)=j(k,\tau-1)$ ( qui est l'{\'e}tiquette
sur $T_k$ du dernier noeud de $T_k$ d'{\'e}tiquette sur $T'$ inf{\'e}rieur {\`a}
$\tau$). On note \[D^\Lambda(T_j)=\left(\sum_{u\in V(T_j)}\sum_l\Lambda(\al_u^l(T_j))-\Lambda\left(\sum_{u\in V(T_j)}\sum_l\al_u^l(T_j)\right)\right),\] 
et on applique la proposition \ref{Pv(T)} qui donne
\bml \sum_{T'\in S(T_v)}S_{\delta_{w_\be}^\Lambda}\rov{\prod_{\substack{u\in 
V_{\beta}(T'),\\e_u<e_w}}}S_{-\delta_{u}^\Lambda(T')+\delta_{w_\be}^\Lambda(T')}=e^{(\Delta_{w_\be}^\Lambda(t,T)-\Delta_{v}^\Lambda(t,T))}\times\\
<\rov{S}[V(T_e(v)),\Lambda],\left(\prod_{T_i\in T_v}e^{{\cal D}^\Lambda(t,T_i)}\right)\sum_{\sigma\in \sigma(X)}(\lov{\prod_{\tau=1}^I}S_{\sum_iX_i^{j(i,\tau)}})>.\end{multline*}
 En appliquant le lemme du
  comportement moyen \ref{lemCMgen}, on a :
\bml \sum_{T'\in S(T_v)}S_{\delta_{w_\be}^\Lambda}\rov{\prod_{\substack{u\in 
V_{\beta}(T'),\\e_u<e_w}}}S_{-\delta_{u}^\Lambda(T')+\delta_{w_\be}^\Lambda(T')}=\\e^{(\Delta_{w_\be}^\Lambda(t,T)-\Delta_{v}^\Lambda(t,T))}<\rov{S}[V(T_e(v)),\Lambda],\left(\prod_{T_i\in
T_v}e^{{\cal D}^\Lambda(T_i)}\right)\prod_{\tau=1}^I\rov{\prod_j}S_{X_i^j})>.\end{multline*}

Ceci   implique donc que :
\bml \sum_{T'\in S(T_v)}(\prod_{u\in 
V(T)}f_u\prod_{i}\exp(L(t,\al_u^i))\f{v}_{\al_u^i})S_{\delta_{w_\be}^\Lambda}\rov{\prod_{\substack{v\in 
V_{\beta}(T'),\\e_u<e_w}}}S_{-\delta_{v}^\Lambda(T')+\delta_{w_\be}^\Lambda(T')}=\\<\rov{\cal S}[V(T_e(v)),\Lambda],e^{-L\left(t,\sum_{u\in V(T_j)}\sum_l\al_u^l(T_j)\right)}\prod_{T_j\in
  T_v}\rov{{\cal
  S}}[V(T_j),\Lambda]>.\end{multline*}
\end{proof}
De m\^eme on montre le 
\begin{lemme}\label{Ltechcompmoyen}
Soit un arbre fruitier $T$ de $\be$ noeuds, $w_\be$ le dernier noeud de $T$. 
\[(\sum_{T'\in S(T)}\lov{\cal S}[V(T'),\Lambda])=\]\[\sum_{T'\subset\neq T}(\rov{\cal S}[V(T'),\Lambda])\prod_{T_j\in T\setminus T'}(e^{-L\left(t,\sum_{u\in V(T_j)}\sum_l\al_u^l(T_j)\right)}\rov{\cal S}[V(T_j),\Lambda]).\]
\end{lemme}
\begin{proof}On applique le th\'eor\`eme \ref{TCMgen}, o\`u on proc\`ede comme dans le lemme \ref{LS(T_v)}.
\end{proof}

\section{Existence formelle}

\begin{definition}{(Les uplets non ordonn{\'e}s)}\label{Duplet}
 
    On note $\ov{Z}^p$  l'ensemble des  p-uplet de $\Z$, modulo les permutations sur les indices, et on note $\ov{\Z}^\N=\bigcup_p\ov{Z}^p$.

Soit $n$ un entier relatif, $\al$ un entier naturel, on d{\'e}finit l'ensemble: 

\[\ov{\Z}^\al(n)=\{\un{k}\in \ov{\Z}^\al \hbox{ tel que }\sum_{i}k_i =n\}.\] 
Soit $\un{n}\in \ov{\Z}^\al$  la norme $\mid.\mid$ est donn{\'e}e par: $\mid
\un{n}\mid=\sum_i\mid n_i \mid$, et le support d'un uplet  
$\un{n}=(n_1,...n_i,..n_\al)$, not{\'e}  $[\un{n}]$, est  le nombre de $\mid n_i\mid $
distinct et non nul.
\end{definition}

Soit $Q$ un polyn{\^o}me de $\R[X_0,...,X_{k}]$. Pour tout $u\in H^s(\T)$, pour   $w=Q(u,\partial_x u,...\partial_{x^{k}} u)$ alors on note:\[\f{w}_n=\sum_\al\sum_{\un{n}\in \ov{\Z}^\al(n)}Q_{\un{n}}\prod_i\f{u}_{n_i}.\]

Donc si on se donne un arbre $T$, son fruit $\al_u(T)$ est un uplet non
 ordonn{\'e} $(\al_u^1,...,\al_u^i)$ et on notera
 $Q_{\al_u(T)}=Q_{(\al_u^1,...,\al_u^i)}$.  Pour le uplet ordonn{\'e}
 $(\al_1,...\al_j)=\al_v(T_e(v))$ ( qui est le fruit au noeud $v$ du sous arbre
 $T_e(v)$ de l'arbre $T$, cf d\'efinition \ref{DTe(v)}), on a donc
 \[Q_{\al_v(T_e(v))}=Q_{(\al_1,...\al_j)}.\]

\begin{lemme}\label{lemformel}
Soit  $P\in \R[X_0,..X_{k}]$ de degr{\'e} minimum sup{\'e}rieur {\`a} 2.
On d\'efinit pour tout $i$:
 \[\Omega_i(u)=\int_T\partial_{X_i}P (u,...,\partial_{x^{k}}u)dx.\]
On note
\[Q(u)=P (u,...,\partial_{x^{k}}u)-\sum_{i=0}^{k}\Omega_i(u)\partial_{x^i}u.\]

Pour un arbre $T$ de $(T_{\alpha,\beta}^e, F(n))$, on note $w_\be$ son
dernier noeud (d'{\'e}tiquette $\be$) et
$S_{\delta_{v}^\Lambda-\delta_{w_\be}^\Lambda(T)}$  l'op{\'e}rateur
associ{\'e} {\`a} la fonction
$\delta_v^\Lambda(T)-\delta_{w_\be}^\Lambda(T)$ (cf. d{\'e}finition
\ref{Ddiv}). Pour tout noeud $v$ de $T$ notons : 
\[f_{v}=Q_{\al_v(T_e(v))}.\] 
Enfin, soit $\chi$ donn{\'e}e par
\bml 
\f{\chi}_{n}(u,t)=e^{L(t,n)}\f{u}_n+\sum_{\al>2}\sum_{\beta=1}^{\mid \alpha \mid-1}\sum_{T\in
  (T_{\alpha,\beta}^e, F_\al(n))}(\\(\prod_{v\in V(T)}f_v\prod_{j}\exp(L(t,\al_v^j))\f{u}_{\al_v^i})S_{\delta_{w_\be}^\Lambda(T)}\rov{\prod_{\substack{v\in 
V_{\beta}(T),\\ e_v<e_w}}}S_{-\delta_{v}^\Lambda(T)+\delta_{w_\be}^\Lambda(T)}). 
\end{multline*} 
Alors   $\chi(v,t)$  
est formellement solution de l'{\'e}quation (\ref{Ebidule}).
\begin{equation}\label{Ebidule}
\bcd
\partial_tu=\sum_{i=0}^{2k+1}\lambda_i(t)\partial_x^iu+Q(u)\\
u(0,x)=v\in H^s(\T)
\ecd
\end{equation}

De plus,  si on note  $q$  le degr{\'e}
   de $P$ en $Y$, il existe un r\'eel $M$ tel que si $v$ est un noeud
   non-r{\'e}sonant, on a
\[\mid f_v\mid\leq M\sum_{\substack{(j_1,...j_{\gamma_v})\in
   \{0,1\}^{\gamma_v}\\ \sum_ij_i \leq q}}\prod_{i=1}^{\gamma_v}\mid \al_v^i(T_e(v))\mid^{j_i},\]
 et si $v$ est un noeud r{\'e}sonant, $f_v=0$
(cf. page \pageref{Dfruit} pour la d{\'e}finition de $\gamma_v$).
Enfin, dans le cas Hamiltonien ( c-a-d $P(u,\partial_xu)=\partial_xQ(u)$),
     si $v$ est un noeud r{\'e}sonant, alors  $f_v=0$.
\end{lemme}

\begin{proof}
Soit

\bml\f{\chi}_n=\sum_{\al,\be}\sum_{T\in
  (T_{\alpha,\beta}^e,
  F_\al(n))}(\prod_{v\in 
V(T)}f_v\prod_{i}\exp(L(t,\al_v^i))\f{v}_{\al_v^i})\\S_{\delta_{w_\be}^\Lambda(T)}\rov{\prod_{\substack{v\in 
V(T),\\e_v<e_w}}}S_{-\delta_{v}^\Lambda(T)+\delta_{w_\be}^\Lambda(T)}.\end{multline*}
 On remarque que
\bml\partial_t\f{\chi}_n=\sum_{\al,\be}\sum_{T\in
  (T_{\alpha,\beta}^e,
  F_\al(n))}(\partial_t(\prod_{v\in 
V(T)}f_v\prod_{i}\exp(L(t,\al_v^i))\f{v}_{\al_v^i})\\S_{\delta_{w_\be}^\Lambda(T)}\rov{\prod_{\substack{v\in 
V(T),\\e_v<e_w}}}S_{-\delta_{v}^\Lambda(T)+\delta_{w_\be}^\Lambda(T)}+ \\(\prod_{v\in 
V(T)}f_v\prod_{i}\exp(L(t,\al_v^i))\f{v}_{\al_v^i})\\\partial_tS_{\delta_{w_\be}^\Lambda(T)}\rov{\prod_{\substack{v\in 
V(T),\\e_v<e_w}}}S_{-\delta_{v}^\Lambda(T)+\delta_{w_\be}^\Lambda(T)}).\end{multline*} 
Ce qui donne d'apr{\`e}s  les notations de la page \pageref{NO},

\bml\partial_t\f{\chi}_n=\sum_{\al,\be}\sum_{T\in (T_{\alpha,\beta}^e,F_\al(n))}(\left(\sum_{v\in 
V(T)}\sum_i\Lambda(\al_v^i(T))-\delta_{w_\be}^\Lambda(T)\right)\rov{{\cal S}}[V(T),\Lambda](v)+\\\sum_{\al,\be}\sum_{T\in (T_{\alpha,\beta}^e,F_\al(n))}(\prod_{v\in 
V(T)}f_v\prod_{i}\exp(L(t,\al_v^i))\f{v}_{\al_v^i})\rov{\prod_{\substack{v\in 
V(T),\\e_v<e_w}}}S_{-\delta_{v}^\Lambda(T)+\delta_{w_\be}^\Lambda(T)}).\end{multline*}
 Pour $T\in (T_{\alpha,\beta}^e,F_\al(n))$, on applique le lemme \ref{LS(T_v)}
au premier noeud $v_1$, d'o{\`u} on d{\'e}duit que

\bml\sum_{T'\in S(v_1(T))}(\prod_{v\in 
V(T)}f_v\prod_{i}\exp(L(t,\al_v^i))\f{v}_{\al_v^i})\rov{\prod_{\substack{v\in 
V(T'),\\e_v<e_w}}}S_{-\delta_{v}^\Lambda(T')+\delta_{w_\be}^\Lambda(T')}=\\\prod_{T_j\in
  v_1(T)}\rov{{\cal
  S}}[V(T_j),\Lambda](v).\end{multline*}
Nous avons donc

\bml\partial_t\f{\chi}_n=\sum_{\al,\be}\sum_{T\in
  (T_{\alpha,\beta}^e,F_\al(n))}\Lambda(n)\rov{{\cal S}}[V(T),\Lambda](v)+\\\left(\sum_\al\sum_{T\in
  (T_{\alpha,1}^e,F_\al(n))}\prod_{v\in V(T)}f_v(T)\prod_i\exp(L(t,\al_{v}^i))\f{v}_{\al_{v}^i}\right)+\\\left(\sum_{\al',\be}\prod_{T'\in
  T_{\al',\be}}\rov{{\cal
  S}}[V(T'),\Lambda](v)\right),\end{multline*}
qui par d{\'e}finition des $f_v$ donne:
\[\partial_t\f{\chi}_{n}=(\Lambda(n))\f{\chi}_{n}+\sum_{\un{n}\in
  \{\ov{\Z}^\N(n)\}}Q_{\un{n}}\prod_i\f{\chi}_{n_i}.\]

\end{proof}  

De m\^eme on d\'emontre la r\'eciproque suivante.
\begin{lemme}
En reprenant les notation du lemme \ref{lemformel}, posons
\[\f{\psi}_{\alpha,n}=\sum_{\beta=1}^{\mid \alpha \mid-1}\sum_{T\in T_{\alpha,\beta}^e, F_\al(T)\in F_\al(n)}\lov{\prod_{v\in 
V_{\beta}(T)}}S_{\delta_{v}^\Lambda(T)}.(-f_{v}) 
\prod_{i}\f{u}_{\al_v^i} .\]
Alors  $\psi(u,t)$, o\`u
$u$ est la solution de   (\ref{Ebidule}), est solution formelle de l'\'equation
lin\'earis\'ee (\ref{Ebidule2}):
\begin{equation}\label{Ebidule2}
\bcd\partial_tv=\sum_{i=0}^{2k+1}\lambda_i(t)\partial_x^iv,\\
v(0,x)=u_0\in H^s(\T).
\ecd
\end{equation}

\end{lemme}
 
\begin{proof}

Formellement, on a
\bml\partial_t\f{\psi}_n=\partial_t\f{u}_n-\sum_{\al>1,\be}\sum_{T\in T_{\alpha,\beta}^e}(\sum_i\Lambda(\al_v^i)-\Lambda(n))\lov{{\cal S}}[V_\be(T),\Lambda]-\\\sum_{\al>1,\be}\sum_{T\in
  T_{\alpha,\beta}^e}\lov{{\cal
  S}}[\{V_\be(T),e_v<\be\},\Lambda](f_{v_\be}\prod_i\f{u}_{\al_{v_\be}^i})+\\\sum_{\al>1,\be}\sum_{T\in
  T_{\alpha,\beta}^e}\sum_{m\in F(T)}(\partial_{\f{u}_m}\lov{{\cal
  S}}[V_\be(T),\Lambda]\partial_t\f{u}_m).\end{multline*}
O\`u  $v_\be$ est le noeud d'\'etiquette  $\be$.
En utilisant le fait que
\[\partial_t\f{u}_m)=(\Lambda(m,t))\f{u}_m+\sum_\al\sum_{\un{m}\in
  \Z^\al(m)}Q_{\un{m}}\prod_i\f{u}_{m_i},\]
on remplace dans l'\'equation et on trouve:

\[\partial_t\f{\psi}_n^{\al,\be}=\sum_{\al,\be}\sum_{T\in
T_{\alpha,\beta}^e}(\Lambda(n))\lov{{\cal S}}[V_\be(T),\Lambda].\]
Donc  $\psi$ est formellement solution de (\ref{Ebidule2}).
\end{proof}

\section{Th\'eor\`eme de convergence}

\begin{definition}
Soit un arbre fruitier $T$ et $U(T)$ un sous-ensemble de $V(T)$ (c'est \`a
dire un ensemble de noeuds de $T$). Pour tout noeud $v$ de $U(T)$, on d\'esigne par $j^v$ un \'el\'ement de
$\R^{\parallel \al_v\parallel}$. Et on note $j^v=(j^v_1,..j^v_i,..j^v_{\parallel
\al_v\parallel})$. On note alors $j(U(T))$ l'ensemble des $j^v$ quand
$v$ d\'ecrit $U(T)$. Soit $J$ un sous ensemble fini de $\R$ et $r$ un
\'el\'ement de $J$:  on dit que $j(U(T))$ est dans $J(r)$ si pour tout $v$ les \'el\'ements de
$j^v$ prennent leurs valeurs dans $J$ et s'il existe un unique indice $i$ et un
unique noeud $v$ tel que $j_i^v=r$.

\end{definition}

 Dans toute la suite, $P(u,\partial_xu)$  est  un
polynome de degr{\'e} $p$. On note respectivement $q_2$ et $q_1$ le
degr{\'e} en $x_0$ et $x_1$ de $P(x_0,x_1)$.   De plus, $f_v$ donn{\'e}
par le lemme formel \ref{lemformel} et $s=\inf(q_1,2)$. Enfin pour
$M(P,\Lambda)$ un r\'eel positif donn\'e on note 
 \[C(M)=
\frac{M(P,\Lambda)}{\epsilon^{q}}(\mid u_0\mid_{s})^{2q_1}(\mid
u_0\mid_{0})^{2q_2}.\]

\begin{lemme}\label{Lt^k}
Soit $r\geq s$, $u$ et $u_0$ des {\'e}l{\'e}ments de $
W_1^{r}(\T)$. Il existe $M(P,\Lambda)>0$ ind{\'e}pendant de $u$ et $T$ tel que pour tout $\epsilon$,  tout arbre $T$  de $\be$ n\oe{}uds et $\al$ fruits, $F(T)\in F(n)$,
si on note $k(T)=\textup{card}(\{v\in V(T):\mid D_v^{\Lambda}(T)\mid\leq C(M)\})\leq \be$,  pour tout $t\leq C(M)$:
\bs\bml\left\vert t^{k(T)}\prod_{v\in V(T)}f_v \prod_{v: \mid
    D_v^\Lambda\mid>C(M)}\frac{1}{\sup(1,\mid
    D_v^\Lambda\mid)}\right\vert\prod_v\prod_i\mid\f{u}_{\al_v^i(T)}\mid\leq \frac{1}{\mid n\mid^{r}}(C(M)t)^{k(T)}\\\times\frac{\mid u_0\mid_r}{\mid u_0\mid_0}\times\left(\sum_{\substack{j(V(T))\in
  \{0,1,2,r\}(r)\\
  \forall v :\sum_{i=1}^{\parallel \al_v\parallel}j_i\leq q_1-1+r }
}\prod_{v\in 
V(T)}\left(\prod_{i=1}^{\parallel
  \al_v\parallel}\mid\frac{\epsilon\mid\al_v^i\mid^{j_i^v}}{\mid
  u_0\mid_{j_i^v} }\f{u}_{\al_v^i}\mid\right)\right) .\end{multline*}\es
Dans le cas Hamiltonien, on a
\bs\bml\left\vert t^{k(T)}\prod_{v\in V(T)}f_v \prod_{v: \mid
    D_v^\Lambda\mid>C(M)}\frac{1}{\sup(C(M),\mid
    D_v^\Lambda\mid)}\right\vert\prod_v\prod_i\mid\f{u}_{\al_v^i(T)}\mid\leq\frac{1}{\mid n\mid^{r}}(C(M)t)^{k(T)}\\\times\frac{\mid u_0\mid_r}{\mid u_0\mid_0}\left(\sum_{\substack{j(V(T))\in
  \{0,1,r\}(r)\\
  \forall v :\sum_{i=1}^{\parallel \al_v\parallel}j_i\leq  r}
}\prod_{v\in 
V(T)}\left(\prod_{i=1}^{\parallel
  \al_v\parallel}\mid\frac{\epsilon\mid\al_v^i\mid^{j_i^v}}{\mid
  u_0\mid_{j_i^v} }\f{u}_{\al_v^i}\mid\right)\right) .\end{multline*}\es
\end{lemme}
\begin{proof}

Pour tout arbre $T$, on a
\bs\bml\left\vert t^{k(T)}\left(\prod_{v\in V(T)}f_v\right)\left(\prod_{v: \mid
    D_v^\Lambda\mid>C}\frac{1}{\sup(C,\mid
    D_v^\Lambda\mid)}\right)\right\vert(t)=\\
\left\vert (Ct)^{k(T)}\prod_{v\in V(T)}\frac{f_v}{\sup(C,\mid
    D_v^\Lambda\mid)}\right\vert(t).\end{multline*}\es

On applique le lemme \ref{lemccv3}, qui donne

\bs\bml\left\vert (Ct)^{k(T)}\prod_{v\in V(T)}\frac{f_v}{\sup(C,\mid
    D_v^\Lambda\mid)}\right\vert(t)\leq
     (Ct)^{k(T)}\frac{e^{4\al}\mid \al_w^i\mid^{r-2}}{(\mid n\mid)^{r} }\left(\prod_{u\in
     V(T)}\frac{1}{\sqrt{C}}\right)\times\\
\left(\prod_{v\in
V(T)} M \sum_{\substack{(j_1,...j_{\parallel \al_v\parallel})\in
  \{0,1,2\}^{\parallel \al_v\parallel}\\ \sum_{i=1}^{\parallel
    \al_v\parallel}j_i\leq q+1}}\prod_{i=1}^{\parallel
\al_v\parallel}\mid\al_v^i\mid^{j_i}\right).\end{multline*}\es

En rempla\c cant $C$ par sa valeur, on termine la preuve du lemme.

\end{proof}
On \'etant maintenant le lemme \ref{Lt^k}, qui porte sur le produit de diviseurs d'un
arbre $T$, au produit d'op\'erateurs de divisions ( cf page \pageref{LS(T_v)}).

\begin{theoreme}{\label{TCgen}} 
Soit $r\geq s$, $u$ et $u_0$ des {\'e}l{\'e}ment de $
W^{r}(\T)$. Il existe $M(P,\Lambda,r)>0$ ind{\'e}pendant de $u$ et $T$ tel que pour tout $\epsilon$,  pour tout arbre $T$ ( $F(T)\in F(n)$),   pour tout $t\leq C(M)$ , il existe
$k\leq \be$  et un fruit  $\al_w^l(T)$ tel que:

\bs\bml\left\vert\rov{\cal S}[V(T),\Lambda]\right\vert(t)\leq\frac{1}{\mid n\mid^{r}}(C(M)t)^k\\\times\frac{\mid u_0\mid_r}{\mid u_0\mid_0}\left(\sum_{\substack{j(V(T))\in \{0,1,2,r\}(r)\\ \sum_{i=1}^{\parallel \al_v\parallel}j_i\leq q_1+r }  }\prod_{v\in 
V(T)}\left(\prod_{i=1}^{\parallel \al_v\parallel}\mid\frac{\epsilon\mid\al_v^i\mid^{j_i^v}}{\mid u_0\mid_{j_i^v} }\f{u}_{\al_v^i}\mid\right)\right) ,\end{multline*}\es
Et dans le cas Hamiltonien:
\bs\bml\left\vert\rov{\cal S}[V(T),\Lambda]\right\vert(t)\leq\frac{1}{\mid n\mid^{r}}(C(M)t)^k\\\times\frac{\mid u_0\mid_r}{\mid u_0\mid_0}\sum_{\substack{j(V(T))\in \{0,1r\}(r)\\ \sum_{i=1}^{\parallel \al_v\parallel}j_i\leq r -1}  }\left(\prod_{v\in 
V(T)}\left(\prod_{i=1}^{\parallel \al_v\parallel}\mid\frac{\epsilon\mid\al_v^i\mid^{j_i^v}}{\mid u_0\mid_{j_i^v} }\f{u}_{\al_v^i}\mid\right)\right) .\end{multline*}\es
\end{theoreme}  
\begin{proof}
La d{\'e}monstration s'obtient en appliquant le  lemme 
\ref{lemaqD} avec
 $k(T)=\textup{card}(\{v\in V(T):\mid D_v^{\Lambda}(T)\mid\leq C(M)\})\leq \be$:
\[\left\vert\rov{S}[V(T),\Lambda]\right\vert(t)\leq \sum_{G\in G(\cal
  T)}\prod_{T'\in G} t^{k(T')}\prod_{v\in V(T'):\mid
  D_v^{\Lambda}(T')\mid> C(M)}\frac{1}{\mid D_v^{\Lambda}(T')\mid},\] (cf. d{\'e}finition \ref{DdecompT} de $G(\cal T)$), puis le lemme
  \ref{Lt^k} sur chacun des $T'$. On obtient la preuve en
  faisant le produit des sous-arbres $T'$ de $T$.
\end{proof}

Ce th{\'e}or{\`e}me montre que le crit{\`e}re de convergence se d{\'e}montre
 ind{\'e}pendemment de la taille de $\mid u\mid_k$, si $t$ est
 suffisement petit. C'est ainsi qu'on d{\'e}montre la convergence de
 l'application formelle $\phi$.
 
Dans le cas g\'en\'eral, on d\'emontre de m\^eme en applicant  le th\'eor\`eme \ref{TAG(T)1}  (au lieu du lemme \ref{lemaqD}),  que:

\begin{theoreme}\label{Tconvgen1}
Soit  $P(u,\partial_xu,...\partial_x^ku)$   un
polynome de degr{\'e} $p$. On note  $q_i$  le
degr{\'e} en $x_i$  de $P(x_0,x_1,...x_k)$.   On suppose que la partie r\'eelle 
$\Re(\Lambda(n,t))$
est positif pour $n$ suffisement grand d\'efinit pour tout  $t\in [0,t_1]$. De
plus, on suppose que les $\lambda_i$ sont $C^1$ et  born\'ees, sur $[0,t_1]$ ($\sup_{t\in [,t_1]}(\sum_i\mid
\lambda_i(t)\mid)\leq \Lambda$). Soit $s=2k+1$, $f_v$ donn{\'e}
par le lemme formel \ref{lemformel}. Soit $r\geq s$, $u$ et $u_0$ des {\'e}l{\'e}ment de $
W^{r}(\T)$. Alors, il existe $M(P,r)>0$ ind{\'e}pendant de $u$ et $T$ tel que pour tout $\epsilon$ un r{\'e}el positif, si on note 

\[C=
\frac{M(P,r)}{\epsilon^{p}}(\mid u_0\mid_{s})^{2p}\Lambda^2,\]
 
 alors pour tout arbre $T$ ( $F(T)\in F(n)$),   pour tout $t\leq \inf(\frac{1}{C},t_1)$ , il existe
$m\leq \be$  et un fruit  $\al_w^l(T)$ tel que

\bs\bml\left\vert\rov{\cal S}[V(T),\Lambda]\right\vert(t)\leq\frac{1}{\mid n\mid^{r}}(Ct)^k\\\times\frac{\mid u_0\mid_r}{\mid u_0\mid_0}\left(\sum_{\substack{j(V(T))\in \{0,1,2,...,r\}(r)\\ \sum_{i=1}^{\parallel \al_v\parallel}j_i\leq p-1+r }  }\prod_{v\in 
V(T)}\left(\prod_{i=1}^{\parallel \al_v\parallel}\mid\frac{\epsilon\mid\al_v^i\mid^{j_i^v}}{\mid u_0\mid_{j_i^v} }\f{u}_{\al_v^i}\mid\right)\right) .\end{multline*}\es
\end{theoreme}
En utilisant le lemme \ref{Ltechcompmoyen} et le th\'eor\`eme  \ref{TAG(T)2} on
en d\'eduit, dans le cas d'un produit \`a gauche ( qui montre l'existence de
$\phi^{-1}$ et nous donnera la non existence de solution dans le cas o\`u $\Re(\Lambda(n,t))$
est n\'egatif), le th\'eor\`eme similaire suivant.

\begin{theoreme}\label{Tconvgen2}
Soit  $P(u,\partial_xu,...\partial_x^ku)$   un
polyn\^ome de degr{\'e} $p$. On note  $q_i$  le
degr{\'e} en $x_i$  de $P(x_0,x_1,...x_k)$.   On suppose que $\Re(\Lambda(n,t))$
est n\'egatif pour $n$ suffisament grand. Soit $s=2k+1$ et  $f_v$ donn{\'e}e
par le lemme formel \ref{lemformel}. Soit $r\geq s$, $u$ et $u_0$ des {\'e}l{\'e}ment de $
W^{r}(\T)$. Il existe $M(P,\Lambda,r)>0$ ind{\'e}pendant de $u$ et $T$ tel que pour tout $\epsilon$ un r{\'e}el positif, si on note 

\[C=
\frac{M(P,\Lambda)}{\epsilon^{p}}(\mid u_0\mid_{s})^{2p},\]
 
 alors pour tout arbre $T$ ( $F(T)\in F(n)$),   pour tout $t\leq \frac{1}{C}$ , il existe
$m\leq \be$  et un fruit  $\al_w^l(T)$ tel que

\begin{small}\bml\left\vert e^{L(n,t)}\lov{\cal S}[V(T),\Lambda]\right\vert(t)\leq\frac{1}{\mid n\mid^{r}}(Ct)^k\\\times\frac{\mid u_0\mid_r}{\mid u_0\mid_0}\left(\sum_{\substack{j(V(T))\in \{0,1,2,...,r\}(r)\\ \sum_{i=1}^{\parallel \al_v\parallel}j_i\leq p-1+r }  }\prod_{v\in 
V(T)}\left(\prod_{i=1}^{\parallel \al_v\parallel}\mid\frac{\epsilon\mid\al_v^i\mid^{j_i^v}}{\mid u_0\mid_{j_i^v} }\f{u}_{\al_v^i}\mid\right)\right) .\end{multline*}\end{small}
\end{theoreme}

Ces diff\'erent r\'esultat suppose que l'on travaille sur $W^r(\T)$. Nous allons
\'etendre notre \'etude aux cas o\`u les condition initiale sont dans
$H^r(\T)$. Le but est alors de montrer que $\phi(u_0,t)$ est aussi dans
$H^k(\T)$. On distingue dans ce cas le cas o\`u $T$ est faiblement r\'esonant
(c-a-d qu'un de ces fruits est \'egal \`a la somme de tout ses fruits). On
proc\`ede de la m\^eme fa\c con que dans la d\'emonstration du lemme \ref{Lt^k} et
du th\'eor\`eme \ref{TCgen}.
\begin{lemme}\label{Lnormedouble}
   Soit  $s=\inf(q_1+\frac{1}{2},\frac{5}{2})$,  $r\geq s$, $u$ et $u_0$ des {\'e}l{\'e}ments de $
H^{r+\frac{1}{2}}(\T)$. Il existe $M(P,\Lambda)>0$ ind{\'e}pendant de $u$ et $T$ tel que pour tout $\epsilon$ un r{\'e}el positif,  pour tout arbre $T$  de $\be$ n\oe{}ud et $\al$ fruits, $F(T)\in F(n)$ tel que tous les
 fruits de $T$ soit distinct de $n$, si on note
 $k(T)=\textup{card}(\{v\in V(T):\mid D_v^{\Lambda}(T)\mid\leq C(M)\})\leq \be$,  pour tout $t\leq C(M)$:
\bs\bml\left\vert t^{k(T)}\prod_{v\in V(T)}f_v \prod_{v: \mid
    D_v^\Lambda\mid>C(M)}\frac{1}{\sup(1,\mid
    D_v^\Lambda\mid)}\right\vert\prod_v\prod_i\mid\f{u}_{\al_v^i(T)}\mid\leq\frac{1}{\mid n\mid^{r}}(C(M)t)^{k(T)}\\\times\frac{\mid u_0\mid_{r-1}}{\mid u_0\mid_0}\times\sum_{\substack{j(V(T))\in
  \{0,\frac{1}{2},2,,2+\frac{1}{2},r-1\}(r-1)\\
  \sum_{i=1}^{\parallel \al_v\parallel}j_i\leq q_1+2+r }
}\left(\prod_{v\in 
V(T)}\left(\prod_{i=1}^{\parallel
  \al_v\parallel}\mid\frac{\epsilon\mid\al_v^i\mid^{j_i^v}}{\mid
  u_0\mid_{j_i} }\f{u}_{\al_v^i}\mid\right)\right) .\end{multline*}\es
Dans le cas Hamiltonien, on a
\bs\bml\left\vert t^{k(T)}\prod_{v\in V(T)}f_v \prod_{v: \mid
    D_v^\Lambda\mid>C(M)}\frac{1}{\sup(C(M),\mid
    D_v^\Lambda\mid)}\right\vert\prod_v\prod_i\mid\f{u}_{\al_v^i(T)}\mid\leq\frac{1}{\mid n\mid^{r}}(C(M)t)^{k(T)}\\\times\frac{\mid u_0\mid_{r-1}}{\mid u_0\mid_0}\sum_{\substack{j(V(T))\in
  \{0,\frac{1}{2},1,1+\frac{1}{2},r-1\}(r-1)\\
  \sum_{i=1}^{\parallel \al_v\parallel}j_i^v\leq r }
}\left(\prod_{v\in 
V(T)}\left(\prod_{i=1}^{\parallel
  \al_v\parallel}\mid\frac{\epsilon\mid\al_v^i\mid^{j_i^v}}{\mid
  u_0\mid_{j_i^v} }\f{u}_{\al_v^i}\mid\right)\right) .\end{multline*}\es
\end{lemme}

\begin{proof}
On note $u$ le premier noeud de $T$. D'apr\`es le
lemme \ref{ARO4}, soit $\mid D_u^\Lambda(T)\mid$ est plus grand que $\frac{n^2}{2}$,
soit il existe trois fruits de $T$, $\al_1,\al_2,\al_3$, distincts du plus grand
des fruits de $T$ et dont le produit est plus
grand que $\frac{n^2}{\al}$. Dans le premier cas, on applique le lemme \ref{Lt^k}  au sous arbre
de $v(T)$ ( les arbres constitu\'es des noeuds sup\'erieurs \`a $u$) qui conclut
le
lemme puisque  $\mid D_u^\Lambda(T)\mid>\frac{n^2}{2}$.
 
Dans le deuxi\`eme cas on applique le lemme \ref{Lt^k} \`a $T$, et on conclut en
 majorant $\mid n\mid$ par $\sqrt{\al}\mid \al_1\al_2\al_3\mid^{\frac{1}{2}}$.
\end{proof}
On en d\'eduiot le
\begin{theoreme}\label{Tnormpar}
Soit    $s=\inf(q_1+\frac{1}{2},\frac{5}{2})$,  $r\geq s$, $u$ et $u_0$ des {\'e}l{\'e}ment de $
H^{r+\frac{1}{2}}(\T)$. Il existe $M(P,\Lambda)>0$ ind{\'e}pendant de $u$ et $T$ tel que pour tout $\epsilon$ un r{\'e}el positif, pour tout arbre $T$ ( $F(T)\in F(n)$) tel que tous les fruits de $T$ sont
distinct de $n$,   pour tout $t\leq C(M)$ , il existe
$k\leq \be$  et un fruit  $\al_w^l(T)$ tel que

\bs\bml\left\vert\rov{\cal S}[V(T),\Lambda]\right\vert(t)\leq\frac{1}{\mid n\mid^{r}}(C(M)t)^k\frac{\mid u_0\mid_{r-1}}{\mid u_0\mid_0}\\\times\sum_{\substack{j(V(T))\in \{0,\frac{1}{2},1,2,1+\frac{1}{2},2+\frac{1}{2},r-1\}(r-1)\\ \sum_{i=1}^{\parallel \al_v\parallel}j_i^v\leq q_1+2+r }  }\left(\prod_{v\in 
V(T)}\left(\prod_{i=1}^{\parallel \al_v\parallel}\mid\frac{\epsilon\mid\al_v^i\mid^{j_i^v}}{\mid u_0\mid_{j_i^v} }\f{u}_{\al_v^i}\mid\right)\right) ,\end{multline*}\es
Dans le cas Hamiltonien, on a
\bs\bml\left\vert\rov{\cal S}[V(T),\Lambda]\right\vert(t)\leq\frac{1}{\mid n\mid^{r}}(C(M)t)^k\frac{\mid u_0\mid_{r-1}}{\mid u_0\mid_0}\\\times\sum_{\substack{j(V(T))\in \{0,\frac{1}{2},1,1+\frac{1}{2},r-1\}(r-1)\\ \sum_{i=1}^{\parallel \al_v\parallel}j_i^v\leq r }  }\left(\prod_{v\in 
V(T)}\left(\prod_{i=1}^{\parallel \al_v\parallel}\mid\frac{\epsilon\mid\al_v^i\mid^{j_i^v}}{\mid u_0\mid_{j_i^v} }\f{u}_{\al_v^i}\mid\right)\right) ,\end{multline*}\es

\end{theoreme}
\begin{proof} 
On proc\`ede comme dans le th\'eor\`eme \ref{TCgen}, en remarquant que pour une
d\'ecomposition $G$ de $G(\cal T )$, si le sous
arbre $T'$ de la d\'ecomposition de $T$, contenant le premier noeud, a un fruit
\'egal \`a $n$, ce n'est pas un fruit de $T$. Il existe donc un autre arbre de
la d\'ecomposition qui a son fruit dans  $F(n)$. En recommencant ce
proc\'ed\'e, on montre qu'il existe un arbre de la
d\'ecompositon de $T$ qui a son fruit dans $F(n)$ et tel que tous ses fruits
soient distincts de $n$. Sur ce sous arbre, on applique le lemme
\ref{Lnormedouble} et sur les autres, le lemme \ref{Lt^k} qui donne la r\'eponse.
\end{proof}

On a aussi le th\'eor\`eme suivant, dans le cas ou $T$ est un arbre dont l'un des fruit a pour valeur la
somme de ses fruits, qu'on note $Rp(n)$ si $F(T)\in F(n)$.
\begin{theoreme}\label{TnormparR(n)}
Soit $T$ tel que son fruit $F(T)$ soit dans $F(n)$ et dont l'un des fruits a
pour valeur $n$. Soit $j(V(T)\in \{0,1,2,1+\frac{1}{2},2+\frac{1}{2},r\}(r)$,
 $N(u_0,j_i^v)$ d\'efinit la norme $\mid u_0\mid_{j_i^v}$ si $j_i^v$ est
distinct de $r$, et $\parallel u_0\parallel_r$ sinon. On a alors
\bs\bml\mid\mid\rov{S}[V(T),\Lambda]\mid(t)\mid\leq \frac{1}{\mid n\mid^{r}}(Ct)^k\frac{\parallel u_0\parallel_r}{\mid u_0\mid_0}\\\times\sum_{\substack{j(V(T))\in
\{0,1,2,1+\frac{1}{2},2+\frac{1}{2},r\}(r)\\
\sum_{i=1}^{\parallel \al_v\parallel}j_i^v\leq q_1+\frac{3}{2}+r }
}\left(\prod_{v\in 
V(T)}\left(\prod_{i=1}^{\parallel
\al_v\parallel}\mid\frac{\epsilon\mid\al_v^i\mid^{j_i^v}}{N( u_0,j_i^v)
}\f{u}_{\al_v^i}\mid\right)\right).\end{multline*}\es

\end{theoreme}
\begin{proof}
 En appliquant le
lemme \ref{Lt^k} sur les arbres de $Rp(n)$, on a 
\bml\mid\mid\rov{S}[V(T),\Lambda]\mid(t)\mid\leq \sum_{T\in Rp(n)}\frac{1}{\mid n\mid^{r}}(Ct)^k\frac{\mid u_0\mid_r}{\mid u_0\mid_0}\\\times\sum_{\substack{j(V(T))\in
\{0,1,2,1+\frac{1}{2},2+\frac{1}{2},r\}(r)\\
\sum_{i=1}^{\parallel \al_v\parallel}j_i\leq q_1+\frac{3}{2}+r }
}\left(\prod_{v\in 
V(T)}\left(\prod_{i=1}^{\parallel
\al_v\parallel}\mid\frac{\epsilon\mid\al_v^i\mid^{j_i}}{\mid u_0\mid_{j_i}
}\f{u}_{\al_v^i}\mid\right)\right).\end{multline*}
Si l'indice $j_i^v$ qui vaut $r$ est tel que $\al_v^i$ est plus grand que
$n$, alors on pose $r=0$ dans la formule et on multiplie par $\frac{1}{\mid
n\mid^r}$, en constatant que pour $j_l^u$ tel que $\al_u^l=n$, on a
$\mid n\mid^r\leq \mid \al_v^i\mid^{j_l^u}\mid \al_u^l\mid^{r-j_l^u}$, ce qui
donne bien la r\'eponse.

Si l'indice $j_i^v$ qui vaut $r$ est tel que $\al_v^i$ est plus petit que
$n$, pour $j_l^u$ tel que $\al_u^l=n$:
$\mid \al_v^i\mid^{r}\leq \mid \al_v^i\mid^{j_l^u}\mid
\al_u^l\mid^{r-j_l^u}$,
qui conclu la preuve du th\'eor\`eme.

\end{proof}
\section{D{\'e}monstration du th{\'e}or{\`e}me d'existence locale}
\begin{theoreme}\label{Treffaible} 
 Soit un polyn{\^o}me $P\in R[X,Y]$ de degr{\'e} $p$ et un r{\'e}el $\lambda$. Soit $q$ le
  degr{\'e} en la variable $Y$ de $P$ et
  $s=\inf(q,2)$, $r\geq s$. Pour tout r{\'e}el $\epsilon$,   $0<\epsilon<1$ , et
  en notant  \[t_0=\left(\frac{M(P,\Lambda)}{\epsilon^q(\mid u_0\mid_s)^{2q}(\mid u_0\mid_0)^{2p}}\right),\]
  il existe  $\Omega(u_0,t)$ une fonction de $[0,t_0]$ dans $\R$ telle que
l'application $\phi$, donn{\'e}e par le lemme de formalisation 
\ref{lemformel} converge analytiquement pour tout $t$ de $[0,t_0]$ et  $\phi\in 
{\cal A}[W^{r}_1(\T)]$. De plus
pour tout entier relatif $n$, on a
\[\f{\phi}_n=\exp\left(i\int_0^t\Omega_n(u_0,y)dy\right)\f{\chi_n},\] o{\`u}
$\Omega$ est  donn{\'e}e par $\chi\in W_1^{r}(\T)$ de la fa\c con suivante : 
\[\Omega_n(u_0,t)=\sum_{T\in (Rp(n)^1)} 
f_v(T)\prod_{i\neq j,\al^j(T)=\sum_i\al^i(T)}\prod_i
\f{\chi}_{\al_v^i(T)}(t,u_0)),\] et est deux fois continuement d\'erivable sur $[0,t_0]$.

\end{theoreme} 
\begin{proof}
\[\sum_{(T,F(T))\in (T_{\al,\be}^e,F_\al(n))}\rov{\cal S}[V(T)]=\sum_{(T,F(T))\in (T_{\al,\be},F_\al(n))}\sum_{T'\in S(T)}\rov{\cal S}[V(T')].\]
(o{\`u} $T_{\al,\be}$ sont les arbres non {\'e}tiquet{\'e}s d'apr{\`e}s les d{\'e}finitions \pageref{Dfruit}).
D'apr{\`e}s le th{\'e}or{\`e}me de convergence \ref{TCgen}, on a
\bml\mid\f{\phi}_n\mid\leq \sum_{\al,\be}\sum_{T\in T_{\al,\be},F(T)\in F_\al(n)}\mid \epsilon^\alpha
 (\frac{\mid u_0\mid_r}{\mid u_0\mid_0}\frac{1}{(\mid
 n\mid)^{r}}\\\left(\sum_{\substack{j(T)\in
  \{0,1,s,r\}(r)\\
  \forall v :\sum_{i=1}^{\parallel \al_v\parallel}j_i\leq q_1+r-1 }
}\prod_{v\in 
V(T)}\left(\prod_{i=1}^{\parallel
  \al_v\parallel}\mid\frac{\epsilon\mid\al_v^i\mid^{j_i^v}}{\mid
  u_0\mid_{j_i^v} }\f{u}_{\al_v^i}\mid\right)\right).\end{multline*}
Par d{\'e}finition, $F(T)\in F(n)$ implique $\sum_v\sum_i\al_v^i=n$. Comme il
existe au plus $5^\be$ arbres distincts de $\be$ noeuds, on a alors
\bml\mid\f{\phi}_n\mid\leq \sum_\be 5^\be\sum_{\al}\sum_{\sum_{i=1}^\al\sum_{e_v=1}^\be\al_v^i=n}\mid \epsilon^\alpha
 (\frac{\mid u_0\mid_r}{\mid u_0\mid_0}\mid\frac{1}{\mid
 n\mid)^{r}}\\\left(\sum_{\substack{j(V(T))\in
  \{0,1,s,r\}(r)\\
  \forall v :\sum_{i=1}^{\parallel \al_v\parallel}j_i^v\leq q_1-1+r }
}\prod_{v\in 
V(T)}\left(\prod_{i=1}^{\parallel
  \al_v\parallel}\mid\frac{\epsilon\mid\al_v^i\mid^{j_i^v}}{\mid
  u_0\mid_{j_i^v} }\f{u}_{\al_v^i}\mid\right)\right).\end{multline*}
En regroupant tous les   $\al_v$ possibles, il vient
\bs\bml\mid\f{\phi}_n\mid\leq \sum_\be 5^\be\sum_{\al}\mid \epsilon^\alpha
 (\frac{\mid u_0\mid_r}{\mid u_0\mid_0}\mid\frac{1}{\mid n\mid^{r}}\\\prod_{v}\sum_{\sum_{i=1}^\al\sum_{v=1}^\be\al_v^i=n}\left(\sum_{\substack{(j_1,...j_{\parallel \al_v\parallel})\in
  \{0,1,s,r\}^{\parallel \al_v\parallel}\\
  \sum_{i=1}^{\parallel \al_v\parallel}j_i\leq q_1-1+r }
}\left(\prod_{i=1}^{\parallel
  \al_v\parallel}\mid\frac{\epsilon\mid\al_v^i\mid^{j_i}}{\mid
  u_0\mid_{j_i} }\f{u}_{\al_v^i}\mid\right)\right).\end{multline*}\es
 Il existe donc un r\'eel 
  $K(q)$ tel que
\[\sum_{n\in \Z}\mid n^r\phi_n\mid\leq (\sum_\be (K(q)5)^\be\sum_{\al>\be}\epsilon^\al)\mid u_0\mid_r).\]
D'o\`u le th\'eor\`eme.

\end{proof}

\begin{theoreme}\label{Texist}
 Soit $q$ le
  degr{\'e} en la variable $Y$ de $P$ et $s=\inf(q,2)$, $r\geq s$. Pour tous r{\'e}els positifs $\epsilon<1$ et  $\epsilon_0<1$ et en notant   \[t_0=\left(\frac{M(P,\Lambda)}{\epsilon_0\epsilon^q(\mid u_0\mid_0)^{2p}(\mid u_0\mid_s)^{2(q)}}\right),\]  il existe une  unique solution  $u(t)\in W_1^{r}$ pour tout $t\in
[0,t_0]$ de l'{\'e}quation perturb{\'e}e: 
\[\bc
\partial_tu=\partial_x^3u+P(u,\partial_xu),\\
u((0,x)=u_0\in W_1^{r}(\T).
\ec\]
\end{theoreme}
\begin{proof}
D'apr{\`e}s le th{\'e}or{\`e}me \ref{Treffaible}, l'application $\phi$ d\'efinit
  dans le th\'eor\`eme \ref{Treffaible} est
  uniform{\'e}ment convergente sur $[0,t_0]$. On peut d{\'e}river sous le signe
  sommes et le lemme formel \ref{lemformel}  montre que
$\phi(v(t),t)\in  W_1^{r}(\T)$  est solution de
l'{\'e}quation perturb{\'e}e (si $v(t)$ est solution de l'{\'e}quation
lin{\'e}aris{\'e}e \ref{Eqexlin} de condition initiale $u_0\in W_1^{r}(\T)$).

 De plus, le th{\'e}or{\`e}me \ref{Treffaible} montre que $\phi$ est
inversible et diff{\'e}rentiable (convergence analytique) sur
$W_1^{r}(\T)$. Soit $u$ une solution de l'{\'e}quation perturb{\'e}e. On note
$v=\phi^{-1}(t,u)$, o{\`u} pour tout $z$ de $W_1^{r}(\T)$,  tout
$t\in [0,t_0]$, on a que $\phi(t,\phi^{-1}(t,z))=\phi^{-1}(t,\phi(t,z))=z$. On en d{\'e}duit,
alors en d{\'e}rivant suivant $t$ que
\[\partial_t\phi^{-1}(t,\phi(t,z))+D_{(t,\phi(t,z))}\phi^{-1}.\partial_t\phi(t,x).\]
Par construction de $\phi$, on a:
\[\partial_t\phi(t,v)+D_{(t,v)}\phi.\partial_x^3v=\partial_x^3\phi(t,v)+P(\phi(t,v),\partial_x\phi(t,v)).\]
Ce qui nous donne
alors\[\partial_x^3v=D_{(t,u)}\phi^{-1}\partial_tu-D_{(t,u}\phi^{-1}\partial_t\phi(t,v).\]
Or on remarque que
\[\partial_tv=\frac{d\phi(t,u(t))}{dt}=D_{(t,u)}\phi^{-1}\partial_tu+\partial_t\phi^{-1}(t,u),\]
qui implique
\[\partial_tv-\partial_x^3v=\partial_t\phi^{-1}(t,u)+D_{(t,u)}\phi^{-1}\partial_t\phi(t,v)=0.\]
Donc si $u$ est solution de l'{\'e}quation nonlin{\'e}aire, alors $v$ est
l'unique 
solution de l'{\'e}quation lin{\'e}aire avec pour condition initiale
$\phi^{-1}(0,u_0)=u_0$ (par construction). Donc $u$ est unique.
\end{proof}
L'unicit{\'e} utilise une propri{\'e}t{\'e} de $\phi$ qui est semblable {\`a}
celle 
d'une forme normale param{\'e}tr{\'e}e suivant le temps.

On d\'emontre aussi dans $H^k(\T)$ le m\^eme type de r\'esultat
\begin{theoreme}\label{TexitsurH} 
 Soit un polyn{\^o}me $P\in R[X,Y]$ de degr{\'e} $p$ et un r{\'e}el $\lambda$. Soit $q$ le
  degr{\'e} en la variable $Y$ de $P$ et
  $s=\inf(q+\frac{1}{2},\frac{5}{2})$, $r\geq s$. Pour tous r{\'e}el positif  $\epsilon<1$, en notant  \[t_0=\left(\frac{M(P,\Lambda)}{\epsilon^q(\mid u_0\mid_s)^{2q}(\mid u_0\mid_0)^{2p}}\right),\]
  il existe  $\Omega(u_0,t)$ une fonction de $[0,t_0]$ dans $\R$ telle que
l'application $\phi$, donn{\'e}e par le lemme de formalisation 
\ref{lemformel} converge analytiquement pour tout $t$ de $[0,t_0]$ et  $\phi\in 
{\cal A}[H^{r+\frac{1}{2}}(\T)]$. De plus
pour tout entier relatif $n$ on a
\[\f{\phi}_n=\exp\left(i\int_0^t\Omega_n(u_0,y)dy\right)\f{\chi_n},\] o{\`u}
$\Omega$ est  donn{\'e}e par $\chi\in H^{r+\frac{1}{2}}(\T)$ de la fa\c con
suivante : 
\[\Omega_n(u_0,t)=\sum_{T\in (Rp(n)^1)} 
f_v(T)\prod_{i\neq j,\al^j(T)=\sum_i\al^i(T)}\prod_i
\f{\chi}_{\al_v^i(T)}(t,u_0)),\] et est deux fois continuement d\'erivables sur $[0,t_0]$.

\end{theoreme}
\begin{proof}
Si on note $Rp(n)$ l'ensemble des arbres fruitiers dont le fruit est dans $F(n)$
et tel qu'il existe un fruit de $T$ qui soit \'egal \`a $n$.
On applique le  theoreme \ref{TnormparR(n)} sur les arbre qui ne
sont pas dans $Rp(n)$, et le th\'eor\`eme \ref{Tnormpar} sur les arbres de
$Rp(n)$. Alors  comme dans la d\'emonstration du th\'eor\`eme \ref{Treffaible}, on a:
\[\sum_n\mid n^{2r}\phi_n\mid\leq K(\parallel u_0\parallel_r+\mid u_0\mid_{r-1}).\]
\end{proof}
On en d\'eduit aussi 
\begin{theoreme}
 Soit $q$ le
  degr{\'e} en la variable $Y$ de $P$ et $s=\inf(q+\frac{1}{2},\frac{5}{2})$, $r\geq s$. Pour tous r{\'e}els positifs $\epsilon<1$ et  $\epsilon_0<1$ et en notant   \[t_0=\left(\frac{M(P,\Lambda)}{\epsilon_0\epsilon^q(\mid u_0\mid_0)^{2p}(\mid u_0\mid_s)^{2(q)}}\right),\] alors il existe une  unique solution  $u(t)\in H^{r+\frac{1}{2}}$ pour tout $t\in
[0,t_0]$ de l'{\'e}quation perturb{\'e}e: 
\[\bc
\partial_tu=\partial_x^3u+P(u,\partial_xu),\\
u((0,x)=u_0\in H^{r+\frac{1}{2}}(\T).
\ec\]\end{theoreme}
Ce qui   d\'emontre bien le th\'eor\`eme \ref{TfnKdVS}.

Dans le cas plus g\'en\'eral, pour $\lambda_i$ des fonctions
 localement $C^1$, on note
$\lambda=(\lambda_1,.....\lambda_m)\in (C^1)^m$. Soit une fonction $F$ de
$(C^1)^m\times [0,t_0]$, differentiable ( au sens de Frechet), on note
$D_\lambda F$ sa differentielle sur $(C^1)^m$ et 
\[\parallel D_\lambda F\parallel=\sum_i\mid
\partial_{\lambda_i}F(\lambda,t)\mid.\]  
On g\'en\'eralise le th\'eor\`eme \ref{TexitsurH} comme suit. 

\begin{theoreme}\label{TexistGen}
Soit $2k+1\leq n$ et $r\geq n$. Si $j$ est le plus grand indice tel que  $\lambda_{2j}(t)$ soit
non globalement 
nul, et si   $\lambda_{2j}(0)$ est strictement n\'egatif ( ou bien $j=0$)
alors il existe une unique solution local $u(\lambda,t)$ de l'\'equation
\[\bc
\partial_tu=\sum_j\lambda_j(t)\partial_x^ju+Q(u),\\
u((0,x)=u_0\in H^{r+\frac{1}{2}}(\T).
\ec\]
pour $t\in [0,t_0]$ o\`u
\[t_0\geq (\sup_{t\in[0,t_0]}(\sum_i\mid\lambda_i(t)\mid))T(P,\parallel
u_0\parallel_k),\]
($T$ est ind\'ependant des $\lambda_i$).
De plus, si on pose $F(\lambda,t)=\Omega(u(\lambda,t))$, on obtient pour $t\in [0,t_0]$,
\[\parallel \partial_\lambda F(\lambda,t)\parallel\leq t_0A(\parallel
u_0\parallel_k)q,\]
et \[\parallel F(\Lambda,t)\parallel\leq\sum_i\mid \lambda_i(t)\mid.\]
\end{theoreme}
\begin{proof}
Il existe $t_1$ tel que
$\sup_{t\in[0,t_1]}(\sum_i\mid\lambda_i(t)\mid+\sum_i\mid\partial_t\lambda_i(t)\mid)
$ existe. On proc\`ede comme dans la d\'emonstration du th\'eor\`eme \ref{Treffaible}, en
appliquant le th\'eor\`eme de convergence \ref{Tconvgen1}. Il existe alors bien  $t_0$
v\'erifiant les hypoth\`eses de l'\'enonc\'e du th\'eor\`eme \ref{TexistGen},
tel que $\Omega(\Lambda,t)$ et $\partial_t\Omega(\Lambda,t)$ existe pour tout $t\leq t_0$. De plus,
$\Omega(\Lambda,t)$ est bien major\'e par une constante ne d\'ependant pas de
$P$, $\lambda$ et $u_0$, qu'on peut choisir tel que
\[(\sup_{t\in[0,t_1]}(\sum_i\mid\Omega_i(t)\mid+\sum_i\mid\partial_t\Omega_i(t)\mid)\leq
(\sup_{t\in[0,t_1]}(\sum_i\mid\lambda_i(t)\mid+\sum_i\mid\partial_t\lambda_i(t)\mid).\]
De m\^eme on montre  que
\[\parallel F(\Lambda,t)\parallel\leq\sum_i\mid \lambda_i\mid.\]

\end{proof}
Le th\'eor\`eme du point fixe implique  que
\begin{theoreme}Soit $2k+1\leq n$ et $r\leq n$. Si $j$ est le plus grand indice tel que  $\Omega_{2j}(u_0)$ soit
non globalement 
nul, et si   $\lambda_{2j}(u_0)$ est strictement n\'egatif ( ou bien $j=0$)
alors il existe une unique solution local $u(t)$ de l'\'equation
\[\bc
\partial_tu=\sum_j\lambda_j\partial_x^ju+P(u),\\
u((0,x)=u_0\in H^{r+\frac{1}{2}}(\T).
\ec\]
\end{theoreme}
Et le th\'eor\`eme r\'eciproque suivant.
\begin{theoreme} Soit $2k+1\leq n$ et $r\leq n$. Si $j$ est le plus grand indice tel que  $\lambda_{2j}(t)$ soit
non globalement 
nul, alors si   $\lambda_{2j}(0)$ est strictement possitif 
alors s'il existe une unique solution local $u(t)$ de l'\'equation
\[\bc
\partial_tu=\sum_j\Omega_j(u(t))\partial_x^ju+Q(u),\\
u((0,x)=u_0\in H^{r+\frac{1}{2}}(\T).
\ec\]
Alors il existe une solution $v(t)$ de l'\'equation:
\[\bc
\partial_tv=\sum_j\Omega_j(u(t))\partial_x^jv,\\
v((0,x)=u_0\in H^{r+\frac{1}{2}}(\T).
\ec\]
et pour tout $\epsilon>0$ il existe $t_0>0$ tel que pour tout $0<t<t_0$ et tout
$k>r$, on a
\[\lim_{n\to \infty}\mid n^ke^{\int_0^t\sum_j\Omega_j(u(t))(in)^j}\f{v}_n\mid\leq \epsilon.\]

\end{theoreme}
\begin{proof} On applique le th\'eor\`eme \ref{Tconvgen2} sur \[\sum_{T\in
\cal{T}, F(T)\in F(n)}\exp(\sum_i\int_0^t\Omega_i(n,s))\lov{S}[V(T),\lambda].\]
\end{proof}
On en d\'eduit finalement :

\begin{theoreme} Si $j$ est le plus grand indice tel que  $\lambda_{2j}(t)$ soit
non globalement 
nul, et si   $\lambda_{2j}(0)$ est strictement possitif, 
alors  l'\'equation:
\[\bc
\partial_tu=\sum_j\Omega_j(u(t))\partial_x^ju+Q(u)\\
u((0,x)=u_0\in H^{r+\frac{1}{2}}(\T);
\ec\]
est mal pos\'ee.
\end{theoreme}
Ce qui d\'emontre  le th\'eor\`eme \ref{Tprincgeneral}.
\section{Contr{\^o}le de $\mid u(t)\mid_{s}$ et $s\leq \frac{5}{2}$}
En vue de d{\'e}montrer l'existence globale des solutions de KdV perturb{\'e},
nous allons d{\'e}montrer une in{\'e}galit{\'e} sur  $\mid
u(t)\mid_{\frac{5}{2}}$ ou $\mid
u(t)\mid_{2}$. On commence  sur un cas particulier, la d{\'e}monstration g{\'e}n{\'e}rale  {\'e}tant identique.
\begin{theoreme}
Soit $u$, definie sur $[0,t_1]$ ( $t_1\geq t_0$), une solution de
\[\bc
\partial_tu=\partial_x^3u+u^2\partial_xu\\
u(0)=u_0\in \cal{H}^{\frac{5}{2}}(\T),
\ec\]
 alors il existe une constante $M$ tel que
pour tout  $t\in [0, t_1]$, et tout $k$, \bs\[\mid u(t)\mid_k\leq \mid u_0\mid_k\left(M\mid u(0)\mid_0^{2}+1\right)+     (M\sup_{t'\in [0,t_1]}\left(\parallel u(t')\parallel_{1}^2\mid u(t')\mid_k\right)\mid t\mid.\]\es
\end{theoreme}

\begin{proof}
Nous avons que:
\bs\[\frac{d\f{u}_n}{dt}=-in^3\f{u}_n+in\left(\sum_{p\in \Z^*}\f{u}_p\f{u}_{-p}\right)\f{u}_n+\sum_{\substack{n_1+n_2+n_3=n,\\(n_1+n_2)(n_1+n_3)(n_2+n_3)\neq 0}}in\f{u}_{n_1}\f{u}_{n_2}\f{u}_{n_3}.\]\es
Posons \[\f{w}_n(t)=\f{u}_n(t)\exp\left(in^3t-in\int_0^t\left(\sum_{p\in \Z^*}\f{u}_p(s)\f{u}_{-p}(s)ds\right)\right),\]
 on a alors:
\bs\[\frac{d\f{w}_n(t)}{dt}=\sum_{\substack{n_1+n_2+n_3=n,\\(n_1+n_2)(n_1+n_3)(n_2+n_3)\neq 0}}in\f{w}_{n_1}(t)\f{w}_{n_2}(t)\f{w}_{n_3}(t)e^{i(n^3-(n_1)^3-(n_2)^3-(n_3)^3)t}.\]\es
Dont on d\'eduit, en int\'egrant que
\bml\left\vert n^{\frac{3}{2}}\f{w}_n(t)- n^{\frac{3}{2}}\f{w}_n(0)\right\vert\leq \mid n^{\frac{3}{2}+1}\mid\times\\\left\vert \int_0^t\sum_{\substack{n_1+n_2+n_3=n,\\(n_1+n_2)(n_1+n_3)(n_2+n_3)\neq 0}}\f{w}_{n_1}(s)\f{w}_{n_2}(s)\f{w}_{n_3}(s)e^{i(n^3-(n_1)^3-(n_2)^3-(n_3)^3)s}ds\right\vert.\end{multline*}
Puisuqe 
\[(n^3-(n_1)^3-(n_2)^3-(n_3)^3)=3(n_1+n_2)(n_2+n_3)(n_1+n_3),\]
 en supposant $\mid n_1\mid>\mid n_2\mid>\mid n_3\mid$, on a  soit $\mid (n^3-(n_1)^3-(n_2)^3-(n_3)^3)\mid>\frac{\mid n_1\mid^2}{8}$, soit $\mid (n^3-(n_1)^3-(n_2)^3-(n_3)^3)\mid\leq\frac{\mid n_1\mid^2}{8}$.
Dans le premier cas nous dirons que $(n_1,n_2,n_3)$ appartiennent {\`a} ${\cal N}_1$ et dans le deuxi{\`e}me cas que $(n_1,n_2,n_3)$ appartiennent {\`a} ${\cal N}_2$. Notons que $\mid n_1\mid\geq \frac{\mid n\mid}{3}$.

Dans le premier cas on int{\`e}gre par partie :
\begin{eqnarray*}
\f{w}_{n_1}(s)\f{w}_{n_2}(s)\f{w}_{n_3}(s)&\to&\partial_s\left(\f{w}_{n_1}(s)\f{w}_{n_2}(s)\f{w}_{n_3}(s)\right)\\
e^{i(n^3-(n_1)^3-(n_2)^3-(n_3)^3)s}&\to&\frac{e^{i(n^3-(n_1)^3-(n_2)^3-(n_3)^3)s}-e^{(i(n^3-(n_1)^3-(n_2)^3-(n_3)^3)t}}{i(n^3-(n_1)^3-(n_2)^3-(n_3)^3)},
\end{eqnarray*}
en remarquant que \bml\partial_s\left(\f{w}_{n_1}(s)\f{w}_{n_2}(s)\f{w}_{n_3}(s)\right)=(\sum_{j=1}^3\frac{\f{w}_{n_1}(s)\f{w}_{n_2}(s)\f{w}_{n_3}(s)}{\f{w}_{n_j}(s)}in_j\times\\\sum_{\substack{n_4+n_5+n_6=n,\\(n_4+n_5)(n_4+n_6)(n_5+n_6)\neq 0}}\f{w}_{n_4}(s)\f{w}_{n_5}(s)\f{w}_{n_6}(s)e^{i((n_j)^3-(n_4)^3-(n_5)^3-(n_6)^3)s}).\end{multline*}
Nous avons donc
\bs\bml\left\vert \int_0^t\sum_{(n_1,n_2,n_3)\in {\cal N}_1}\f{w}_{n_1}(s)\f{w}_{n_2}(s)\f{w}_{n_3}(s)e^{i(n^3-(n_1)^3-(n_2)^3-(n_3)^3)s}ds\right\vert\leq\\\sum_{\substack{n_1+n_2+n_3=n,\\(n_1+n_2)(n_1+n_3)(n_2+n_3)\neq 0}}32\left\vert \frac{\f{w}_{n_1}(0)\f{w}_{n_2}(0)\f{w}_{n_3}(0)}{n^2}\right\vert+\\\sum_{\substack{n_1+n_2+n_3=n,\\(n_1+n_2)(n_1+n_3)(n_2+n_3)\neq 0}}\frac{8}{(n_1)^2}\int_0^t\mid(\sum_{j=1}^3\frac{\f{w}_{n_1}(s)\f{w}_{n_2}(s)\f{w}_{n_3}(s)}{\f{w}_{n_j}(s)}in_j\times\\\sum_{\substack{n_4+n_5+n_6=n,\\(n_4+n_5)(n_4+n_6)(n_5+n_6)\neq 0}}\f{w}_{n_4}(s)\f{w}_{n_5}(s)\f{w}_{n_6}(s)e^{i((n_j)^3-(n_4)^3-(n_5)^3-(n_6)^3)s})\mid ds.\end{multline*}\es
Il suffit de remarquer que 
\bml\mid n\mid^{\frac{3}{2}+1}\sum_{\substack{n_1+n_2+n_3=n,\\(n_1+n_2)(n_1+n_3)(n_2+n_3)\neq 0}}\left\vert \frac{\f{w}_{n_1}(0)\f{w}_{n_2}(0)\f{w}_{n_3}(0)}{n^2}\right\vert\leq \\\sum_{\substack{n_1+n_2+n_3=n,\\(n_1+n_2)(n_1+n_3)(n_2+n_3)\neq 0}}\left\vert(n_1)^{\frac{3}{2}}\f{w}_{n_1}(0)\f{w}_{n_2}(0)\f{w}_{n_3}(0)\right\vert\leq \mid u(0)\mid_{\frac{3}{2}}(\mid u(0)\mid_0)^2,\end{multline*}
 et aussi que 
\bs\bml\mid n\mid^{\frac{3}{2}+1}\sum_{\substack{n_1+n_2+n_3=n,\\(n_1+n_2)(n_1+n_3)(n_2+n_3)\neq 0}}\frac{8}{(n_1)^2}\int_0^t\mid(\sum_{j=1}^3\frac{\f{w}_{n_1}(s)\f{w}_{n_2}(s)\f{w}_{n_3}(s)}{\f{w}_{n_j}(s)}in_j\times\\\sum_{\substack{n_4+n_5+n_6=n,\\(n_4+n_5)(n_4+n_6)(n_5+n_6)\neq 0}}\f{w}_{n_4}(s)\f{w}_{n_5}(s)\f{w}_{n_6}(s)e^{i((n_j)^3-(n_4)^3-(n_5)^3-(n_6)^3)s})\mid ds\leq \\\sum_{n_1+n_2+n_3+n_4+n_5=n}t\sup_{t'\in [0,t]}\left(\left\vert (n_1)^{\frac{3}{2}}\f{w}_{n_1}(t')\f{w}_{n_2}(t')\f{w}_{n_3}(t')\f{w}_{n_4}(t')\f{w}_{n_5}(t')\right\vert\right),\end{multline*}\es
pour avoir l'in{\'e}galit{\'e}:
\bs\bml\mid n^{\frac{3}{2}+1}\mid\times\left\vert \int_0^t\sum_{(n_1,n_2,n_3)\in {\cal N}_1}\f{w}_{n_1}(s)\f{w}_{n_2}(s)\f{w}_{n_3}(s)e^{i(n^3-(n_1)^3-(n_2)^3-(n_3)^3)s}ds\right\vert\leq \\32 \mid u(0)\mid_{\frac{3}{2}}(\mid u(0)\mid_0)^2+t 16\sup_{t'\in [0,t]}(\mid u(t')\mid_{\frac{3}{2}}(\mid u(t')\mid_0)^4).\end{multline*}\es
Dans le cas o{\`u} $(n_1,n_2,n_3)$ sont dans ${\cal N}_2$, on a  $\mid
(n^3-(n_1)^3-(n_2)^3-(n_3)^3)\mid\leq\frac{\mid n_1\mid^2}{8}$, donc soit $\mid
n_2\mid>\frac{\mid n_1\mid }{2}$ et $\mid n_3\mid>\frac{\mid n_1\mid }{2}$, soit
$\mid n_2\mid>\frac{\mid n_1\mid }{2}$ et $\mid n_3\mid>\frac{\mid n_2\mid
}{2}$. On peut donc supposer que $\mid n_2\mid>\mid n_3\mid>\frac{\mid n_1\mid
}{4}$. On a donc:
\bml\mid n^{\frac{3}{2}+1}\mid\left\vert \int_0^t\sum_{(n_1,n_2,n_3)\in {\cal N}_2}\f{w}_{n_1}(s)\f{w}_{n_2}(s)\f{w}_{n_3}(s)e^{i(n^3-(n_1)^3-(n_2)^3-(n_3)^3)s}ds\right\vert\leq \\\sum_{\substack{n_1+n_2+n_3=n,\\\mid n_1\mid>\mid n_2\mid>\mid n_3\mid>\frac{\mid n_1\mid }{4}}}\int_0^t\left\vert  (n_1)^{\frac{3}{2}+1}\f{w}_{n_1}(s)\f{w}_{n_2}(s)\f{w}_{n_3}(s)\right\vert ds.\end{multline*} On conclut en remarquant finalement que 
\[\sum_{n_2=\frac{\mid n_1\mid }{4}}^{\mid n_1\mid}\mid \f{w}_{n_2}\mid\leq \parallel w\parallel_1\sqrt{\sum_{n_2=\frac{\mid n_1\mid }{4}}^{\mid n_1\mid}\frac{1}{n_2^2}},\] qui implique que
\[\sum_{n_2=\frac{\mid n_1\mid }{4}}^{\mid n_1\mid}\mid \f{w}_{n_2}\mid\leq \parallel w\parallel_1\sqrt{\frac{8}{\mid n_1\mid}},\] et qui donne:
\bs\[\sum_{\substack{n_1+n_2+n_3=n,\\\mid n_1\mid>\mid n_2\mid>\mid n_3\mid>\frac{\mid n_1\mid }{4}}}\int_0^t\left\vert  (n_1)^{\frac{3}{2}+1}\f{w}_{n_1}(s)\f{w}_{n_2}(s)\f{w}_{n_3}(s)\right\vert ds\leq t\sup_{t'\in[0,t]}(\mid u(t')\mid_{\frac{3}{2}}\parallel u(t')\parallel_1^2).\]\es

\end{proof}

\begin{theoreme}\label{Texist2}
Soit $P$ un polyn{\^o}me de degr{\'e} $q$:
\[\bc
\partial_tu=\partial_x^3u+\partial_xP(u)\\
u(0)=u_0\in \cal{H}^{2}(\T).
\ec\]
 S'il existe une solution  $u\in {\cal H}^{\frac{5}{2}}(\T)$ definie sur
$[0,t_1]$ ( $t_1\geq t_0$) alors pour tout  $t\in [0, t_1]$, et tout $k\leq
\frac{5}{2}$, on a \begin{multline*}\mid u(t)\mid_k\leq \mid u_0\mid_k\left(M(P,\lambda)\mid u(0)\mid_0^{q}+1\right)\\+ (M(P,\lambda))^2\sup_{t'\in [0,t_1]}\left(\mid u(t')\mid_0^{q}\parallel u(t')\parallel_{1}^3\mid u(t')\mid_k\right)\mid t\mid.\end{multline*}
\end{theoreme}
\begin{proof}
 Si on note
\[\f{P}_n(u)=\sum_{\al=2}^q\sum_{\un{n}\in\Z^\al(n)}P_{\un{n}}\prod_{j=1}^\al\f{u}_{n_j},\]

on pose   \[\exp{\left(-in^3t-in\int_0^t\sum_{\al>1}\sum_{\un{m}\in\Z^\al(0)}P_{(\un{m},n)}\prod_j\f{u}_{n_j}(s)ds\right)}\f{w}_n(t)=\f{u}_n(t).\] Alors
\bs\[\frac{d\f{w}_n(t)}{dt}=\sum_{\substack{\un{n}\in \ov{Z}^\al(n)\\\forall j n_j\neq n}}inP_{\un{n}}\exp(i(\sum_jn_j^3-n^3)s)\prod_j\f{w}_{n_j}(t)\mid,\]\es
d'o\`u
\[\mid w_n(t)-w_n(0)\mid\leq \mid n\int_0^t \sum_{\al}\sum_{\substack{\un{n}\in \ov{Z}^\al(n)\\\forall j n_j\neq n}}P_{\un{n}}\exp(i(\sum_jn_j^3-n^3)s)\prod_j\f{w}_{n_j}(s)ds\mid.\] 

D'apr{\`e}s la proposition \ref{ARO4}, on distingue deux cas
\begin{description}
\item[(1)]  
 soit il existe $\mid n_1\mid=\sup(\mid n_i\mid)\geq \mid n_2\mid\geq \mid n_3\mid\geq \mid n_4\mid$ tels que  $\mid n_2n_3n_4\mid\geq \frac{\mid n_1\mid ^{2}}{\al^2}$ et $\mid n_2\mid>\frac{\mid n_1\mid^{\frac{2}{3}}}{\al^{\frac{1}{3}}}$
\item[(2)]soit $\mid \sum_jn_j^3-n^3 \mid\geq \frac{\sup(\mid n_i\mid )^2}{\al^3}$. 
\end{description}
Dans le deuxi{\`e}me cas, on int{\`e}gre par partie, 
\begin{eqnarray*}
\prod_j\f{w}_{n_j}(s)&\to&\partial_s\left(\prod_j\f{w}_{n_j}(s)\right)\\
e^{s(\sum_jn_j^3-n^3)}&\to&\frac{e^{s(\sum_jn_j^3-n^3)}-e^{t(\sum_jn_j^3-n^3)}}{i(\sum_jn_j^3-n^3))},
\end{eqnarray*}

ce qui donne 
\bs\bml\left\vert n\int_0^t P_{\un{n}}\exp(i(\sum_jn_j^3-n^3)s)\prod_j\f{w}_{n_j}(s)ds\right\vert\leq \left\vert \frac{\al^3}{n}P_{\un{n}}\prod_j\f{w}_{n_j}(0)\right\vert+\\\int_0^t\sum_j\al^3 \mid P_{\un{n}}\prod_{k\neq j}\f{w}_{n_k}(t)\sum_{\be}\sum_{\un{m}\in \ov{Z}^\be(n_j)}P_{\un{m}}\times\\\times\left(\exp(i(\sum_jn_j^3-n^3)s)-1)\exp(i(\sum_jm_j^3-n_j^3)s\right)\prod_j\f{w}_{m_j}(s)ds\mid.\end{multline*}\es

Dans le premier cas, on observe plus attentivement les diviseurs :
\bs\bml(n_1)^3+(n_2)^3+(n_3)^3-(n_1+n_2+n_3)^3+(n_1+n_2+n_3)^3+\sum_{i=4}^\al (n_i)^3-(\sum_{i=1}^\al n_i)^3=\\3(n_1+n_2)(n_1+ n_3)(n_2+n_3)+\sum_{i=4}^\al (n_i)^3+(n_1+n_2+n_3)^3-(\sum_{i=1}^\al n_i)^3.\end{multline*}\es
Supposons que $\mid n_4\mid<\frac{\mid n_1\mid^{\frac{2}{3}}}{\al^2}$, alors comme
 \bs\bml\sum_{i=4}^\al (n_i)^3+(n_1+n_2+n_3)^3-(\sum_{i=1}^\al
n_i)^3=\\3(n_1+n_2+n_3)^2(\sum_{i=4}^\al n_i)+3(\sum_{i=4}^\al
n_i)^2(n_1+n_2+n_3)+\sum_{i=4}^\al (n_i)^3-(\sum_{i=4}^\al (n_i))^3,
\end{multline*}\es
nous avons que  soit $\mid n_2\mid>\frac{\mid n_1\mid}{2}$ et $\mid n_3\mid>\frac{\mid n_1\mid}{4\al}$, soit $\mid n_2\mid>\frac{\mid n_1\mid}{2}$ et $\mid n_3\mid\leq \frac{\mid n_1\mid}{4\al}$ donc $ \mid n_4\mid>\frac{\mid n_1+n_2\mid}{2\al}$, soit $\mid n_2\mid\leq\frac{\mid n_1\mid}{2}$ donc $ \mid n_4\mid>\frac{\mid n_2+n_3\mid}{2\al}$. Dans ce dernier cas on d{\'e}duit de  $\frac{\mid n_1\mid^{\frac{2}{3}}}{\al^2} >\mid n_4\mid>\frac{\mid n_2+n_3\mid}{2\al}$, que $\mid n_2\mid>\mid n_3\mid>\mid n_2\mid-\frac{\mid n_1\mid^{\frac{2}{3}}}{\al}$. Et du fait que $\mid n_2n_3n_4\mid\geq \frac{\mid n_1\mid ^{2}}{\al^2}$, on a que $\mid n_4\mid>\frac{\mid n_1\mid^2}{\mid n_3\mid^2 4\al}$.
Nous avons donc quatre possibilit{\'e}s:
\begin{description}
\item[(i)]$\mid n_1\mid>\mid n_2\mid>\mid n_3\mid>\mid n_4\mid>\frac{\mid n_1\mid^{\frac{2}{3}}}{\al}$,
\item[(ii)] $\mid n_1\mid>\mid n_2\mid>\mid n_3\mid>\frac{\mid n_1\mid}{4\al}$,
\item[(iii)]$\mid n_1\mid>\mid n_2\mid>\frac{\mid n_1\mid}{2}$ et $\mid n_3\mid>\mid n_4\mid>\mid n_3\mid>\frac{\mid n_1+n_2\mid}{2\al}$,
\item[(iv)] $\mid n_4\mid>\frac{\mid n_1\mid^2}{\mid n_3\mid^2 4\al}$, $\mid n_1\mid>\mid n_2\mid>\frac{\mid n_1\mid^{\frac{2}{3}}}{\al^{\frac{1}{3}}}$ et $\mid n_2\mid>\mid n_3\mid>\mid n_2\mid-\frac{\mid n_1\mid^{\frac{2}{3}}}{\al}$.
\end{description}

Dans le cas (i)
\[\sum_{m>\frac{\mid n_1\mid^{\frac{2}{3}}}{\al}}\mid \f{w}_m\mid\leq \parallel w\parallel_1\sqrt{\sum_{m>\frac{\mid n_1\mid^{\frac{2}{3}}}{\al}}\frac{1}{m^2}}\leq \parallel w\parallel_1\frac{\sqrt{2\al}}{\mid n_1\mid^{\frac{1}{3}}},\]
donc il existe $M$ ind{\'e}pendant de $w$ tel que
\[\sum_{\substack{(n_1,n_2,...n_\al)\\\mid n_1\mid>\mid n_2\mid>\mid n_3\mid>\mid n_4\mid>\frac{\mid n_1\mid^{\frac{2}{3}}}{\al}}}\mid n_1\mid^{\frac{3}{2}+1}\left\vert\prod_j\f{w}_{n_j}\right\vert\leq M\mid w\mid_{\frac{3}{2}}\mid w\mid_0^{\al-4}\parallel w\parallel_1^3.\]
Dans le cas (ii)
\[\sum_{m>\frac{\mid n_1\mid}{4\al}}\mid \f{w}_m\mid\leq \parallel w\parallel_1\sqrt{\sum_{m>\frac{\mid n_1\mid}{4\al}}\frac{1}{m^2}}\leq  \parallel w\parallel_1\frac{\sqrt{2\al}}{\mid n_1\mid^{\frac{1}{2}}}\]
donc il existe $M$ ind{\'e}pendant de $w$ tel que
\[\sum_{\substack{(n_1,n_2,...n_\al)\\ \mid n_1\mid>\mid n_2\mid>\mid n_3\mid>\frac{\mid n_1\mid}{4\al}}}\mid n_1\mid^{\frac{3}{2}+1}\left\vert\prod_j\f{w}_{n_j}\right\vert\leq M\mid w\mid_{\frac{3}{2}}\mid w\mid_0^{\al-3}\parallel w\parallel_1^2.\]
Dans le cas (iii)
\[\sum_{m>\frac{\mid n_1+n_2\mid}{2\al}}\mid \f{w}_m\mid\leq \parallel w\parallel_1\sqrt{\sum_{m>\frac{\mid n_1+n_2\mid}{2\al}}\frac{1}{m^2}}\leq  \parallel w\parallel_1\frac{\sqrt{2\al}}{\left(\mid n_1\mid-\mid n_2\mid\right)^{\frac{1}{2}}}\]
et
\[\sum_{\mid n_1\mid> m>\frac{\mid n_1\mid}{2}}\f{w}_m\frac{1}{\mid n_1\mid-\mid m\mid}\leq \parallel w\parallel_1 \frac{4}{\mid n_1\mid}.\]
Donc il existe $M$ ind{\'e}pendant de $w$ tel que
\[\sum_{\substack{(n_1,n_2,...n_\al)\\ \mid n_1\mid>\mid n_2\mid>\frac{\mid n_1\mid}{2}\\ \mid n_3\mid>\mid n_4\mid>\mid n_3\mid>\frac{\mid n_1+n_2\mid}{2\al}}}\mid n_1\mid^{\frac{3}{2}+1}\left\vert\prod_j\f{w}_{n_j}\right\vert\leq M\mid w\mid_{\frac{3}{2}}\mid w\mid_0^{\al-4}\parallel w\parallel_1^3.\]
De m{\^e}me dans le cas (iv), on a 
\[\sum_{m>\frac{(n_1)^2}{(n_3)^2}}\mid \f{w}_m\mid\leq 2 \parallel w\parallel_1\left\vert\frac{n_2}{n_1}\right\vert,\]
 
\[\sum_{\mid n_2\mid >m>\mid n_2\mid-\frac{\mid n_1\mid^{\frac{2}{3}}}{\al}}\mid m\f{w}_m\mid\leq 2 \parallel w\parallel_1\sqrt{\frac{\mid n_1\mid^{\frac{2}{3}} }{\al}},\] et aussi
\[\sum_{m>\frac{\mid n_1\mid^{\frac{2}{3}}}{\al^{\frac{1}{3}}}}\mid \f{w}_m\mid\leq 2 \parallel w\parallel_1\frac{\al}{\mid n_1\mid^{\frac{1}{3}}}.\]
Donc il existe $M$ ind{\'e}pendant de $w$ tel que

\[\sum_{\substack{(n_1,n_2,...n_\al)\\ \mid n_4\mid>\frac{\mid n_1\mid^2}{\mid n_3\mid^2 4\al},\\ \mid n_1\mid>\mid n_2\mid>\frac{\mid n_1\mid^{\frac{2}{3}}}{\al^{\frac{1}{3}}}, \\ \mid n_2\mid>\mid n_3\mid>\mid n_2\mid-\frac{\mid n_1\mid^{\frac{2}{3}}}{\al}}}\mid n_1\mid^{\frac{3}{2}+1}\left\vert\prod_j\f{w}_{n_j}\right\vert\leq M\mid w\mid_{\frac{3}{2}}\mid w\mid_0^{\al-4}\parallel w\parallel_1^3.\]
  En regroupant le premier et le deuxi{\`e}me cas, et en faisant la somme sur les $n$, puis en passant au sup, on a donc l'in{\'e}galit{\'e} souhait{\'e}e:
\[\mid n^{k}(w_n(t)-w_n(0))\mid\leq\mid M t\sup_{t'\in[0,t]}\left(\parallel u(t')\parallel_1^3\mid u(t')\mid^{q-3}\mid w_n(t')\mid_k\right) M(P)+\]\[M(P)\mid u(0)\mid^{2\al}\mid w_n(0)\mid_k.\]
\end{proof}
\section{D\'emonstration des lemmes techniques}
\begin{definition}
 Soit un arbre $T$, on d{\'e}finit pour tout noeud $u$ de $T$, $s_u(T)$ par:

Si $u$ est un noeud non r{\'e}sonant: $s_u(T)=(\mid p_u(T)\mid)$.
Si $u$ est un noeud r{\'e}sonant
 $s_u(T)=0$. 
\end{definition}
\begin{lemme}\label{Ltruc1}
Soit $f_v$ d{\'e}finie dans le lemme formel \ref{lemformel}, alors pour
tout arbre $T$ de $\be$ n\oe{}uds, tel que $F(T)\in F(n)$ (la somme
de ses fruits vaut $n$):
\[\mid\prod_{v\in V(T)}\frac{\mid f_v(T)\mid}{\mid
  s_v(T)\mid}\mid\leq\frac{1}{\mid n\mid}\prod_{v\in
  V(T)}M\left(\sum_{\substack{(j_1,...j_{\parallel\al_v(T)\parallel})\in
  \{0,1\}^{\parallel \al_v(T)\parallel}\\ \sum j_i\leq q+1}}\prod_{i}\mid\al_v^i(T)\mid^{j_i}\right)\]
($\parallel.\parallel$ est donn{\'e}e par la d{\'e}finition \ref{Dfruit}, page
  \pageref{Dfruit}).
De plus dans le cas Hamiltonien:
\[\mid\prod_{v\in V(T)}\frac{\mid f_v(T)\mid}{\mid
  s_v(T)\mid}\mid\leq\prod_{v\in
  V(T)}M.\]
 \end{lemme}
\begin{proof}

Dans le cas Hamiltonien, il existe une constante $M$ d\'ependant de $P$ ( la
perturbation), tel que pour tout arbre $T$, $\mid f_v(T)\mid =M\mid p_v(T)\mid$, qui donne bien ma
r\'eponse. 

Dans le
cas g\'en\'eral, il existe une constante $M$ d\'ependant de la
perturbation $P$, tel que pour tout arbre $T$, on a
\[\mid f_v(T)\mid\leq M(\prod_{w\in A_v}\mid p_w(T)\mid)(\prod_j\mid
\al_v^j(T)\mid),\]
ce qui conclut la preuve.
\end{proof}

\begin{definition}
Soit un r\'eel positif $C$, un noeud d'un arbre $T$ est dit dans $(\cal R)_C(T)$, s'il est r\'esonant ( c.\`a.d qu'il
existe un fruit de $\al_u(T_e(u))$ qui a pour valeur $p_u(T)$), si $\mid
p_u(T)\mid=\sup_j(\mid \al_u^j(T_e(u))\mid)$ et si $\mid
p_u(T)\mid\geq \al(u)\sqrt C$.
\end{definition} Pour simplifier les d\'emonstrations, on suppose que
$\Lambda(x)=x^3$ et  on note $\al(u)=\gamma_u+1$.
\begin{lemme}\label{LA}
Soit $T$ un arbre fruitier.
Supposons que pour tout $u$ de $T$:
\[\mid D_u^\Lambda(T)\mid<\frac{\mid p_u(T)\mid \sqrt C}{2}.\]
Alors :
\begin{description}
\item[(1)] soit $\mid p_u(T)\mid<\mid \al(u)\mid\sqrt{C}$, 
\item[(2)]soit il existe trois fruits de $T_e(u)$ tel que:
\[\mid\al_u^1(T_e(u))\al_u^2(T_e(u))\al_u^3(T_e(u))\mid\geq \frac{\mid p_u(T)\mid^3}{(\al(u))^3};\]

\item[(3)] soit $u$ est dans $(\cal R)_C(T)$.
\end{description}
\end{lemme}
\begin{proof}
Soit $w$ le noeud parmis les noeuds adjacent \`a $u$, o\`u  $u$ v\'erifit :
\[\mid D_w^\Lambda(T)\mid=\sup_{v_in A_u\cup\{u\}}(\mid D_v^\Lambda(T)\mid).\]
Alors 
\[\mid D_w^\Lambda(T)\mid\geq \frac{\mid \sum_{v\in
A_u}D_v^\Lambda(T)-D_u^\Lambda(T)\mid}{\al(u)}.\]
Par d\'efinition,
\[\mid \sum_{v\in
A_u}D_v^\Lambda(T)-D_u^\Lambda(T)\mid=\mid
\sum_j\left(\al_u^j(T_e(u))\right)^3-\left(p_u(T)\right)^3\mid.\]
On pose $q_j=\al_u^j(T_e(u))$, et on applique le lemme arithm\'etique \ref{ARO4} qui
implique que :
\begin{description}
\item[(i)]soit $\mid \sum_{v\in
A_u}D_v^\Lambda(T)-D_u^\Lambda(T)\mid\geq \frac{(\sup_j(\mid
\al_u^j(T_e(u))\mid))^2}{2}$; 
\item[(ii)]soit il existe trois fruits de $T_e(u)$ tel que:
\[\mid\al_u^1(T_e(u))\al_u^2(T_e(u))\al_u^3(T_e(u))\mid\geq \frac{\mid (\sup_j(\mid
\al_u^j(T_e(u))\mid)\mid^3}{(\al(u))^3},\]
\item[(iii)] soit il existe un fruit $\al_u^j(T_e(u))$ tel que  $p_u(T)=\al_u^j(T_e(u))$
et $\mid p_u(T)\mid=\sup_j(\mid \al_u^j(T_e(u))\mid)$.
\end{description}

Donc si $(1)$ n'est pas v\'erifi\'es, le  cas $(iii)$ implique que $u$ est dans
$(\cal R)_C$, donc que $(3)$ est v\'erifi\'e. Et le cas $(i)$ implique que 
\[\mid D_w^\Lambda(T)\mid\geq\frac{(\sup_j(\mid
\al_u^j(T_e(u))\mid))^2}{2\al(u)}\geq \frac{\mid p_u(T)\mid}{2\sqrt{C}},\]
ce qui contredit les hypoth\`eses.
Le cas $(ii)$ implique que $(2)$ est v\'erifi\'e.

\end{proof}

\begin{lemme}\label{LB}
Soit un arbre fruitier $T$ tel que tout noeud $u$ de $V(T)$ ne soit pas dans
$(\cal R)_C(T)$ et tel que pour tout $u$ de $T$, on a
\[\mid D_u^\Lambda(T)\mid<\frac{\mid p_u(T)\mid \sqrt C}{2}.\]
Alors, si on note $u_1$ le premier noeud de $T$, on a
\[\left\vert \prod_{u\in V(T)}p_u(T)\right\vert\leq  \left\vert\prod_{u\in
V(T):u>u_1} p_u(T)\right\vert^{\frac{1}{3}} \left(\prod_{u\in V(T)}\mid \prod_i\al_u^i(T)\mid^{\frac{1}{3}}\al(u)^2 \sqrt C \right).\]
\end{lemme}
\begin{proof}
On applique le lemme \ref{LA}, \`a $T$, donc ou bien $\mid p_u(T)\mid\leq
\al(u)\sqrt C$, ou bien $\mid p_u(T)\mid\leq \al(u)\prod_j\mid
\al_u^j(T_e(u))\mid^{\frac{1}{3}}$. Comme les fruits et les $p_u(T)$ sont des
entiers, on a donc
\[\left\vert p_u(T)\right\vert\leq \left\vert \prod_{w\in A_u}\mid p_w(T)\mid^{\frac{1}{3}}\prod_j\mid \al_u^j(T)\mid^{\frac{1}{3}}\al(u)^2
\sqrt C \right\vert;\]
ce qui conclut la preuve.
\end{proof}
On en d\'eduit alors le 
\begin{lemme}\label{LC}
Soit un arbre fruitier $T$ tel que tout noeud $u$ de $V(T)$ ne soit pas dans
$(\cal R)_C(T)$ et tel que pour tout $u$ de $T$:
\[\mid D_u^\Lambda(T)\mid<\frac{\mid p_u(T)\mid \sqrt C}{2}.\]
Alors, si on note $u_1$ le premier noeud de $T$,
\[\left\vert \prod_{u\in V(T)}p_u(T)\right\vert\leq  \prod_{u\in
V(T)}\left\vert \prod_i\mid \al_u^i(T)\mid^{\frac{1}{2}}\al(u)^{3} (\sqrt C)^{\frac{3}{2}} \right\vert.\]
\end{lemme}
\begin{proof}
Pour tout $u$ de $T$, on note  $l(u,T)$  le nombre de noeud reliant la racine de
$T$ \`a $u$ ( le nombre de noeud inferieur \`a $u$ au sens de la racine).
Nous allons montrer par r\'ecurence sur le nombre de noeud de $T$ que l'on a :
\[\left\vert \prod_{u\in V(T)}p_u(T)\right\vert\leq  \prod_{u\in
V(T)}\left\vert \mid\prod_i
\al_u^i(T)\mid^{\sum_{j=1}^{l(u,T)}\frac{1}{3^j}}\al(u)^{2\sum_{j=0}^{l(u,T)}\frac{1}{3^j}}
(\sqrt C)^{\sum_{j=0}^{l(u,T)}\frac{1}{3^j}} \right\vert.\]
Si $T$ n'a qu'un seul noeud, le lemme \ref{LB} donne la r\'eponse.
Supposons la propri\'e t\'e vrai pour tout arbre de moins de $\be$ noeuds.

Soit $T$ un arbre de $\be $ noeud, le lemme \ref{LB} implique que :
\[\left\vert \prod_{u\in V(T)}p_u(T)\right\vert\leq  \prod_{u\in
V(T):u>u_1}\mid p_u(T)\mid^{\frac{1}{3}} \left\vert \prod_{u\in V(T)}\mid
\al_u(T)^{\frac{1}{3}}\al(u)^{2} (\sqrt C)\mid \right\vert.\]

Or \[\prod_{u\in
V(T):u>u_1}\mid p_u(T)\mid=\prod_{w\in A_{u_1}}\prod_{u\in
V(T(w))}\mid p_u(T)\mid.\]
On applique alors l'hypoth\`ese de r\'ecurence sur les $T(w)$, en remarquant que\\
$\al_u(T(w))=\al_u(T)$, $l(u,T(w))=l(u,T)-1$, on termine la r\'ecurence.
\end{proof}
\begin{lemme}\label{LD}
Soit un arbre $T$, $u$ un noeud de $T$ dans $(\cal R)_C(T)$, soit   $w$ le noeud
adjacent inferieur \`a $u$. Si $w$ n'est pas dans  $(\cal R)_C$, sur $T$, alors $w$ n'est pas un noeud de $(\cal R)_C(T')$ avec
$T'=T\setminus u$.
\end{lemme}
\begin{proof}
Si $w$ est un noeud de $(\cal R)_C(T')$ avec
$T'=T\setminus u$, alors il existe un fruit de $u$ sur $T_e(u)$ qui prend pour
valeur $p_w$. Et ce fruit est le plus grand des fruits de $u$ sur $T_e(u)$,
ce qui est impossible.
\end{proof}
Soit $\cal X(T)$ l'ensemble des noeuds de $T$ qui sont dans $(\cal R)_C(T)$.
Soit $\cal Y(T)$ l'ensemble des noeuds de $T$ qui ne sont pas dans $(\cal R)_C(T)$.
Soit $T(\cal X)=T\setminus \cal Y(T)$ et $T(\cal Y)=T\setminus \cal X(T)$
Soit $\cal Z(T)$ l'ensemble des noeuds $u$ de $T$ qui v\'erifit  $\mid
D_u^\Lambda(T)\mid\geq \frac{\mid p_u(T)\mid\sqrt{C}}{2}$.
Soit  $\cal W(T)$  l'ensemble des noeuds $u$ de $T$ qui sont dans $\mid
D_u^\Lambda(T)\mid< \frac{\mid p_u(T)\mid\sqrt{C}}{2}$
Soit $T_{\cal Z}= T\setminus\cal W$ et $T_{\cal W}= T\setminus\cal Z$. On note
$T_{\cal W}(\cal X)$  resp. $T_{\cal W}(\cal Y)$ sont les sous arbres
de $T_{\cal W}$ dont les noeuds sont dans $(\cal R)_C$ resp. ne sont pas dans $(\cal R)_C$.
Alors,
\bs\bml\prod_{u\in V(T)}\frac{\mid p_u(T)\mid }{\sup(t,\mid D_u^\Lambda(T)\mid)}=\left(\prod_{u\in V(T_{\cal W}(\cal X))}\frac{\mid p_u(T_{\cal W}(\cal X))\mid
}{\sup(t,\mid D_u^\Lambda(T_{\cal W}(\cal X))\mid)}\right)\times\\\times\left(\prod_{u\in
V(T_{\cal W}(\cal Y))}\frac{\mid p_u(T_{\cal W}(\cal Y))\mid }{\sup(t,\mid
D_u^\Lambda(T_{\cal W}(\cal Y))\mid)}\right)\left(\prod_{u\in V(T_{\cal Z})}\frac{\mid p_u(T_{\cal Z})\mid }{\sup(t,\mid D_u^\Lambda(T_{\cal Z})\mid)}\right).\end{multline*}\es
\begin{lemme}\label{LE}
Soit un arbre fruitier $T$, sans noeud r\'esonant. Il existe une injection $\theta$ de $V(T_{\cal
W}(\cal X))$ dans $V(T_{\cal Z})$ tel que pour tout $u$ de $V(T_{\cal
W}(\cal X))$:
\[\mid
D_{\theta(u)}^\Lambda(T)\mid\geq \frac{\mid p_u(T)p_{\theta(u)}(T)\mid}{2\al_u(T_{\cal
W}(\cal X))}.\]
\end{lemme}
\begin{proof}
D'apr\`es le lemme \ref{LD}, pour tout noeud $w$ de $V(T_{\cal Z})$ le noeud
adjacent sup\'erieur \`a $w$ sur $T$ est dans $T_{\cal Z}$ et 
$\sup_j(\mid \al_w^j(T_{\cal
W}(\cal X))\mid)=\sup_j(\mid \al_u^j(T)\mid)$. 

Donc en proc\'edant comme dans
le lemme \ref{LA}, on montre que $\mid D_u^\Lambda(T)\mid\geq \frac{\mid p_u(T)\mid^2}{2}$.
\end{proof}

On en d\'eduit alors le
\begin{lemme}\label{lemccv0}
Soit un entier $n$, un arbre $T$ v{\'e}rifiant $F(T)\in F(n)$ de $\al$
 fruits et $\be$ n\oe{}ud. Soit un r{\'e}el
 positif $C\geq 1$. Il existe un fruit $\al_w^j(T)$ tel que pour tout r\'eel
 positif $r$, on ait:
\bs\bml\prod_{u\in V(T)}\mid\frac{s_u(T)}{\sup(C,\mid D^\Lambda_u\mid)}\mid\leq
 \frac{e^{2\al}\mid \al_w^j\mid^{r-1}}{\mid n\mid^r}
\left(\prod_{u\in V(T)}\mid\prod_i
\al_{u}^{i}(T)\mid^{1}(T)\right)\\\left(\prod_{u\in
     V(T): \mid D_u^\Lambda\mid<\frac{\mid p_v\mid\sqrt{C}}{2}}\sqrt{C}\inf(\mid\frac{1}{D^\Lambda_u}\mid,\frac{1}{C})\right) \left(\prod_{u\in
 V(T): \mid D_u^\Lambda\mid\geq\frac{\mid p_v\mid\sqrt{C}}{2}}\mid\frac{1}{\sqrt{C}}\mid\right).\end{multline*}\es
\end{lemme}

\begin{proof}
Soit $\theta$ donn\'e par le lemme \ref{LE}. On note pour un arbre $T$:
\begin{eqnarray*}
A&=& W(T)\cap Y(T)\\
B&=& W(T)\cap X(T)\\
C&=& Z(T)\setminus \theta (B)\\
D&=& Z(T)\cap \theta (B)
\end{eqnarray*}

Alors:
\bs\bml\prod_{u\in V(T)}\mid\frac{s_u(T)}{\sup(C,\mid D^\Lambda_u\mid)}\mid=\prod_{u\in W(T)}\frac{1}{\sup(C,\mid D^\Lambda_u\mid)}\prod_{u\in A}\mid
p_u\mid\\\left(\prod_{u\in B}\mid p_u\mid\prod_{u\in D}\frac{\mid p_u\mid}{\mid
D_u^\Lambda(T)\mid}\right)\prod_{u\in C}\frac{\mid p_u\mid}{\mid
D_u^\Lambda(T)\mid}.\end{multline*}\es
Or d'apr\`es le lemme \ref{LE},
\[(\prod_{u\in B}\mid p_u\mid\prod_{u\in D}\frac{\mid p_u\mid}{\mid
D_u^\Lambda(T)\mid}=\prod_{u\in B}\frac{\mid p_up_{\theta(u)}\mid}{\mid
D_{\theta(u)}^\Lambda(T)\mid}\leq 2\prod_{u\in B}\al(u).\]
Donc
\[(\prod_{u\in B}\mid p_u\mid\prod_{u\in D}\frac{\mid p_u\mid}{\mid
D_u^\Lambda(T)\mid}\leq 2\prod_{u\in B}\al(u)\sqrt C \prod_{u\in D}\frac{1}{\sqrt{C}}.\]
Par d\'efinition de $C$,
\[\prod_{u\in C}\frac{\mid p_u\mid}{\mid
D_u^\Lambda(T)\mid}\leq \prod_{u\in C}\frac{2}{\sqrt{C}}.\]
D'apr\`es le lemme \ref{LC}, si on pose $\{u: u\geq A\}$ l'ensemble compos\'e des noeuds
de $T$  plus grand ou \'egaux aux noeuds de $A$, on a
\[\prod_{u\in A}\mid
p_u\mid\leq \left(\prod_{u\in V(T): u\geq A}C\al(u)^2\prod_i\mid
\al_{u}^{i}(T)\mid^{\frac{1}{2}}(T)\right).\]
Si on pose $\al_w^j(T)$ le fruit le plus grand en valeur absolu de $T$, s'il
existe un noeud $v$ plus petit que $w$ ( $w\geq w$), qui soit dans $A$. On pose
$z$ le plus grand de ces noeuds; alors d'apr\`es le lemme \ref{LA}, appliqu\'e au
sous arbre $T'=T\setminus\{B\cup C\cup D\}$, il existe trois fruits de $T_e'(z)$,
$\al_z^1,\al_z^2,\al_z^3$ tels que 
\[\mid \al_z^1\al_z^2\al_z^3 \mid\geq\frac{\mid \al_w^j(T)\mid^2 }{\al^3},\]
o\`u bien $ \mid \al_w^j(T)\mid\leq \al$. Donc
\[\mid \al_w^j(T)\mid ^r\leq \al^3 \mid \al_w^j(T)\mid
^{r-\frac{1}{2}}\prod_{\substack{u\in V(T),i\\
(u,i)\neq (w,j)}}\mid
\al_{u}^{i}(T)\mid^{\frac{1}{2}}.\]
S'il n'existe pas de noeud $v$ plus petit que $w$, qui soit dans
$A$, alors $w\not \in \{u: u\geq A\}$.
Donc
\bml\mid \al_w^j(T)\mid ^r\left(\prod_{u\in V(T): u\geq A}C\al(u)^2\prod_i\mid
\al_{u}^{i}(T)\mid^{\frac{1}{2}}(T)\right)\leq\\
\mid \al_w^j(T)\mid ^{r-1}\left(\prod_{u\in V(T)}C\al(u)^2\prod_i\mid
\al_{u}^{i}(T)\mid^{1}(T)\right).\end{multline*}
Donc:
\bs\bml\prod_{u\in V(T)}\mid\frac{s_u(T)}{\sup(C,\mid D^\Lambda_u\mid)}\mid\leq\frac{\al^r}{\mid n\mid^r}\left(\prod_{u\in W(T)}\frac{1}{\sup(C,\mid D^\Lambda_u\mid)}\mid\right)(2\prod_{u\in B}\al(u)\sqrt C) \times\\\prod_{u\in D}\frac{1}{\sqrt{C}}\prod_{u\in C}\frac{2}{\sqrt{C}} \mid \al_w^j(T)\mid ^r\left(\prod_{u\in V(T): u\geq A}C\al(u)^2\prod_i\mid
\al_{u}^{i}(T)\mid^{\frac{1}{2}}(T)\right).\end{multline*}\es
D'o\`u le lemme.

\end{proof}

 Il nous reste {\`a} appliquer ce lemme sur $\phi$.
\begin{lemme}\label{lemccv3}
Soit $f_v$ d{\'e}finie dans le lemme formel \ref{lemformel},  $q$ le degr{\'e}
en $Y$ de $P$, $p$ le dergr{\'e} de $p$ et  soit   $T$ un arbre de $\al$ fruits, $F(T)\in F(n)$, alors en reprenant les notations du lemme \ref{lemccv0}, on a
 
\bs\bml\left\vert\prod_{v\in 
V(T)}\inf(\mid\frac{1}{D_v^\Lambda}\mid,\frac{1}{C})f_{v} 
\prod_{v\in V(T)}\f{u}_{\al_v^i}\right\vert\leq      \frac{e^{4\al}\mid \al_w^i\mid^{r-1}}{(\mid n\mid)^{r} }\\\left(\prod_{u\in
     V(T): \mid D_u^\Lambda\mid<\frac{\mid p_v\mid\sqrt{C}}{2}}\sqrt{C}\inf(\mid\frac{1}{D^\Lambda_u}\mid,\frac{1}{C})\right) \left(\prod_{u\in
 V(T): \mid D_u^\Lambda\mid\geq\frac{\mid p_v\mid\sqrt{C}}{2}}\mid\frac{1}{\sqrt{C}}\mid\right)\times\\
\left(\prod_{u\in
V(T)} M \sum_{\substack{(j_1,...j_{\parallel \al_v\parallel})\in
  \{0,1,2\}^{\parallel \al_v\parallel}\\ \sum_{i=1}^{\parallel
    \al_v\parallel}j_i\leq q+1}}\prod_{i=1}^{\parallel
\al_v\parallel}\mid\al_v^i\mid^{j_i}\right),\end{multline*}\es

\end{lemme}
 
\begin{proof} On applique le lemme \ref{Ltruc1} et  le lemme \ref{lemccv0}.
 
\end{proof}  
Dans le cas g\'en\'eral, on d\'emontre le lemme suivant de la m\^eme fa\c con en applicant la
proposition \ref{ARO5}.
\begin{lemme}
Soit $\Lambda(x,t)=\sum_{j=1}^{k}\lambda_j(t)(ix)^j$. Soit $T$ un arbre
fruitier, $t_0$ un r\'eel positif fix\'e  tel que les
$\lambda_i(t)$ soit $C^1$ sur $[0,t_0]$.
Supposons que pour tout $u$ de $T$ on ait :
\[\mid D_u^\Lambda(T,t)\mid<\frac{\mid p_u(T)\mid^{k-2} \sqrt C}{2}.\]
Alors :
\begin{description}
\item[(1)] soit $\mid p_u(T)\mid<\mid \al(u)\mid\sqrt{C}$. 
\item[(2)]soit il existe trois fruits de $T_e(u)$ tel que
\[\mid\al_u^1(T_e(u))\al_u^2(T_e(u))\al_u^3(T_e(u))\mid\geq \frac{\mid p_u(T)\mid^3}{(\al(u))^3}.\]

\item[(3)] soit $u$ est dans $(\cal R)_C$.
\end{description}
\end{lemme}

De plus  on a le
\begin{lemme}
Soit $\Lambda(x,t)=\sum_{j=1}^{k}\lambda_j(t)(ix)^j$. Soit $T$ un arbre
fruitier, $t_0$ un r\'eel positif fix\'e, tel que les
$\lambda_i(t)$ soit $C^1$ sur $[0,t_0]$.
Si pour le noeud $u$ de $T$, on a :
\[\mid D_u^\Lambda(T,t)\mid\geq\frac{\mid p_u(T)\mid^{k-2} \sqrt C}{2},\]
alors
\begin{description}
\item[(1)] soit $\mid p_u(T)\mid<\mid \al(T)\mid\sqrt{C}$, 
\item[(2)]soit il existe trois fruits de $T(u)$ tel que:
\[\mid\al_u^1(T_e(u))\al_u^2(T_e(u))\al_u^3(T_e(u))\mid\geq \frac{\mid p_u(T)\mid^3}{(\al(u))^3},\]

\item[(3)] soit il exist un fruit $\al_v^i(T(u)$ de $T(u)$ tel que
$\al_v^i(T(u)=p_u(T)$.
\item[(4)] soit \[\frac{\mid \partial_tD_u^\Lambda(T,t)\mid}{\mid
D_u^\Lambda(T,t)\mid}\leq\frac{\sum_i\mid \lambda_i'(t)\mid }{\mid
\lambda_{k}\mid}.\]
\end{description}
\end{lemme}

\begin{appendix} 
\section{Propri{\'e}t{\'e}s arithm{\'e}tiques} 

\begin{lemme}[De division]\label{lemDiv}  
 
Soit les $q_i$, $\alpha\geq 2$ entiers relatifs v{\'e}rifiants les 
hypoth{\`e}ses:  
\begin{description}
\item[(A)] $\mid \sum_i q_i\mid \neq \sup_i\mid q_i\mid$ et 
$m= \sup_i\mid q_i\mid=\mid q_j\mid$, 
\item[(B)] $\sum_{i\neq j}\mid q_i \mid \leq \frac{(m)^{\frac{k-1}{k}}}{2^k
\alpha^{\frac{k-1}{k}}}$, avec $k \geq 2$.
\end{description} 

Alors on a
$\mid \sum_{i=1}^{\alpha} q_i^k-(\sum_{i=1}^{\alpha} q_i)^k  \mid \geq \frac{1}{2}\mid 
m^{k-1} \mid$. 
 
\end{lemme} 
 \begin{proof} 
 
Les hypoth{\`e}ses impliquent que  $\sum_{i \neq j}q_i \neq 0$.
On note $n = \sum_i q_i$,
et on pose $j$ donn{\'e} par $m=\mid q_j\mid$, alors
\[\sum_{i=1}^{\alpha} q_i^k-(\sum_{i=1}^{\alpha} q_i)^k=\sum_{s=1}^k C_k^s 
q_j^{k-s}(\sum_{i=1, i\neq j}^{\alpha} q_i)^s+\sum_{i 
\not= j} q_i^k-(n-q_j)^k 
.\] 
Donc 
\bs\[
\sum_{i=1}^{\alpha} q_i^k-(\sum_{i=1}^{\alpha} q_i)^k =q_j^k(\sum_{i=1: i\neq j}^{\alpha} q_i)+\sum_{s=2}^k 
C_k^s q_j^{k-s}(\sum_{i=1, i\neq j}^{\alpha} q_j)^s+
\sum_{i \not= j} q_i^k-(\sum_{i \not= j} q_i)^k.\]\es
\bs\bml\mid\sum_{i=1}^{\alpha} q_i^k-(\sum_{i=1}^{\alpha} q_i)^k \mid \geq\mid 
m^{k-1}(\sum_{i=1:i\neq j}^{\alpha} q_i)\mid
 -\sum_{s=2}^k C_k^s\mid q_j^{k-s} (\sum_{i=1: i\neq j}^{\alpha}  q_i 
)^s\mid\\
-\mid\sum_{i=1}^{\alpha} 
\mid q_i\mid^k+(\sum_{i\neq j} \mid q_i\mid)^k 
\mid  .\end{multline*}\es

Or d'apres les hypoth{\`e}ses (B): 
\[\sum_{s=2}^k C_k^s\mid q_j^{k-s} (\sum_{i=1: i\neq j}^{\alpha}  q_i 
)^s\mid \leq (\sum_{s=2}^k C_k^s)\mid q_j^{k-2} (\sum_{i\neq j}  q_i)^2\mid\] 
et 
$\sum_{i=1}^{\alpha} 
\mid q_i\mid^k\leq \frac{m^{k-1}}{2^k\alpha^k}\alpha$ 
et 
$(\sum_{i\neq j} \mid q_i\mid)^k\leq 
(m^{\frac{k-1}{k}}\frac{\alpha}{2^k\alpha})^k,$
alors 
\bs\[ 
\mid\sum_{i=1}^{\alpha} q_i^k-(\sum_{i=1}^{\alpha} q_i)^k \mid \geq \mid m^{k-1}\mid(\sum_{i=1:i\neq j}^\alpha q_i)\mid-
(\sum_{s=2}^k C_k^s)\mid m^{k-2} (\sum_{i\neq j}  q_i)^2\mid-\frac{1}{4}\mid m \mid^{k-1} 
,\]\es 
\[\mid\sum_{i=1}^{\alpha} q_i^k-(\sum_{i=1}^{\alpha} q_i)^k \mid \geq \mid m^{k-1}(\sum_{i\neq j} 
q_i)\mid-\frac{1}{2}\mid m^{k-1} (\sum_{i=1:i\neq j}^\alpha q_i)\mid \] 
\end{proof}
\begin{proposition}\label{ARO4}
Soit les $q_i$, $\alpha$ entiers relatifs tel que $q_1$ $q_2$ $q_3$ soient les trois plus grands en valeur absolue class{\'e} par ordre d{\'e}croissant. Si un des trois est plus petit que $\frac{\mid q_1\mid}{4}$ alors soit  $q_3$ plus grand que $\frac{\mid q_1\mid^{\frac{2}{3}}}{\al}$ et il existe $q_4$ tel que $\mid q_2q_3q_4\mid>\frac{\mid q_1\mid^2}{\al^2}$, soit  
\[\mid \sum_{i=1}^{\alpha} q_i^3-(\sum_{i=1}^{\alpha} q_i)^3  \mid \geq \frac{1}{2}\mid 
q_1^{2} \mid.\]  
 
\end{proposition}

\begin{proof}
\[\sum_{i=1}^{\alpha} q_i^3-(\sum_{i=1}^{\alpha} q_i)^3 =q_1^3+q_2^3+(\sum_{i=3}^{\alpha} q_i)^3-( q_1+q_2+\sum_{i=3}^{\alpha} q_i)^3+\sum_{i=3}^{\alpha} q_i^3-(\sum_{i=3}^{\alpha} q_i)^3.\] On obtient alors
\[\sum_{i=1}^{\alpha} q_i^3-(\sum_{i=1}^{\alpha} q_i)^3 =3(q_1+q_2)(q_2+\sum_{i=3}^{\alpha} q_i)(q_1+\sum_{i=3}^{\alpha} q_i)+\sum_{i=3}^{\alpha} q_i^3-(\sum_{i=3}^{\alpha} q_i)^3.\] Si $q_1+q_2=0$ on applique le lemme \ref{lemDiv}, sinon on remarque que si $\mid q_3\mid \leq \frac{\mid q_1\mid}{4}$ et que tous les autres $q_i$ sont major{\'e}s en valeur absolu par $\frac{\mid q_1\mid^{\frac{2}{3}}}{\al}$ alors $\mid (q_1+q_2)(q_2+\sum_{i=3}^{\alpha} q_i)(q_1+\sum_{i=3}^{\alpha} q_i)\mid \geq \frac{\mid q_1\mid^2}{4}$. En applicant le lemme \ref{lemDiv} sur $\sum_{i=3}^{\alpha} q_i^3-(\sum_{i=3}^{\alpha} q_i)^3$, on trouve qu'il est d'ordre $\mid 3(q_3)(\sum_{i=4}^{\alpha} q_i)(\sum_{i=3}^{\alpha} q_i)\mid +\mid\frac{q_1^2}{\al^2}\mid$. Or $\mid 3(q_3)(\sum_{i=4}^{\alpha} q_i)(\sum_{i=3}^{\alpha} q_i)\mid \leq 3(q_3)^2(\al-3)^2\sup_{i\geq 4}(\mid q_i\mid)$. Donc il faut que $ (q_3)^2(\al-3)^2\sup_{i\geq 4}(\mid q_i\mid)\geq \frac{\mid q_1\mid}{4}$. 
\end{proof}
Qui se g\'en\'eralise ainsi:
\begin{proposition}\label{ARO5}
Soit les $q_i$, $\alpha$ entiers relatifs tel que $q_1$ $q_2$ $q_\al$ soient
class{\'e} par ordre d{\'e}croissant ( en valeur absolue) tel que
pour tout $i$, $q_i\neq \sum_jq_j$. Soit $\lambda_1,...\lambda_k$, $k$ nombre
complex avec $\lambda_k$ non nul.  Alors soit on a $\mid q_2\mid>\frac{\mid
q_1\mid^{\frac{(k-1)}{k}}}{\al^2}$, soit  $\mid q_1\mid^{k-1}\leq
\frac{\sum_i\mid\lambda_i\mid}{\mid \lambda_k\mid}$, soit
\[\mid \sum_{j=1}^k\lambda_j\sum_{i=1}^{\alpha} q_i^j-(\sum_{i=1}^{\alpha} q_i)^j  \mid \geq \frac{1}{2}\mid 
q_1 \mid^{k-1}.\]  
 
\end{proposition}

\begin{proof} Si 
$\mid q_2\mid_leq \frac{\mid
q_1\mid^{\frac{(k-1)}{k}}}{\al^2}$ alors d'apr\`es la proposition \ref{lemDiv} soit $\mid q_1\mid^{k-1}\leq
\frac{\sum_i\mid\lambda_i\mid}{\mid \lambda_k\mid}$, soit on a 
$\mid \sum_{j=1}^k\lambda_j\sum_{i=1}^{\alpha} q_i^j-(\sum_{i=1}^{\alpha} q_i)^j  \mid \geq \frac{1}{2}\mid 
q_1 \mid^{k-1}$.   
\end{proof}
Nous allons maintenant d{\'e}montrer les lemmes qui relient les produits
d'op{\'e}rateurs aux diviseurs. Dans un premier temps on suppose que $\Lambda$
est ind\'ependant du temps.
\begin{lemme}\label{lemaqD1}
Soit un arbre fruitier $T$ de $\be\geq 2$ noeuds, on note $C(T)=\{v\in V(T): (\mid D_v^\Lambda(T)\mid \leq \frac{1}{t} \hbox{ ou } A_v=0)\hbox{ et } \forall u\leq v, \mid D_u^\Lambda(T)\mid > \frac{1}{t}\}$ (cela signifie que $v$ est le noeud tel qu'il n'existe aucun noeud inf{\'e}rieur au sens de la racine   qui ont un diviseur plus petit que $\frac{1}{t}$ et que le diviseur en $w(T)$ est plus petit que $\frac{1}{t}$ sauf si $w(T)$ est un noeud extr{\^e}me). Si pour le premier noeud de $T$, $v_1$, $\mid D_{v_1}^\Lambda(T)\mid > \frac{1}{t}$, alors on a:
\bs\bml\sum_{S(T)}\left(\rov{S}[V(T),\Lambda]\right)(t)=\frac{1}{D_{v_1}^\Lambda(T)}\prod_{T_i\in v_1(T)}\sum_{S(T_i)}\left(\rov{S}[V(T_i),\Lambda]\right)(t)+\\\frac{1}{D_{v_1}^\Lambda(T)}\sum_{v\in C(T)}D_{v}^\Lambda(T)\sum_{S(T)}\left(\rov{S}[V(T),\Lambda]\right)(t) -\frac{1}{D_{v_1}^\Lambda(T)}\sum_{v\in C(T)}\left(\rov{S}[V(T)\setminus\{ v\},\Lambda]\right)(t).\end{multline*}\es
\end{lemme}
\begin{proof}
La d{\'e}monstration se fait par r{\'e}curence sur le nombre de noeud de $T$:

Si $T$ poss{\`e}de deux noeuds, $v_1$ et $v_2$:
\bs\[\sum_{S(T)}\left(\rov{S}[V(T),\Lambda]\right)(t)=e^{t(\sum_l\Lambda(\al_{v_1}^l(T))-D_{v_1}^\Lambda(T))}\int_0^te^{sD_{v_1}^\Lambda(T)}e^{-sD_{v_2}^\Lambda(T)}ds\int_0^se^{xD_{v_2}^\Lambda(T)}dx.\]\es
On int{\`e}gre par partie:
\begin{eqnarray*}
e^{sD_{v_1}^\Lambda(T)}&\to&\frac{e^{sD_{v_1}^\Lambda(T)}}{D_{v_1}^\Lambda(T)}\\
e^{-sD_{v_2}^\Lambda(T)}ds\int_0^se^{xD_{v_2}^\Lambda(T)}dx&\to&-D_{v_2}^\Lambda(T)e^{-sD_{v_2}^\Lambda(T)}ds\int_0^se^{xD_{v_2}^\Lambda(T)}dx+1,
\end{eqnarray*}
qui donne:
\bs\bml e^{t(\sum_l\Lambda(\al_{v_1}^l(T))-D_{v_1}^\Lambda(T))}\int_0^te^{sD_{v_1}^\Lambda(T)}e^{-sD_{v_2}^\Lambda(T)}ds\int_0^se^{xD_{v_2}^\Lambda(T)}dx=\\\frac{1}{D_{v_1}^\Lambda(T)}e^{t(\sum_l\Lambda(\al_{v_1}^l(T))-D_{v_2}^\Lambda(T))}\int_0^te^{xD_{v_2}^\Lambda(T)}dx-\\\frac{1}{D_{v_1}^\Lambda(T)}e^{t(\sum_l\Lambda(\al_{v_1}^l(T))-D_{v_1}^\Lambda(T))}\int_0^te^{sD_{v_1}^\Lambda(T)}+\\\frac{D_{v_2}^\Lambda(T)}{D_{v_1}^\Lambda(T)}e^{t(\sum_l\Lambda(\al_{v_1}^l(T))-D_{v_1}^\Lambda(T))}\int_0^te^{sD_{v_1}^\Lambda(T)}e^{-sD_{v_2}^\Lambda(T)}ds\int_0^se^{xD_{v_2}^\Lambda(T)}dx.\end{multline*}\es
Comme $v_2\in C(T)$, on a bien la r{\'e}ponse. 
Supposons la propri{\'e}t{\'e} vrai pour tout arbre de $\be-1\geq 2$ noeuds:
Le lemme \ref{LS(T_v)} appliqu{\'e} au premier noeud donne:
\bs\bml(\sum_{T'\in S(T)}\rov{\cal S}[V(T'),\Lambda])=\\<\rov{\cal S}[V(T_e(v_1)),\Lambda],\prod_{T_j\in v_1(T)}\left(e^{-t\Lambda\left(\sum_{u\in V(T_j)}\sum_l\al_u^l(T_j)\right)}\sum_{T'\in S(T_j)}\rov{\cal S}[V(T'),\Lambda]\right)>.\end{multline*}\es
On note $v(k)$ le premier noeud de $T_k$ ( remarquons qu'il est inveriant sur $S(T_k)$). On int{\`e}gre par partie:
\[
e^{sD_{v_1}^\Lambda(T)}\to\frac{e^{sD_{v_1}^\Lambda(T)}}{D_{v_1}^\Lambda(T)}\]
et 
\bs\bml\prod_{T_j\in v_1(T)}\left(e^{-t\sum_{u\in V(T_j)}\sum_l\Lambda\left(\al_u^l(T_j)\right)}\sum_{T'\in S(T_j)}\rov{\cal S}[V(T'),\Lambda]\right)\to\\
-\sum_{T_j\in v_1(T)}\left(\sum_{u\in V(T_j)}\sum_l\Lambda\left(\al_u^l(T_j)\right)-\Lambda\left(\sum_{u\in V(T_j)}\sum_l\al_u^l(T_j)\right)\right)\times\\\prod_{T_j\in v_1(T)}\left(e^{-t\sum_{u\in V(T_j)}\sum_l\Lambda\left(\al_u^l(T_j)\right)}\sum_{T'\in S(T_j)}\rov{\cal S}[V(T'),\Lambda]\right)+\\\sum_{T_k\in v_1(T)}\left(e^{-t\sum_{u\in V(T_k)}\sum_l\Lambda\left(\al_u^l(T_j)\right)}\sum_{T'\in S(T_k)}\rov{\cal S}[V(T')\setminus \{v(k)\},\Lambda]\right)\times\\\prod_{T_j\in v_1(T)\setminus\{T_k\}}\left(e^{-t\sum_{u\in V(T_j)}\sum_l\Lambda\left(\al_u^l(T_j)\right)}\sum_{T'\in S(T_j)}\rov{\cal S}[V(T'),\Lambda]\right).\end{multline*}\es
Enfin sur tous les arbre $T_k$ tel que $v(k)\not\in C(T)$ on applique l'hypoth{\`e}se de r{\'e}curence sur $T_k$ dans:
\[-e^{t\left(\sum_{u\in V(T)}\sum_l\Lambda\left(\al_u^l(T)\right)-D_{v_1}^\Lambda(T)\right)}\int_0^te^{D_{v_1}^\Lambda(T)s}\frac{D_{v(k)}^\Lambda(T_k)}{D_{v_1}^\Lambda(T)}\times\]\[\prod_{T_j\in v_1(T)}\left(e^{-t\sum_{u\in V(T_j)}\sum_l\Lambda\left(\al_u^l(T_j)\right)}\sum_{T'\in S(T_j)}\rov{\cal S}[V(T'),\Lambda]\right)s)ds,\]
qui donne bien la r{\'e}ponse en remarquant que $D_{v}^\Lambda(T_k)=D_{v}^\Lambda(T)$.

\end{proof}

\begin{lemme}\label{lemaqD}
Soit $C$ et $\Lambda$ des r{\'e}els positifs $T\in T_{\al,\be}^e$ un arbre, alors  pour tout $t\leq \frac{1}{C}$ il existe
$k(T)=Card(\{v\in V(T):\mid D_v^{\Lambda}(T)\mid\leq C\})\leq \be$:
\[\left\vert\sum_{T'\in S(T)}\rov{S}[V(T'),\Lambda]\right\vert\leq \sum_{G\in G(\cal T)}\prod_{T'\in G} t^{k(T')}\prod_{v\in V(T'):\mid D_v^{\Lambda}(T')\mid> C}\frac{1}{\mid D_v^{\Lambda}(T')\mid},\]
( pour la d{\'e}finition de $G(\cal T)$ cf \ref{DdecompT} \pageref{DdecompT}).
\end{lemme}
\begin{proof}
Nous allons d{\'e}montrer cette propri{\'e}t{\'e} par r{\'e}curence sur le nombre de noeuds de $T$ . Si $T$ poss{\`e}de un noeuds, dans le cas o{\`u} $\mid D_v^\Lambda(T)\mid >\frac{1}{t}$, la prori{\'e}t{\'e} est vraie. Si $\mid D_v^\Lambda(T)\mid \leq \frac{1}{t}$, alors:
\[\mid e^{t(\sum_l\Lambda(\al_{v}^l(T))-D_{v}^\Lambda(T))}\int_0^te^{sD_{v}^\Lambda(T)}ds\mid\leq t.\]
Supposons la propri{\'e}t{\'e} vraie pour tout arbre de $\be-1\geq 1$ noeuds. Soit un arbre $T$ de $\be$ noeuds et $C=\frac{1}{t}$. Si pour le premier noeud de $T$, $v_1$, $\mid D_{v_1}^\Lambda(T)\mid \leq \frac{1}{t}$, on applique le lemme \ref{LS(T_v)} au premier noeud et l'hypoth{\`e}se de r{\'e}curence sur les $T_i\in v_1(T)$, qui donne:
\bs\[\left\vert\sum_{T'\in S(T)}\rov{S}[V(T'),\Lambda]\right\vert\leq\int_0^t\prod_{T_i\in v_1(T)}\left\vert\sum_{G\in G(\cal T_i)}\prod_{T'\in G} s^{k(T')}\prod_{\substack{v\in V(T'):\\\mid D_v^{\Lambda}(T')\mid> C}}\frac{1}{\mid D_v^{\Lambda}(T')\mid}\right\vert ds,\]\es
 qui donne la r{\'e}ponse. Dans le cas o{\`u}  
$\mid D_{v_1}^\Lambda(T)\mid > \frac{1}{t}$, si pour tout $w\in C(T)$ on a $\mid D_w^\Lambda\mid\leq \frac{1}{t}$, on applique le lemme \ref{lemaqD1}:
\bs\bml\sum_{S(T)}\left(\rov{S}[V(T),\Lambda]\right)(t)=\frac{1}{D_{v_1}^\Lambda(T)}\prod_{T_i\in v_1(T)}\sum_{S(T_i)}\left(\rov{S}[V(T_i),\Lambda]\right)(t)+\\\frac{1}{D_{v_1}^\Lambda(T)}\sum_{v\in C(T)}D_{v}^\Lambda(T)\sum_{S(T)}\left(\rov{S}[V(T),\Lambda]\right)(t) -\frac{1}{D_{v_1}^\Lambda(T)}\sum_{v\in C(T)}\left(\rov{S}[V(T)\setminus\{ v\},\Lambda]\right)(t).\end{multline*}\es Sur le premier terme, on applique alors l'hypoth{\`e}se de r{\'e}curence sur les $T_i$. Sur le deuxi{\`e}me terme, on applique le lemme \ref{LS(T_v)} au premier noeud et l'hypoth{\`e}se de r{\'e}curence sur les $T_i\in v_1(T)$, qui donne:
\bs\begin{eqnarray*}\left\vert\sum_{T'\in
S(T)}\rov{S}[V(T'),\Lambda]\right\vert&\leq&\int_0^t\prod_{T_i\in
v_1(T)}\left\vert\sum_{G\in G(\cal T_i)}\prod_{T'\in G} s^{k(T')}\prod_{\substack{v\in
V(T'):\\\mid D_v^{\Lambda}(T')\mid> C}}\frac{1}{\mid
D_v^{\Lambda}(T')\mid}\right\vert(s)ds\\
&\leq&t\prod_{T_i\in v_1(T)}\left\vert\sum_{G\in G(\cal T_i)}\prod_{T'\in G} s^{k(T')}\prod_{v\in V(T'):\mid D_v^{\Lambda}(T')\mid> C}\frac{1}{\mid D_v^{\Lambda}(T')\mid}\right\vert.\end{eqnarray*}\es
Or $\frac{t}{D_v^{\Lambda}(T_i)}\leq 1$ quand $v\in C(T)$. Pour le troisi{\`e}me terme on applique l'hypoth{\`e}se de r{\'e}curence et le fait que $\frac{1}{D_{v_1}^{\Lambda}(T_i)}\leq t$.
 Si il existe  $w\in C(T)$ tel que$\mid D_w^\Lambda\mid> \frac{1}{t}$, alors au dernier noeud multiple $u$ inferieur {\`a} w, on applique le lemme \ref{LS(T_v)}:
\bml\sum_{T'\in S(T)}(\rov{ S}[V(T),\Lambda])=\\<\sum_{T'\in S(T_e(u))}\rov{ S}[V(T)\Lambda],\prod_{T''\in T_u}e^{-t\sum_{z\in V(T'')}\sum_l\al_z^l(T'')}(\rov{ S}[V(T''),\Lambda])>.\end{multline*}
Si on note $T_i$ l'arbre de $T_u$ qui contient $w$ ( c'est un tronc) on a:
\[(\rov{ S}[V(T_i),\Lambda])=<\rov{\prod_{z<w}}S_{D_z^\Lambda},\frac{1-e^{-tD_z^\Lambda}}{D_z^\Lambda}>,\] qui donne en rempla\c cant:
\bml\sum_{T'\in S(T)}(\rov{ S}[V(T),\Lambda])=\\\frac{1}{D_z^\Lambda}\left(\sum_{T'\in S(T)}(\rov{ S}[V(T)\setminus \{w\},\Lambda])-\sum_{T'\in S(T(r,w))}(\rov{ S}[V(T)\setminus \{w\},\Lambda])\right).\end{multline*} Il nous reste {\`a} appliquer l'hypoth{\`e}se de r{\'e}curence sur $T\setminus T(w)$ pour le premier terme et sur le sous arbre $T(r,w)$ pour le deuxi{\`e}me terme. 

Le lemme se conclu en remarquant que pour $C<\frac{1}{t}$ si $C\leq \mid D_z^\Lambda\mid \leq \frac{1}{t}$, on a tout simplement que $t\leq \mid D_z^\Lambda\mid $ et il suffit de remplacer pour avoir la r{\'e}ponse.
\end{proof}

De m\^eme on d\'emontre dans le cas o\`u les $\lambda_i$ sont variable:
\begin{lemme}
Soit un tronc $T$ et $f$ une fonction d\'erivable de $\R$ dans $\C$. Soit
\[\Lambda(t,x)=\sum_{j=0}^p\lambda_{2j}(t)(ix)^{2j}+\sum_{j=0}^q\lambda_{2j+1}(ix)^{2j+1}.\]
On suppose que les $\lambda_i$ sont $C^1$. Soit $t_0$ tel que pour tout $t\in
[0,t_0]$, $\lambda_{2p}(t)>0$, si on note:
\[M(\Lambda,t)=e^{\sup_{s\in [0,t]x\in \R}\Re(\Lambda(t,x))},\]
 alors: 
\bml\left\vert<\sum_{T'\in S(T)}\rov{S}[V(T'),\Lambda],f>\prod_{v\in
V(T),i}e^{L(t,\sum_{v,i}\al_v^i))}\right\vert\leq \\\sum_{G\in G(\cal T)}\prod_{T'\in G}
\sup_{x\in [0,t]}\{\prod_{v\in V(T')}M(\Lambda,t)((\prod_ie^{L(t,\al_v^i)})\frac{1}{\mid D_v^{\Lambda}(T',x)\mid+\frac{1}{t}})\times\\(\mid f(t)\mid(1+(\sum_{v\in V(T')}\frac{\mid \frac{dD_v^{\Lambda}(T',x)}{dx}\mid}{\mid D_v^{\Lambda}(T',x)\mid+\frac{1}{t}}) +t\mid f'(t)\mid)\}.\end{multline*}
 
\end{lemme}
\begin{proof}
On se donne
$D_{v_1}^\Lambda(s,T)=D_{v_1}^\Lambda(s,T)^-+D_{v_1}^\Lambda(s,T))^+$, o\`u :
 pour tout $s\in [0,t]$  $\mid 
D_{v_1}^\Lambda(s,T)^-\mid\leq \frac{2}{t}$, et pour tout $s\in [0,t]$ $\mid 
D_{v_1}^\Lambda(s,T))^+\mid> \frac{1}{t}$. De plus $D_{v_1}^\Lambda(s,T)^+$ est $C$.

On montre le lemme par r\'ecurence sur le nombre de noeud $T$:

A l'ordre 1 si  $\mid 
D_{v_1}^\Lambda(s,T))\mid\leq \frac{1}{t}$ pour tout $s\in [0,t]$ la
propri\'et\'e est vraie. S'il existe $s$, tel que $\mid 
D_{v_1}^\Lambda(t,T))\mid>\frac{1}{t}$ alors:
\bs\bml e^{(\sum_lL(\al_{v_1}^l(T))-{\cal
D}_{v_1}^\Lambda(t,T))}\int_0^te^{{\cal D}_{v_1}^\Lambda(s,T)}f(s)ds=\\
e^{(\sum_lL(\al_{v_1}^l(T))-{\cal
D}_{v_1}^\Lambda(t,T))}\int_0^t\frac{D_{v_1}^\Lambda(s,T)-D_{v_1}^\Lambda(s,T)^-}{D_{v_1}^\Lambda(s,T)^+}e^{{\cal D}_{v_1}^\Lambda(s,T)}f(s)ds=\\
\frac{e^{(\sum_lL(\al_{v_1}^l(T))}f(t)}{D_{v_1}^\Lambda(t,T)^+}-e^{(\sum_lL(\al_{v_1}^l(T))-{\cal
D}_{v_1}^\Lambda(t,T))}\int_0^t\frac{e^{{\cal
D}_{v_1}^\Lambda(s,T)}}{D_{v_1}^\Lambda(s,T))^+}f'(s)ds+\\
e^{(\sum_lL(\al_{v_1}^l(T))-{\cal
D}_{v_1}^\Lambda(t,T))}\int_0^t\frac{e^{{\cal
D}_{v_1}^\Lambda(s,T)}}{D_{v_1}^\Lambda(s,T))}\frac{f(s)\frac{dD_{v_1}^\Lambda(s,T))^+}{dt}}{D_{v_1}^\Lambda(s,T)^+)}ds-\\e^{(\sum_lL(\al_{v_1}^l(T))-{\cal
D}_{v_1}^\Lambda(t,T))}\int_0^t\frac{D_{v_1}^\Lambda(s,T)^-}{D_{v_1}^\Lambda(s,T)^+}e^{{\cal D}_{v_1}^\Lambda(s,T)}f(s)ds.\end{multline*}\es
Du fait que $\lambda_{2p}(t)>0$, alors pour $x$ suffisement grand, pour tout
$s\leq t$: $\Re(\Lambda(t,x)-\Lambda(s,x))\leq 0$ qui conclu le lemme a l'ordre
1. 

A l'odre $n$:
\bml<\sum_{T'\in S(T)}\rov{S}[V(T'),\Lambda],f>=\\e^{(\sum_lL(\al_{v_1}^l(T))-{\cal
D}_{v_1}^\Lambda(t,T))}\int_0^te^{{\cal D}_{v_1}^\Lambda(s,T)}....e^{(-{\cal
D}_{v_i}^\Lambda(t,T))}\int_0^{s_i}e^{{\cal
D}_{v_i}^\Lambda(s_{i+1},T)}...f(s_n)ds....ds_n.\end{multline*}
On suppose que $i$ est le plus grand indice tel qu'il existe $s$, $\mid 
D_{v_i}^\Lambda(s,T))\mid>\frac{1}{t}$ alors en int\'egrant par  partie suivant
la variable $s_{i-1}$ on a
\bs\bml e^{(\sum_lL(\al_{v_1}^l(T))-{\cal
D}_{v_1}^\Lambda(t,T))}\int_0^te^{{\cal D}_{v_1}^\Lambda(s,T)}....e^{-{\cal
D}_{v_{i}}^\Lambda(t,T))}\int_0^{s_{i}}e^{{\cal
D}_{v_i}^\Lambda(s_{i+1},T)}...f(s_n)ds....ds_n=\\
e^{(\sum_lL(\al_{v_1}^l(T))-{\cal
D}_{v_1}^\Lambda(t,T))}\int_0^te^{{\cal D}_{v_1}^\Lambda(s,T)}....\\
e^{-{\cal
D}_{v_{i-1}}^\Lambda(t,T))}\int_0^{s_{i-1}}e^{{\cal
D}_{v_{i-1}}^\Lambda(s_{i},T)}\frac{1}{D_{v_{i}}^\Lambda(s_{i},T)^+}e^{-{\cal
D}_{v_{i+1}}^\Lambda(s_i,T))}\int_0^{s_{i}}...f(s_n)ds....ds_n-\\
e^{(\sum_lL(\al_{v_1}^l(T))-{\cal
D}_{v_1}^\Lambda(t,T))}\int_0^te^{{\cal D}_{v_1}^\Lambda(s,T)}....\\
e^{-{\cal
D}_{v_{i}}^\Lambda(s_{i-1},T))}\int_0^{s_{i}}e^{{\cal
D}_{v_{i}}^\Lambda(s_{i+1},T)}\frac{d}{ds_{i+1}}\frac{1}{D_{v_{i}}^\Lambda(s_{i+1},T)^+}e^{-{\cal
D}_{v_{i+1}}^\Lambda(s_{i+1},T))}\int_0^{s_{i+1}}...f(s_n)ds....ds_n-\\
e^{(\sum_lL(\al_{v_1}^l(T))-{\cal
D}_{v_1}^\Lambda(t,T))}\int_0^te^{{\cal D}_{v_1}^\Lambda(s,T)}....e^{-{\cal
D}_{v_{i}}^\Lambda(s_{i-1},T))}\int_0^{s_{i}}\frac{D_{v_{i}}^\Lambda(s_{i+1},T)^-}{D_{v_{i}}^\Lambda(s_{i+1},T)^+}e^{{\cal
D}_{v_i}^\Lambda(s_{i+1},T)}...f(s_n)ds....ds_n.\end{multline*}\es
On pose alors \[g=\frac{1}{D_{v_{i}}^\Lambda(s_{i},T)}e^{-{\cal
D}_{v_{i-1}}^\Lambda(s_i,T))}\int_0^{s_{i}}...f(s_n)ds_{i+1}....ds_n\]
 Soit $T''$ le sous arbre de $T$
auquel on retire les noeuds sup\'erieur \`a $v_{i-1}$,  on applique
l'hypoth\`ese de r\'ecurence sur:
\[\left\vert<\sum_{T'\in S(T'')}\rov{S}[V(T'),\Lambda],g>\prod_{v\in
V(T),i}e^{L(t,\sum_{v,i}\al_v^i))}\right\vert.\]
De plus on pose:
\[h=\int_0^{s_{i}}e^{{\cal
D}_{v_{i}}^\Lambda(s_{i+1},T)}\frac{d}{ds_i}\frac{1}{D_{v_{i}}^\Lambda(s_{i+1},T)}e^{-{\cal
D}_{v_{i+1}}^\Lambda(s_{i+1},T))}\int_0^{s_{i+1}}...f(s_n)ds_{i+2}....ds_n,\]
et on applique l'hypoth\`ese de r\'ecurence \`a  $T(v_{i})$, qui conclu le
lemme dans ce cas l\`a.

\end{proof}
Qui implique les th\'eor\`emes suivants:
\begin{theoreme}\label{TAG(T)1}
Soit $C$ un r\'eel positif. Soit
$\Lambda(t,x)=\sum_{j=0}^p\lambda_{2j}(t)(ix)^{2j}+\sum_{j=0}^q\lambda_{2j+1}(ix)^{2j+1}$.
Soit $t_0$ tel que pour tout $t\in [0,t_0]$, $\lambda_{2p}(t)>0$, alors il
existe $M(\lambda),P)$ tel que:
Soit $T\in T_{\al,\be}^e$ un arbre, alors  pour tout $t\leq \frac{1}{C}$ il existe
$k(T)=Card(\{v\in V(T):\mid D_v^{\Lambda}(T)\mid\leq C\})\leq \be$:
\bs\bml\left\vert\sum_{T'\in S(T)}\rov{S}[V(T'),\Lambda]\prod_{v\in
V(T)}\prod_ie^{L(t,\al_v^i)}\right\vert\leq \\\sum_{G\in G(\cal T)}\prod_{T'\in G}
t^{k(T')}\sup_{x\in [0,t]}\prod_{v\in V(T'):\mid D_v^{\Lambda}(T')\mid> C}\frac{1}{\mid D_v^{\Lambda}(T',x)\mid},\end{multline*}\es
( pour la d{\'e}finition de $G(\cal T)$ cf \ref{DdecompT} \pageref{DdecompT}).

\end{theoreme}


De m\^eme on a la propri\'et\'e r\'eciproque:
\begin{theoreme}\label{TAG(T)2}
Soit $C$ un r\'eel positif. Soit
$\Lambda(t,x)=\sum_{j=0}^p\lambda_{2j}(t)(ix)^{2j}+\sum_{j=0}^q\lambda_{2j+1}(ix)^{2j+1}$.
Soit $t_0$ tel que pour tout $t\in [0,t_0]$, $\lambda_{2p}(t)<0$, alors il
existe $M(\lambda),P)$ tel que:
Soit $T\in T_{\al,\be}^e$ un arbre, alors  pour tout $t\leq \frac{1}{C}$ il existe
$k(T)=Card(\{v\in V(T):\mid D_v^{\Lambda}(T)\mid\leq C\})\leq \be$:
\bml\left\vert e^{L(t,\sum_{v,i}\al_v^i))}\sum_{T'\in S(T)}\lov{S}[V(T'),\Lambda]\right\vert\leq \\\sum_{G\in G(\cal T)}\prod_{T'\in G} t^{k(T')}\prod_{v\in V(T'):\mid D_v^{\Lambda}(T')\mid> C}\frac{1}{\mid D_v^{\Lambda}(T')\mid},\end{multline*}
( pour la d{\'e}finition de $G(\cal T)$ cf \ref{DdecompT} \pageref{DdecompT}).

\end{theoreme}
 \section{Formalisme des cha{\^\i}nes} 

\begin{definition}\label{Dch1} 
 
Soit $(X)$ un ensemble de $I$ sous ensemble ordonn{\'e} fini $X_i$ d'un corps commutatif 
$K$. On appelle maillon, les {\'e}l{\'e}ments de la cha{\^\i}ne 
{\'e}l{\'e}mentaire $X_i$, et on note $X_i^j$  les {\'e}l{\'e}ments de $X_i$ 
ordonn{\'e}s par l'indice j, et $I_i$  le nombre de maillons de la 
cha{\^\i}ne {\'e}l{\'e}mentaire $X_i$. Une chaine $(X)$ est donc une partie 
de  $I$ sous ensemble ordonn{\'e} finie $X_i$ d'un corps commutatif 
$K$. 
 
On d{\'e}finit l'application injective $\sigma$: 
\[\sigma:\{1,...,\sum_iI_i=J\}\longrightarrow \prod_i\{1,...I_i\}\] 
\[k \longrightarrow (j(1,k),...,j(i,k),j(I,k))\] 
  qui v{\'e}rifie: 
\[\hbox{ si } k \geq k'\Rightarrow  \forall i,j(i,k) \geq j(i,k').\] 
 
On dit alors que $\sigma$ est une application respectant la non commutation 
des maillons sur $X$. 
\end{definition}

$\sigma(X)$ d{\'e}signe l'ensemble des applications respectant la non 
commutation des maillons sur $(X)$. $J$ le nombre d'{\'e}l{\'e}ments de $(X)$ est aussi appel{\'e} la taille de $(X)$, etrespectivement $I_i$ la taille de $X_i$.

\begin{definition}\label{Dchres} 
 
Une cha{\^\i}ne $X$ est dite r{\'e}sonante si 
\[\sum_{i,j}X_i^j=0.\] 
De plus elle est  sous-r{\'e}sonante si il existe $\sigma$ de $\sigma(X)$, $k$ et $p$ deux entiers tels que:
\[\sum_{i} 
X_{i}^{j(i,k)}=\sum_{i} 
X_{i}^{j(i,p)}\] 
et \[\forall \sigma \in \sigma(X),\forall n \leq I \hbox{  }\sum_{i} 
X_{i}^{j(i,n)} \neq 0.\] 

\end{definition} 
 \begin{definition}  
  
Un chainon est un couple  $((X),\sigma) $, o{\`u} $\sigma$ est un  
{\'e}l{\'e}ment de $\sigma(X)$. L'ensemble des chainons est alors not{\'e}:  
\[((X), \sigma(X))\] 
De plus on dit que deux chainons $A=((X), \sigma)$ et $B=((Y),  
\sigma')$ sont {\'e}gaux si $(X)=(Y)$ et si pour tout $i, j$,  
$X_i^{j(i,n)}=Y_i^{j'(i,n)}$; o{\`u} $j(i,n)$ et $j'(i,n)$ sont les  
applications d{\'e}duitent respectivement de $ \sigma$ et $ \sigma'$.  
\end{definition}

Soit, une cha{\^\i}ne $(X)$ d{\'e}termin{\'e}e par ses cha{\^\i}nes 
{\'e}l{\'e}mentaires de taille $I_i$, $G_{\un{X}}$ d{\'e}signe l'ensemble 
des multi-entiers $\un{j}$ de $\N^I$ tels que: 
\[\forall i, j_i \in \{0,1,.. I_i\}.\] 
On rappelle que par convention: 
\[\forall i, X_i^0=0.\] 
 
\begin{definition}\label{Dchgen} 
 
 
On note ainsi $G_X$ des applications distinctes qui {\`a} $i\in\{1,...I\}$ associe $j(i)\in \{0,1,...I_i\}$.
\end{definition}

On  suppose que $S$ est un op{\'e}rateur de type d{\'e}fini en \ref{Ds}: 

 On donne d'abord deux notations pour les produits d'op{\'e}rateurs:
\[\lov{\prod_{i=1}^n}S_{w_i}=\lov{\prod_{i}}S_{w_i}=S_{w_n}....S_{w_1}\]
et inversement \[\rov{\prod_{i=1}^n}S_{w_i}=\rov{\prod_{i}}S_{w_i}=S_{w_1}....S_{w_n}.\]  
\begin{lemme}{Du comportement moyen g{\'e}n{\'e}ralis{\'e}}\label{lemCMgen} 
 
Soit $X$ une cha{\^\i}ne quelconque de taille $I$: 
 \[\sum_{\sigma \in \sigma(X)}(\lov{\prod_{\tau =1}^I}S_{\sum_{i} 
   X_{i}^{j(i,\tau)}}) =\prod_{i}(\lov{\prod_j}  S_{X_i^j})\] 
\end{lemme}

 \begin{proof} 
La d{\'e}monstration s'effectue par r{\'e}currence. 
A l'ordre 2:
\begin{eqnarray*} 
\partial_t(S_{w_1+w_2}S_{w_1}+S_{w_1+w_2}S_{w_2})&=&(w_1+w_2)(S_{w_1+w_2}S_{w_1}+S_{w_1+w_2}S_{w_2})+\\ 
&&S_{w_1}+S_{w_2}\\ 
\partial_t(S_{w_1})(S_{w_2})&=&(w_1+w_2)(S_{w_1})(S_{w_2})+S_{w_1+w_2}S_{w_2})+\\ 
&&S_{w_1}+S_{w_2} .
\end{eqnarray*} 
Or au moins $w_1$ ou $w_2$ on le m{\^e}me signe que $w_1+w_2$ d'o{\`u} 
\[(S_{w_1})(S_{w_2})(t_{w_1+w_2})=(S_{w_1+w_2}S_{w_1}+S_{w_1+w_2}S_{w_2})(t_{w_1+w_2})=0.\] 
 
Reste donc {\`a} montrer que la propri{\'e}rt{\'e} est vraie {\`a} tout 
ordre. Supposons la vrai \` a l'ordre $n-1$: 
\begin{eqnarray*} 
\sum_{\sigma \in \sigma(X)}(\lov{\prod_{\tau =1}^I}S_{\sum_{i} X_{i}^{j(i,\tau)}} 
)&=&\sum_p(S_{\sum_{i} X_{i}^{I_i}}\sum_{\sigma \in 
    \sigma(X\setminus X_p^{I_p})}\lov{\prod_{\tau 
    =1}^I}S_{\sum_{i} X_{i}^{j(i,\tau)}})\\ 
\partial_t\sum_{\sigma \in \sigma(X)}\lov{\prod_{\tau =1}^I}S_{\sum_{i} X_{i}^{j(i,\tau)}} 
&=&\sum_p\partial_t(S_{\sum_{i} X_{i}^{I_i}}(\lov{\prod_{j<I_p} 
  S_{X_p^j}})\prod_{i\neq p}(\lov{\prod_j  S_{X_i^j}})). 
\end{eqnarray*} 
Or 
 
\begin{eqnarray*} 
\partial_t(S_{\sum_{i} X_{i}^{I_i}}(\lov{\prod_{j<I_p} 
  S_{X_p^j}})\prod_{i\neq p}(\lov{\prod_j 
  S_{X_i^j}}))&=& (\sum_{i} X_{i}^{I_i})(S_{\sum_{i} X_{i}^{I_i}}(\lov{\prod_{j<I_p} 
  S_{X_p^j}})\times\\ 
&&\prod_{i\neq p}(\lov{\prod_j 
  S_{X_i^j}}))+(\lov{\prod_{j<I_p} 
  S_{X_p^j}})\prod_{i\neq p}(\lov{\prod_j 
  S_{X_i^j}}))\\ 
  \partial_t\sum_{\sigma \in \sigma(X)}\lov{\prod_{\tau =1}^I}S_{\sum_{i} X_{i}^{j(i,\tau)}}&=& 
  (\sum_{i} X_{i}^{I_i})\sum_{\sigma \in 
  \sigma(X)}\lov{\prod_{\tau =1}^I}S_{\sum_{i} 
  X_{i}^{j(i,\tau)}}+\\ 
&& \sum_p\partial_t(\lov{\prod_{j<I_p} 
  S_{X_p^j}})\prod_{i\neq p}(\lov{\prod_j  S_{X_i^j}})). 
\end{eqnarray*} 
Or on remarque que: 
\begin{eqnarray*} 
\partial_t\prod_{i}(\lov{\prod_j  S_{X_i^j}})&=& 
\sum_p\partial_t(\lov{\prod_j  S_{X_p^j}})\prod_{i\neq p}(\lov{\prod_j  S_{X_i^j}}) 
\\ 
\partial_t\prod_{i}(\lov{\prod_j  S_{X_i^j}})&=&(\sum_{i} X_{i}^{I_i})\prod_{i}(\lov{\prod_j  S_{X_i^j}})+ 
 \sum_p(\lov{\prod_{j<I_p}  S_{X_p^j}})\prod_{i\neq p}(\lov{\prod_j  S_{X_i^j}}) .
 \end{eqnarray*} 
Ainsi l'{\'e}galit{\'e} est d{\'e}duite car au mois un des $ \Re(X_{i}^{I_i})+f(X_{i}^{I_i})$ est de m{\^e}me signe que $\Re(\sum_{i} X_{i}^{I_i})+f(\sum_{i} X_{i}^{I_i})$ 
\end{proof}  
 
Soit  $(X)$ d\'efinit par ses $i$ chaines \'el\'ementaires de taille  $I_i$,
alors
$G_{\un{X}}$ est l'ensemble des multis entiers  $\un{j}$ de $\N^I$ tel que: 
$\forall i, j_i \in \{0,1,.. I_i\}.$ 
On pose que  $\forall i, X_i^0=0$. 
 
\begin{definition}\label{Dchgen} 
 
Pour tout  $\un{j}$ de $G_{\un{X}}$; 
$\sum_i X_i^{j_i}$
est un g\'en\'erateur de $(X)$. L'ensemble des g\'en\'erateurs de  $(X)$ est
$\{\bigcup_{\un{j}\in G_{\un{X}}}\sum_i X_i^{j_i} \}.$ 
Alors il existe une bijection entre $G_X$ et
$\prod_i\{1,...I_i\}$. On peu voir que le nombre de g\'en\'erateurs de $(X)$ est
donn\'e par$\prod_i (I_i).$ 
\end{definition}

\begin{theoreme}\label{TCMgen}

Soit   $X$ une chaine arbitraire de taille  $I$, $h$ une fonction integrable: 
\begin{multline*}
\sum_{\sigma \in \sigma(X)}<\lov{\prod_{\tau =1}^I}S_{\sum_{i}  
  X_{i}^{j(i,\tau)}+\delta},h>
  =\\
(\sum_{\un{j}\in G_X}\prod_{i}(\lov{\prod_{j>j(i)}} 
  S_{Y(\un{j})_i^j})<S_{\un{j}}, 
  (\prod_{i}\rov{\prod_{j<j(i)}}(-1)^{\un{j}} 
  S_{Y_{i}^{j}(\un{j})})(h))>,\end{multline*} 
o\`u $Y_{i}^{j}(\un{j})= (-X_{i}^{j(i)}+X_{i}^{j})signe(j-j(i))$,
$(-1)^{\un{j}}=\prod_i(-1)^{j(i)}$ et
$S_{\un{j}}=S_{\delta+\sum_iX_i^{j(i)}}$.
\end{theoreme} 
\begin{proof}
\begin{eqnarray*} 
\partial_tS_{w_1+\delta}S_{\delta}&=&(w_1+\delta)S_{w_1+\delta}S_{\delta}+S_{\delta}\\ 
\partial_t((S_{w_1})S_{\delta}-S_{w_1+\delta}(S_{w_1})Id)&=& (w_1+\delta)((S_{w_1})S_{\delta}-S_{w_1+\delta}(S_{\delta})Id)\\ 
&&+S_{\delta}+(S_{w_1})Id-(S_{w_1})Id\\ 
\end{eqnarray*} 
Alors les fonctions $(S_{w_1+\delta}S_{\delta})$ et
$((S_{w_1})S_{\delta}-S_{w_1+\delta}(S_{X})Id)$ sont solutions de la m\^eme  ODE
at de m\^eme condition initiale: 
\[(S_{w_1+\delta}S_{\delta})(t_w)=0\] 
\[((S_{w_1})S_{\delta}-S_{w_1+\delta}(S_{X})Id)(t_w)=((S_{w_1})(S_{\delta})(t_w)=0\] 
car au moins   $w_1$ ou  $\delta$ est de m\^eme signe que $w_1+\delta$. 

On montre alors par r\'ecurence: 
\[\partial_t\sum_{\sigma \in \sigma(X)}(\lov{\prod_{\tau =1}^I}S_{\sum_{i} X_{i}^{j(i,\tau)}+\delta}) 
=(\sum_{i} X_{i}^{J_i}+\delta)\sum_{\sigma \in \sigma(X)}(\lov{\prod_{\tau =1}^I}S_{\sum_{i} X_{i}^{j(i,\tau)}+\delta})+\]\[ 
\sum_{\sigma \in \sigma(X)}(\lov{\prod_{\tau =1}^{I-1}}S_{\sum_{i} X_{i}^{j(i,\tau)}+\delta}) 
\] 
Or par r\'ecurence: 
 
\[\sum_{\sigma \in \sigma(X)}(\lov{\prod_{\tau =1}^{I-1}}S_{\sum_{i} X_{i}^{j(i,\tau)}+\delta})= 
\sum_i\sum_{S_{\un{X\setminus X_i^{J_i}}}}\prod_{i}(\lov{\prod_{j}}S_{Y(\un{j})_i^j}) 
(S_{\un{j}}\prod_{i}\rov{\prod_{j}}(-1)^{\un{j}}S_{Y_i^{j}(S_{\un{j}})})\] 
On en d\'eduit:

\[\partial_t\sum_{S_{\un{X}}}\prod_{i}(\lov{\prod_{j}}S_{Y(\un{j})_i^j}) 
(S_{\un{j}}\prod_{i}\rov{\prod_{j}}(-1)^{\un{j}}S_{Y_i^{j}(S_{\un{j}})})=\]\[(\sum_{i} X_{i}^{J_i}+\delta)\sum_{S_{\un{X}}}\prod_{i}(\lov{\prod_{j}}S_{Y(\un{j})_i^j}) 
(S_{\un{j}}\prod_{i}\rov{\prod_{j}}(-1)^{\un{j}}S_{Y_i^{j}(S_{\un{j}})}) 
+\]\[ \sum_i\sum_{S_{\un{X\setminus X_i^{J_i}}}}\prod_{i}(\lov{\prod_{j}}S_{Y(\un{j})_i^j}) 
(S_{\un{j}}\prod_{i}\lov{\prod_{j}}S_{Y_i^{j}(S_{\un{j}})})\sum_{S_{\un{X}}}\prod_{i}(\lov{\prod_{j}}S_{Y(\un{j})_i^j}) 
(\prod_{i}\rov{\prod_{j}}(-1)^{\un{j}}S_{Y_i^{j}(S_{\un{j}})}) 
\] 
  
De m\^eme la condition initiale est donn\'e par
  
\[\sum_{S_{\un{X}}}\prod_{i}(\lov{\prod_{j}} S_{Y(\un{j})_i^j}) 
 (S_{\un{j}}\prod_{i}\rov{\prod_{j}}(-1)^{\un{j}}S_{Y_i^{j}(S_{\un{j}})})(t_{(\sum_{i} X_{i}^{J_i})})=0.\] 
Nous devons maintenant montrer
\[\sum_{S_{\un{X}}}\prod_{i}(\lov{\prod_{j}}S_{Y(\un{j})_i^j}) 
(\prod_{i}\rov{\prod_{j}}(-1)^{\un{j}}S_{Y_i^{j}(S_{\un{j}})})  =0.\]
On le prouve par r\'ecurence, car pour tout $i$
\[\sum_{j(i)=0}^{J_i}(\lov{\prod_{j}}S_{Y(\un{j})_i^j})\rov{\prod_{j}}(-1)^{j(i)}S_{Y_i^{j}(S_{\un{j}})})= \sum_{j(i)=1}^{J_i}(\lov{\prod_{j}}S_{Y(\un{j})_i^j})\rov{\prod_{j}}(-1)^{j(i)}S_{Y_i^{j}(S_{\un{j}})})+(\prod_jS_{X_i^j})\]
Donc nous pouvons voir que:
\[(\prod_jS_{X_i^j})=(\prod_jS_{X_i^j-X_i^1+X_i^1}S_{X_i^1}),\]
si on note
 $\delta={X_i^1}$,
en applicant l'hypoth\`ese de r\'ecuence:
\[(\prod_jS_{X_i^j})=(\prod_jS_{X_i^j-X_i^1+X_i^1}S_{X_i^1})=\sum_{j(i)=0}^{J_i-1}(\lov{\prod_{j}}S_{\overline{Z}(S_{\un{j}})_i^j})\rov{\prod_{j}}(-1)^{j(i)}S_{Z_{i(k)}^{j(k)-j}(S_{\un{j}})})\]
avec $Z_{i(k)}^{j(k)-j}(S_{\un{j}})=X_i^j-X_i^1-(X_i^{j(i)}-X_i)$ et
$\overline{Z}(S_{\un{j}})_i^j=-X_i^j+X_i^1+(X_i^{j(i)}-X_i)$. On conclu que:
\[\sum_{j(i)=1}^{J_i}(\lov{\prod_{j}}S_{Y(\un{j})_i^j})\rov{\prod_{j}}(-1)^{j(i)}S_{Y_i^{j}(S_{\un{j}})})+(\prod_jS_{X_i^j})=0\]
et gr\`ace \{a la prorpi\'et\'e du produit:
\[\sum_{S_{\un{X}}}\prod_{i}(\lov{\prod_{j}}S_{Y(\un{j})_i^j}) 
(\prod_{i}\rov{\prod_{j}}(-1)^{\un{j}}S_{Y_i^{j}(S_{\un{j}})})=0.\]

\end{proof}
 \end{appendix} 
\small{
}

\end{document}